\documentclass[10pt]{article}
\usepackage[margin=1in]{geometry}  
\usepackage{graphicx} 
\usepackage{amssymb}             
\usepackage{amsmath}   
\usepackage{slashed}            
\usepackage{amsfonts}              
\usepackage{amsthm} 
\usepackage{breqn}
\usepackage{verbatim}              
\usepackage{enumitem}
\usepackage{fancyhdr}
\theoremstyle{definition}
\newtheorem{thm}{Theorem}[section]

\newtheorem{mydef}{Definition}[section]
\newtheorem{lem}[thm]{Lemma}
\newtheorem{prop}[thm]{Proposition}
\newtheorem{cor}[thm]{Corollary}

\newtheorem*{rmk}{Remark}

\newtheorem*{pf}{Proof}

\newcommand{\RR}{\mathbb{R}}      

\newcommand{\q}{\quad}
\newcommand{\p}{\partial}
\newcommand{\DD}{\mathcal{D}}
\newcommand{\curl}{\text{curl}\,}
\newcommand{\lei}{L^{2}(\mathcal{D}_{t})}
\newcommand{\leb}{L^{2}(\p\mathcal{D}_{t})}
\newcommand{\lli}{L^{2}(\Omega)}
\newcommand{\llb}{L^{2}(\p\Omega)}
\newcommand{\den}{{\bar\rho}_{0}}
\newcommand{\Dt}{{\p}_{t}+v^{k}{\p}_{k}}

\newcommand{\noi}{\noindent}
\newcommand{\nab}{\nabla}

\newcommand{\lap}{\Delta}
\newcommand{\di}{\text{div}\,}

\newcommand{\cnab}{\overline{\nab}}
\newcommand{\cp}{\overline{\partial}}
\newcommand{\dx}{\,dx}
\newcommand{\vol}{\text{vol}\,}
\newcommand{\lee}{\langle\langle}
\newcommand{\ree}{\rangle\rangle}
\newcommand{\symdot}{\tilde{\cdot}}
\newcommand{\linf}{L^{\infty}(\p\Omega)}
\newcommand{\kk}{\kappa}

\newcommand{\e}{\mathfrak{e}}
\newcommand{\lliw}{L^2_w(\Omega)}
\newcommand{\llbw}{L^2_w(\p\Omega)}
\numberwithin{equation}{section}

\begin{document}
\title{On the Motion of a Compressible Gravity Water Wave with Vorticity}
\author{Chenyun Luo}
\date{}
\maketitle
\begin{abstract}
 We prove a priori estimates for the compressible Euler equations modeling the motion of a liquid with moving physical vacuum boundary with unbounded initial domain. The liquid is under influence of gravity but without surface tension. Our fluid is not assumed to be irrotational. But the physical sign condition needs to be assumed on the free boundary. We generalize the method used in \cite{LL} to prove the energy estimates in an unbounded domain up to arbitrary order. In addition to that, the a priori energy estimates are in fact uniform in the sound speed $\kk$. As a consequence, we obtain the convergence of solutions of compressible Euler equations with a free boundary to solutions of the incompressible equations, generalizing the result of \cite{LL} to when you have an unbounded domain. On the other hand, we prove that there are initial data satisfying the compatibility condition in some weighted Sobolev spaces, and this will propagate within a short time interval, which is essential for proving long time existence for slightly compressible irrotational water waves.
\end{abstract}
\tableofcontents
\section{Introduction}
We consider the compressible water wave problem in $\mathbb{R}^n, n=2,3$. We use the notation $\DD_t:=\{(x',x_n)\in\RR^n: x_n\leq \Sigma(t,x')\}$ to represent the domain occupied by the fluid at each fixed time $t$, whose boundary $\p\DD_t=\{(x',x_n):x_n=\Sigma(t,x')\}$ moves with the velocity of the fluid. Under this setting, the motion of the fluid is described by the Euler equations
\begin{align}
\begin{cases}
D_tv := \p_t v+\nab_v v = -\frac{1}{\rho}\p p - g\e_n,\q x\in\DD \label{e}\\
D_t\rho+\rho\di v=0\q x\in\DD,
\end{cases}
\end{align}
with the initial and boundary condition
\begin{align}
\begin{cases}
\{x:(0,x)\in \DD\} =\DD_0\\
v=v_0,\rho=\rho_0 \q\text{on}\,\{0\}\times\DD_0
\end{cases}
\begin{cases} \label{ibcond}
D_t|_{\p\DD}\in T(\p\DD)\\
p|_{\p\DD}=0
\end{cases}
\end{align}
where   
$\DD : = \cup_{0\leq t\leq T}\{t\}\times \DD_t$,
$g>0$ is the gravity constant and $\rho$ denotes the density of the fluid and the equation of the state is given by 
$$
p = p(\rho),\q p'(\rho)>0,\q\text{for}\q \rho\geq \den ,
$$
where $\den:=\rho|_{\p\DD}>0$ is a constant (for simplicity, we set $\den=1$), which is in the case of a liquid.

We prove the energy estimates for the local (in time) solutions of system (\ref{e})-(\ref{ibcond}), taking prescribed initial data, such that for every fixed time $t\in[0,T]$, $|v(t,x)|\to 0$, $|v_t(t,x)|\to 0$, and $\Sigma(t,x') \to \{(x',0):x'\in\RR^{n-1}\}$ as $|x|\to \infty$.  In fact, we are able to show that there exist initial data satisfying the compatibility condition \eqref{cpt cond} in some weighted Sobolev spaces with weight $w(x) = (1+|x|^2)^\mu, \mu\geq 2$. This implies that our data is at least of $O(|x|^{-2})$ as $|x|\to\infty$.

 We introduce the enthalpy $h$ to be a function of the density, i.e., $h(\rho)=\int_{1}^{\rho}p'(\lambda)\lambda^{-1}\,d\lambda$. Since $\rho \geq \den=1$ can then be thought as a function of $h$, we define $e(h)=\log \rho(h)$. Under these new variables, (\ref{e})-(\ref{ibcond}) can be re-expressed as
\begin{align}
\begin{cases}
D_t v = -\p h -  g\e_n,\q \text{in}\, \mathcal{D}\\
\di v = -D_t e(h)=-e'(h)D_t h.\q \text{in}\, \mathcal{D} \label{EE}
\end{cases}
\end{align}
Together with initial and boundary conditions
\begin{align}
\begin{cases}
\{x:(0,x)\in\DD\} = \DD_0,\\
v=v_0,h=h_0 \q \text{on}\,\{0\}\times\DD_0. \label{EIBC}
\end{cases}
\begin{cases}
D_t|_{\p \DD}\in T(\p\DD),\\
h = 0 \q \text{on}\,\p\DD.
\end{cases}
\end{align}
(\ref{EE}) looks exactly like the incompressible Euler equations, where $h$ takes the position of $p$ and $\di v$ is no longer $0$ but determined by $h$. In addition, we take the gravity constant $g=1$. On the other hand, we would like to impose the following natural conditions on $e(h)$: For each fixed $r\geq 1$, there exists a constant $c_0$ such that
\begin{align}
|e^{(k)}(h)|\leq c_0,\q \text{and}\,\, |e^{(k)}(h)| \leq c_0|e'(h)|^k \leq c_0|e'(h)|, \q \text{if}\,\, k\leq r+1.
\label{e_kk}
\end{align}
\indent In order for the initial boundary problem (\ref{EE})-(\ref{EIBC}) to be solvable the initial data has to satisfy certain compatibility conditions at the boundary. By the second equation in (\ref{e}),(\ref{ibcond}) implies that $\di v|_{\p\DD}=0$. We must therefore have $h_0|_{\p\DD_{0}}=0$ and $\di v_0|_{\p\DD_0}=0$, which is the zero-th compatibility condition. Furthermore, $m$-th order compatibility condition can be expressed as
\begin{align}
(\Dt)^j h|_{\{0\}\times\p\DD_{0}}=0\q j = 0,\cdots,m.
\label{cpt cond}
\end{align}
In \cite{LL}, we have proved that for each fixed $m$,  there exists initial data satisfying $m$-th order compatibility condition if the sound speed $c(t,x)$ is sufficiently large. In addition, the energies $E_r$, defined as (\ref{Er}), are bounded uniformly at time $0$, regardless of the sound speed.

\indent Let $N$ be the exterior unit normal to the free surface $\p\DD_t$. We will prove a priori bounds for (\ref{EE})-(\ref{EIBC}) in Sobolev spaces under the assumption
\begin{align}
\nab_{N}h\leq -\epsilon<0\q \text{on}\,\,\p\DD_t,\label{taylor-sign cond}
\end{align}
where $\nab_{N}=N^i\p_i$ and $\epsilon>0$ is a constant. (\ref{taylor-sign cond}) is a natural physical condition. It says that the pressure and hence the density is larger in the interior than at the boundary. The system (\ref{e})-(\ref{ibcond}) is ill-posed in absence of (\ref{taylor-sign cond}), an easy counter-example can be found in \cite{CL} and \cite{Eb2}. Furthermore, if the fluid is assumed to be incompressible ($\di v=0$) and irrotatonal ($\curl v=0$), (\ref{taylor-sign cond}) can in fact be proved via strong maximum principle. In addition, under the presence of the gravity, it can be shown that $-\nab_N p$ can actually be bounded uniformly by some positive constants from below (see Wu \cite{W1,W2}). In fact, we show that the assumption is plausible by proving \eqref{taylor-sign cond} in the case when the liquid is slightly compressible and irrotational (Section \ref{section 7}). Heuristically, in the Lagrangian coordinates (where $x:=x(t,y)$, $\frac{dx(t,y)}{dt}=v(t,x(t,y))$, see Section \ref{section 2}), we have $x_{tt} = v_t$ and so 
$$
-\nab_N h = x_{tt}\cdot N + N \cdot e_n,
$$
and because $v_t = x_{tt}$ decays to $0$ at infinity, we conclude $-\nab_N h\geq \epsilon>0$ for some $\epsilon>0$ pointwisely. We shall discuss more about this in Section \ref{section 7} (see the remark after Theorem \ref{thm 7.2}).  But (\ref{taylor-sign cond}) needs to be assumed if the fluid is rotational and without surface tension.

\indent Euler equations involving free-boundary has been studied intensively by many authors. The first break through in solving the well-posedness for the incompressible and irrotational water wave problem for general data came in the work of Wu \cite{W1,W2} who solved the problem in both two and three dimensions. For the general incompressible problem with nonvanishing curl Christodoulou and Lindblad \cite{CL} were the first to obtain the energy estimates assuming the physical sign condition. In addition, Zhang and Zhang \cite{ZZ} generalized Wu's work to incompressible water wave with nonvanishing curl. For the compressible problem, Lindblad \cite{L} later proved local well-posedness for the general problem modeling the motion of a liquid in a bounded domain by Nash-Moser iteration, and this result was generalized to the case of an unbounded domain by Trakhinin \cite{Tr}. But these results do not contain a priori estimates for the solution due to the loss of regularity on the moving boundary.

Very recently, together with Lindblad, we obtained a new type of a priori energy estimates for the compressible Euler equations with free boundary in a bounded domain, which are uniform in the sound speed \cite{LL}. This, in fact, leads to the convergence of solutions of compressible Euler equations with a free boundary to solutions of the incompressible equations in a bounded domain. In other words, we proved the so-called \textit{incompressibe limit problem} for compressible free boundary Euler equations.  It is worth mentioning here that the incompressible limit in $\RR^n$ or $\mathbb{T}^n$ was established in \cite{Eb0,  KM1, KM2, Majda, MS}, and there are further works (e.g., \cite{Alazard, Cheng1, Cheng2, DE, Eb, S}) that treat the incompressible limit in domains with fixed boundary.  

The goal of this paper is to generalize the above results to compressible water waves, i.e., the fluid domain becomes unbounded and diffeomorphic to the half space. To our knowledge, these results appear to be the first that concern a priori energy bounds for a compressible water wave. Furthermore, the incompressible limit allows one to approximate a slightly compressible water wave by an incompressible water wave, for which the long time existence is well-known, e.g., \cite{GMS, HIT,IT,IP,W3,W4}.  In addition, we show that the a priori energy estimates can also be generalized to weighted $L^2$-Sobolev spaces, which is an essential first step for proving long time existence also for compressible water waves.

\subsection{Energy conservation and higher order energies}
\indent The boundary conditions $p|_{\p\DD_t}=0$ and $\rho|_{\p\DD_t}=1$ leads to that the zero-th order energy is conserved, i.e., let
\begin{align}
E_0(t)= \frac{1}{2}\int_{\DD_t}\rho|v|^2\,dx +\int_{\DD_t}\rho Q(\rho)\,dx+\int_{\DD_t\cap\{x_n>0\}}x_n\dx-\int_{\DD_t^c\cap \{x_n<0\}}x_n\dx +\int_{\DD_t}(\rho-1)x_n\dx\label{E_0}
\end{align}
where $Q(\rho)=\int_{1}^{\rho}p(\lambda)\lambda^{-2}\,d\lambda $. These integrals are bounded here because of the decay properties of our functions involved.\\

A direct computation yields 
\begin{multline*}
\frac{d}{dt}E_0(t)
=-\int_{\DD_t}(\p_i p) v^i\dx - \int_{\DD_t}\rho(\p_ix_n)v^i\dx+\int_{\DD_t}p(\rho)D_t\rho \rho^{-1}\dx + \int_{\p\DD_t}s\cdot v_N\,dS + \int_{\DD_t}(\p_t\rho)x_n \dx\\
= (\int_{\DD_t} p (\di v)\dx + \int_{\DD_t}p D_t\rho \rho^{-1}\dx) + (\int_{\DD_t}(v\cdot\p\rho+\rho\di v)x_n\dx+\int_{\DD_t}(\p_t\rho)x_n \dx) = 0.
\end{multline*}

The higher order energies $E_r(t)$ are defined in a similar fashion, but instead of using the regular inner product, we introduce a positive definite quadratic form $Q$ which, when restricted to the boundary, is the inner product of the tangential components, i.e., $Q(\alpha,\beta)=\Pi \alpha\cdot\Pi\beta$, where $\alpha$ and $\beta$ are $(0,r)$ tensors. To be more specific, we define
\begin{align}
Q(\alpha,\beta)= q^{i_1j_1}\cdots q^{i_rj_r}\alpha_{i_1\cdots i_r}\beta_{j_1\cdots j_r},
\end{align}
where 
\begin{align*}
q^{ij}=\delta^{ij}-\eta(d)^2\mathcal{N}^i\mathcal{N}^j,\\
d(x) = dist(x,\p\DD_t),\\
\mathcal{N}^i=-\delta^{ij}\p_j d.
\end{align*}
Here $\eta$ is a smooth cut-off function satisfying $0\leq \eta(d)\leq1$, $\eta(d)=1$ when $d\leq\frac{d_0}{4}$ and $\eta(d)=0$ when $d>\frac{d_0}{2}$. $d_0$ is a fixed number that is smaller than the injective radius $l_0$, which is defined to be the largest number $l_0$ such that the map
\begin{align}
\p\DD_{t}\times(-l_0,l_0)\to\{x:dist(x,\p\DD_t)<l_0\}, \label{inj rad}
\end{align}
given by
\begin{align}
(\bar{x},l)\to x=\bar{x}+l\mathcal{N}(\bar{x}),
\end{align}
is an injection.\\
\indent The  higher order energies we propose are
\begin{equation}
E_r = \sum_{s+k=r} E_{s,k}+K_r+\sum_{j\leq r+1}W_{j}^2\label{Er},\q r\geq 2,\q E_r^*=\sum_{r'\leq r} E_{r'},
\end{equation}
 where
\begin{dmath}
E_{s,k}(t)=\frac{1}{2}\int_{\DD_t}\rho \delta^{ij}Q(\p^sD_t^kv_i,\p^sD_t^{k}v_j)\dx+\frac{1}{2}\int_{\DD_t}\rho e'(h)Q(\p^sD_t^k h,\p^s D_t^kh)\dx\\
+\frac{1}{2}\int_{\p\DD_t}\rho Q(\p^sD_t^k h,\p^sD_t^k h)\nu\,dS,  \label{Esk}
\end{dmath}
where $\nu = (-\nab_{N}h)^{-1}$ and
\begin{align}
K_r(t) &= \int_{\DD_t}\rho|\p^{r-1}\curl v|^2\dx \label{K_r},\\
W_r(t) &= \frac{1}{2}||\sqrt{e'(h)}D_t^rh||_{\lei}+\frac{1}{2}||\nab D_t^{r-1}h||_{\lei} .\label{Er time}
\end{align}

Here $W_r$ is the (higher order) energy for the wave equation
\begin{align}
D_t^2e(h) -\lap h=(\p_i v^j)(\p_j v^i),\label{waveequation}
\end{align}
which is obtained by commuting divergence through the first equation of (\ref{e}) using
\begin{equation}
[D_t, \p_i] = -(\p_iv^j)\p_j.\label{lowestcommutator}
\end{equation}

The energies $E_r$ defined above in fact control all components of
\begin{align*}
||v||_{r,0} := \sum_{k+s=r, k<r} ||\p^s D_t^k v||_{\lei},\\
||h||_r := \sum_{k+s=r, k<r}||\p^s D_t^k h||_{\lei} + ||\sqrt{e'(h)}D_t^rh||_{\lei},\\
||D_th||_{r,1} = \sum_{k+s=r,k<r-1}||\p^sD_t^{k+1} h||_{\lei} +||\sqrt{e'(h)}\nab D_t^rh||_{\lei}+||e'(h)D_t^{r+1}h||_{\lei},\\
\lee h\ree_r :=\sum_{k+s=r} ||\p^s D_t^kh||_{\leb},
\end{align*}
in the interior and on the boundary (section \ref{section 5}).  Although $E_r$ only controls the tangential components, the fact that
we also control the divergence $\sum_{j\leq r+1}W_{j}^2$ (through $\di v=-D_t e(h)$) and the curl $K_r$ allows us to control all components. In fact, by a Hodge type decomposition
\begin{align}
|\p v| \lesssim |\overline{\p}v|+|\di v|+|\curl v| \label{intro hodge},
\end{align}
where the tangential derivatives are given by $\overline{\p}h=\Pi\p h$. In addition, if $|\nabla_{N} h|\geq \epsilon>0$ then the boundary term gives an estimate
for the regularity of the boundary. In fact, one can show that if $q$ vanishes on the boundary then
\begin{equation}
\Pi \p^r q = (\overline{\p}^{r-2}\theta)\nab_{\mathcal{N}} q+O(\p^{r-1}q)+O(\overline{\p}^{r-3}\theta),\label{projest0}
\end{equation}
where $\theta$  is the second fundamental form of the boundary and $\overline{\p}$ stand for tangential derivatives, so
\begin{equation}
\| \overline{\p}^{r-2}\theta\|_{\leb}^2 \leq \frac{C}{\epsilon} E_r^*
+ C \sum_{r'\leq r-1} \| \p^{r'} h\|_{\leb}^2.\label{projest}
\end{equation}
Now, because of the estimates \eqref{intro hodge}-\eqref{projest}, using elliptic estimates (section 3) one can show that
\begin{align}
||v||_{r,0}^2+||h||_r^2 \leq C_r(K,M,c_0, E_{r-1}^*)E_r^* \label{intro int  est},\\
||D_th||_{r,1}^2+\lee h\ree_r^2 \leq C_r(K,M,c_0,\frac{1}{\epsilon},  E_{r-1}^*)E_r^* \label{intro bdy est},
\end{align}
for some continuous function $C_r$. In fact, we use many of such functions throughout this paper, but we shall not distinguish them unless otherwise specified, i.e., $C_r$ would always denote continuous functions depend on constants $K,M,c_0,\frac{1}{\epsilon}$ and the energies $E_{r-1}^*$.

\subsection{The main results}
Sections \ref{section 4}-\ref{section 5} are devoted to prove a priori energy estimates implying that the energies $E_r$ remain bounded as long as certain a priori assumptions are true. To be specific, we show
\thm \label{perp main} Let $v,h$ be the solutions for \eqref{EE}-\eqref{EIBC} and $E_r $ be the energy defined as (\ref{Er}), then for each fixed integer $r\geq 1$
\begin{align}
|\frac{dE_r(t)}{dt}|\leq C_r(K,\frac{1}{\epsilon},M, c_0, E_{r-1}^*)E_r^*(t) \label{intro main1}
\end{align}
holds, where $C_r$ are continuous functions and $E_r^*=\sum_{s=0}^rE_s$, provided \eqref{e_kk} a priori assumptions
\begin{align}
|\theta|+\frac{1}{l_0}\leq K,\q\text{on}\,\p\DD_t\label{geometry_bound}\\ 
-\nab_{N}h\geq \epsilon>0,\q\text{on}\,\p\DD_t\label{RT sign}\\
1\leq |\rho| \leq M,\q\text{in}\,\DD_t \label{bootstrap rho}\\
|\p^j \curl_{ij} v| \leq M,\q\text{in}\,\DD_t \label{bootstrap curl}\\
|\p v|+|\p h|+|\p^2 h| + |\p D_th|\leq M,\q\text{in}\,\DD_t \label{v,h}\\
|e'(h)D_th|+|e'(h)D_t^2h|\leq M.\q\text{in}\,\DD_t \label{D_te(h)}
\end{align}
The bounds (\ref{geometry_bound}) gives us control of geometry of the free surface $\p\DD_t$. A bound for the second fundamental form $\theta$ gives a bound for the curvature of $\p\DD_t$, and a lower bound for the injectivity radius of the exponential map $l_0$ measures how far off the surface is from self-intersecting.  In the case when $\DD_t$ is unbounded, the uniform a priori bounds for $|D_th|$ and $ |D_t^2h|$ are weakened to the bounds \eqref{D_te(h)}, and we need them to hold uniformly to pass to the incompressible limit.

\rmk $r$ will be used to denote an integer throughout this manuscript. In particular, we will not use fractional Sobolev norms. 

\rmk The assumption \eqref{D_te(h)} can, in fact, be relaxed to $|\sqrt{e'(h)}D_th|+|\sqrt{e'(h)}D_t^2h|\leq M$. However, we will no longer be able to pass \eqref{intro main1} to the incompressible limit.

In Section 6, We show that the energy bounds \eqref{intro main1} remain valid uniformly as the sound speed goes to infinity. For physical reasons, the sound speed is defined by
$
c(t,x) = \sqrt{p'(\rho)}.
$
In this paper, the sound speed $\kk$ is defined by viewing $\{p_{\kk}(\rho)\}$ as a family parametrized by $\kk\in\RR^+$, such that for each $\kk$ we have
$
\kk:=p_{\kk}'(\rho)|_{\rho=1}.
$  
We consider the compressible Euler equations depend on $\kk$:

\begin{align}
\begin{cases}
D_t v_{\kk} = -\p h_{\kk} - \e_n,\\
\di v_{\kk} = -D_te_{\kk}(h) \label{intro E kk}.
\end{cases}
\end{align}
Here, we further assume that $e_{\kk}(h)$ satisfies:
\begin{align}
e_{\kk}(h)\to 0,\q\text{as}\,\,\kk\to\infty,
\label{e_kk a}
\end{align}
and for each fixed $r\geq 1$, there exists a constant $c_0$ such that
\begin{align}
|e_{\kk}^{(k)}(h)|\leq c_0,\q \text{and}\,\, |e_{\kk}^{(k)}(h)| \leq c_0|e_{\kk}'(h)|^k \leq c_0|e_{\kk}'(h)|, \q \text{if}\,\, k\leq r+1.
\label{e_kk b}
\end{align}
Under this setting, we show
\thm \label{intro uniform kk}
Let $(v_{\kk}, h_{\kk})$ solves \eqref{intro E kk}.  Let $\widetilde{E}_{r}$ be defined as $\widetilde{E}_{r}=\sum_{s+k=r}E_{s,k}+K_{r}+\sum_{j\leq r+1}\widetilde{W}_j$, where
$$
\widetilde{W}_j = \frac{1}{2}||e_{\kk}'(h)D_t^jh_{\kk}||_{\lli}+\frac{1}{2}||\sqrt{e_{\kk}'(h)}\nab D_t^{j-1}h_{\kk}||_{\lli}.
$$
If, in addition, the physical sign condition holds, i.e., 
$$
-\nab_N h_{\kk}\geq \epsilon>0,
$$
then there exists $T>0$, independent of $\kk$, such that for any smooth solutions of \eqref{intro E kk} for $0\leq t\leq T$ satisfies
\begin{align}
\widetilde{E}_{r,\kk}^*(t)\leq 2\widetilde{E}_{r,\kk}^*(0),\q \text{whenever}\,\,r> n/2+3/2
\end{align}
and this estimate can be carried over to the case when $\kk = \infty$, i.e., the energy estimates for the incompressible Euler equations.

Theorem \ref{intro uniform kk} is a direct consequence of the a priori energy bounds \eqref{intro main1} are uniform in $\kk$ via Gronwall's lemma. Moreover, \eqref{intro main1} remains valid since that our estimates do not depend on the lower bound of $e_{\kk}^{(k)}(h)$, which goes to $0$ as $\kk\to\infty$, and the elliptic estimates \eqref{intro int  est}-\eqref{intro bdy est} can be carried to the incompressible case apart from the term $||\p D_t^{k} h||_{\lei}, 0\leq k\leq r-1$. But this can be bounded via $||\lap D_t^kh||_{\lli}$, given that $D_t^kh$ decays sufficiently fast at infinity (see Section \ref{section 6}.1). 

In addition, apart from the coefficient in front of the highest order time derivative our energy does not depend in crucial way on $\kk$ but uniformly (as $\kk\to \infty$) control the corresponding norms of all but the highest order time derivative. This leads to that the a priori
$L^\infty$ bounds also hold uniformly and the norms are bounded uniformly up to a fixed time. The convergence of solutions for the compressible Euler equations to the solution for the incompressible equations then follows from Arzela-Ascoli theorem. 

\thm \label{intro convergence} Let $u_0$ be a divergence free vector field such that its corresponding pressure $p_0$, defined by $\lap p_0 = -(\p_i u_0^k)(\p_k u_0^i)$ and $p_0\big|_{\p\DD_0}=0$, satisfies the physical condition $-\nab_N p_0\big|_{\p\DD_0} \geq \epsilon >0$. Let $(u, p)$ be the solution of the incompressible free boundary Euler equations with data $u_0$, i.e.
$$
\rho_0 D_t u = -\p p,\q \di u=0,\qquad p|_{\p\DD_0}=0,\quad u|_{t=0}=u_0
$$
with the constant density $\rho_0 =1$. Furthermore, let $(v_{\kk}, h_{\kk})$ be the solution for the compressible Euler equations \eqref{intro E kk}, with the density function $\rho_{\kk} :h\rightarrow \rho_{\kk}(h)$, and the initial data $v_{0\kk}$ and $h_{\kk}|_{t=0} = h_{0\kk}$, satisfying the compatibility condition \eqref{cpt cond} up to order $r+1$, as well as the physical sign condition \eqref{taylor-sign cond}. Suppose that $\rho_{\kk} \to \rho_0=1$, $v_{0\kk}\to u_0$ and $h_{0\kk}\to p_0$ as $\kk\to\infty$, such that $E_{r,\kk}^*(0)$ is bounded uniformly independent of $\kk$, then
$$
(v_{\kk},h_{\kk})\to (u,p).
$$

\rmk It is essential to make sure that the data satisfies \eqref{cpt cond} in Theorem \ref{intro convergence}. A good example is that if one starts with $v_0=u_0$, where $\di u_0=0$ and $h_0=0$ (e.g., $\rho_0=1$), then it is easy to see, after taking divergence on both sides of first equation of \eqref{intro E kk}, $D_t^2 h|_{t=0}=(\p u_0)\cdot(\p u_0)$, and this would in general contradicts that $D_t^2 h=0$ on $\p\DD_t$. To prevent this from happening, we give data in terms of enthalpy $h$, and hence $\rho \sim 1+h/\kk$. We are then able to construct initial data via solving a system of elliptic equations.

In Section \ref{section 7} we prove that there exist initial data satisfying the conditions in Theorem \ref{intro convergence} in weighted Sobolev spaces $H^{r+1}_w(\Omega)$ with $w(x)=(1+|x|^2)^{\mu},\mu\geq 2$ when $\kk$ is sufficiently large.  In particular, we prove:
\thm \label{convergence of data} Let $u_0$ and $p_0$ are the initial data for the incompressible Euler equations defined in Theorem \ref{intro convergence}, and we further assume $u_0\in H_w^s(\DD_0)$ for $s\geq r+1$, where $r> \frac{n}{2}+\frac{3}{2}$, if $\DD_0$ is unbounded, diffeomprphic to the half space. Let $\rho_\kappa(h)\sim\rho_0+h/\kappa$, then there exists initial data $v_{0\kk}$ and $h_{0\kk}$ satisfying the compatibility condition \eqref{cpt cond} up to order $r+1$, such that $v_{0\kk}\to u_0$,  $h_{0\kk}\to p_0$ as $\kk\to \infty$, and $E_{r,\kk}^*(0)$ (and hence $\widetilde{E}_{r,\kk}^*(0)$) is uniformly bounded for all $\kk$.

\rmk Theorem \ref{convergence of data} implies that we need the initial data to be in $H^5_w$ when $n=3$ and in $H^4_w$ when $n=2$. 

In addition to this, we show that the physical sign condition \eqref{taylor-sign cond} can be verified via the maximum principle when the liquid is assumed to be irrotational.  Finally, Section 8 is devoted to prove the weighted energy estimates for the compressible Euler equations, as an analogue to Theorem \ref{perp main}. 

\thm \label{intro weighted estimate}
Let $E_{w,r}$ be defined as \eqref{Er} in $L^2_w(\DD_t)$ with weight $w(x)=(1+|x|^2)^{\mu},\mu\geq 2$. Then for each fixed $r\geq 1$, we have
\begin{align}
|\frac{dE_{w,r}(t)}{dt}|\leq C_r(K,\frac{1}{\epsilon},M, c_0, E_{w, r-1}^*)E_{w,r}^*(t),
\end{align}
provided the \eqref{e_kk} and a priori assumptions \eqref{geometry_bound}-\eqref{D_te(h)}. 

Theorem \ref{intro weighted estimate} suggests that one should also be able to prove the (local) existence for localized solutions, given localized initial data constructed in Theorem \ref{convergence of data}, which serves as a good preparation for proving dispersive estimates and hence long time existence for a compressible water wave.

\section{The Lagrangian coordinates}\label{section 2}
\indent Let us first introduce Lagrangian coordinate, under which the boundary becomes fixed. Let $\Omega=\{x\in\RR^n:x_n\leq 0\}$ be the half space in $\RR^n$, and let $f_0:\Omega\to\DD_0$ to be a diffeomorphism. The Lagrangian coordinate $(t,y)$ where $x=x(t,y)=f_t(y)$ are given by solving
\begin{align}
\frac{dx}{dt}=v(t,x(t,y)), \q x(0,y)=f_0(y),\q y\in\Omega \label{change}
\end{align}
The boundary becomes fixed in the new coordinate, and we introduce the notation
\begin{align}
D_t = \frac{\p}{\p t}|_{y=\text{constant}} = \frac{\p}{\p t}|_{x=\text{constant}} + v^k\frac{\p}{\p x^k}, \label{lag}
\end{align}
to be the material derivative and
\begin{align*}
\p_i = \frac{\p}{\p x^i} = \frac{\p y^a}{\p x^i}\frac{\p}{\p y^a}.
\end{align*}
Due to (\ref{lag}), we shall also call $D_t$ as the time derivative as well by slightly abuse of terminology.

\indent Sometimes it is convenient to work in the Eulerian coordinate $(t,x)$, and sometimes it is easier to work in the Lagrangian coordinate $(t,y)$. In the Lagrangian coordinate the partial derivative $\p_t=D_t$ has more direct significance than it in the Eulerian frame. However, this is not true for spatial derivatives $\p_i$. The notion of space derivative that plays a more significant role in the Lagrangian coordinate is that the covariant differentiation with respect to the metric $g_{ab}(t,y)$. We shall not involve covariant derivatives in our energy; instead, we use the regular Eulerian spatial derivatives. We will work mostly in the Lagrangian coordinate in this paper. However, our statements are coordinate independent.\\
\indent The Euclidean metric $\delta_{ij}$ in $\DD_{t}$ induces a metric
\begin{align}
g_{ab}(t,y)=\delta_{ij}\frac{\p x^i}{\p y^a}\frac{\p x^j}{\p y^b}, \label{g}
\end{align} 
in $\Omega$ for each fixed $t$. We will denote covariant differentiation in the $y_{a}$-coordinate by $\nab_a$, $a=1,\cdots,n$, and the differentiation in the $x_i$-coordinate by $\p_i$, $i=1,\cdots,n$. Here, we use the convention that differentiation with respect to Eulerian coordinates is denoted by letters $i,j,k,l$ and with respect to Lagrangian coordinate is denoted by $a,b,c,d$.  \\

\indent The regularity of the boundary is measured by the regularity of the normal, let $N^a$ to be the unit normal to $\p\Omega$, 
$$
g_{ab}N^aN^b=1,
$$
and let $N_a=g_{ab}N^b$ denote the unit co-normal, $g^{ab}N_aN_b=1$. The induced metric $\gamma$ on the tangent space to the boundary $T(\p\Omega)$ extended to be $0$ on the orthogonal complement in $T(\Omega)$ is given by
$$
\gamma_{ab}=g_{ab}-N_aN_b,\q \gamma^{ab}=g^{ac}g^{bd}\gamma_{cd}=g^{ab}-N^aN^b.
$$
The orthogonal projection of an $(0,r)$ tensor $S$ to the boundary is given by
$$
(\Pi S)_{a_1,\cdots,a_r}=\gamma_{a_1}^{b_1}\cdots\gamma_{a_r}^{b_r}S_{b_1,\cdots,b_r},
$$
where $\gamma_{a}^{b}=g^{bc}\gamma_{ac}=\delta_{a}^{b}-N_aN^b$. In particular, the covariant differentiation on the boundary $\overline{\nab}$ is given by
$$
\overline{\nab}S=\Pi \nab S.
$$
Note that $\overline{\nab}$ is invariantly defined since the projection and $\nab$ are. The second fundamental form of the boundary $\theta$ is given by $\theta_{ab}=(\cnab N)_{ab}$, and the mean curvature of the boundary $\sigma=tr\theta=g^{ab}\theta_{ab}$. \\  
\indent It is now important to compute time derivative of the metric $D_tg$, as well as the normal $D_tN$, as well as the time derivative of corresponding measures. 

\thm
Let $x=f_t(y)=x(t,y)$ be the change of variable given by (\ref{change}) and $g$ be the metric given by (\ref{g}),and $\gamma_{ab} = g_{ab}-N_a N_b$, where $N_a= g_{ab}N^b$ is the co-normal to $\p\Omega$, set
\begin{align}
v_a(t,y) = v_i(t,x)\frac{\p x^i}{\p y^a},\\
 u^a = g^{ab}u_b,\\
d\mu_g, \q \text{volume element with respect to the metric}\,\, g ,\\
d\mu_\gamma,\q\text{surface element with respect to the metric}\,\, \gamma.\label{dmu}
\end{align}
Then
\begin{align}
D_t g_{ab} = \nab_a v_b+\nab_b v_a,\label{Dtg}\\
D_t g^{ab} = -g^{ac}g^{bd}D_tg_{cd},\label{Dtg_inverse}\\
D_t N_a = -\frac{1}{2}N_a(D_tg^{cd})N_cN_d,\label{DtN}\\
D_t d\mu_g = \di v\,d\mu_g,\label{dg}\\
D_t d\mu_\gamma =(\sigma v \cdot N)\,d\mu_\gamma. \label{T2'}
\end{align}
\begin{proof}
The detail proof can be found in \cite{LL}.
\end{proof}

\section{Basic estimates on the domain with free boundary}
Most of the results in this section will be stated in a coordinate-independent fashion. Throughout this section, $\nab$ will refer to covariant derivative with respect to the metric $g_{ij}$ in $\Omega$, and $\cnab$ will refer to covariant differentiation on $\p\Omega$ with respect to the induced metric $\gamma_{ij}=g_{ij}-N_iN_j$. Hence, in this section, $\Omega$ will be used to denote a general domain with smooth boundary. In addition, we shall assume the normal $N$ to $\p\Omega$ is extended to a vector field in the interior of $\Omega$ satisfying $g_{ij}N^iN^j\leq 1$ by the same way introduced in lemma 2.1. 
\subsection{Elliptic estimates}
\mydef\label{tensor} Let $u:\Omega\subset\RR^n\to\RR^n$ be a smooth vector field, and $\beta_k=\beta_{Ik}=\nab_{I}^{r}u_k$ be the $(0,r)$-tensor defined based on $u_k$, where $\nab_{I}^{r}=\nab_{i_1}\cdots\nab_{i_r}$ and $I=(i_1,\cdots,i_r)$ is the set of indices. Let $\di \beta_k=\nab_i\beta^{i}=\nab^r\di u$ and $\curl \beta=\nab_i\beta_j-\nab_j\beta_i=\nab^r\curl u_{ij}$.
\mydef (Norms) If $|I|=|J|=r$, let $g^{IJ}=g^{i_1j_1}\cdots g^{i_r j_r}$ and $\gamma^{IJ}=\gamma^{i_1j_1}\cdots\gamma^{i_rj_r}$. If $\alpha$, $\beta$ are $(0,r)$ tensors, let $\langle\alpha,\beta\rangle=g^{IJ}\alpha_I \beta_J$ and $|\alpha|^2=\langle\alpha,\alpha\rangle$. If $(\Pi\beta)_{I}=\gamma_{I}^{J}\beta_J$ is the projection, then $\langle \Pi \alpha,\Pi\beta\rangle=\gamma^{IJ}\alpha_I\beta_J$.Let
\begin{align*}
||\beta||_{L^2(\Omega)}=(\int_{\Omega}|\beta|^2\,d\mu_g)^{\frac{1}{2}},\\
||\beta||_{L^2(\p\Omega)}=(\int_{\p\Omega}|\beta|^2\,d\mu_\gamma)^{\frac{1}{2}},\\
||\Pi\beta||_{L^2(\p\Omega)}=(\int_{\p\Omega}|\Pi\beta|^2\,d\mu_\gamma)^{\frac{1}{2}}.
\end{align*}
\indent We now state the following Hodge-type decomposition theorem, which serves as a main ingredient for proving the elliptic estimates. 
\thm \label{hodge}(Hodge-decomposition) Let $\beta$ be defined in Definition \ref{tensor}. If $|\theta|+|\frac{1}{l_0}|\leq K$, where $\theta$ is the second fundamental form and $l_0$ is the injective radius defined in (\ref{inj rad}), then
\begin{align}
|\nab\beta|^2\lesssim g^{ij}\gamma^{kl}\gamma^{IJ}\nab_k\beta_{Ii}\nab_l\beta_{Jj}+|\di \beta|^2+|\curl\beta|^2,\\\label{Hodge}
\int_{\Omega}|\nab\beta|^2\,d\mu_g\lesssim\int_{\Omega}(N^iN^jg^{kl}\gamma^{IJ}\nab_k\beta_{Ii}\nab_l\beta_{Jj}+|\di \beta|^2+|\curl\beta|^2+K^2|\beta|^2)\,d\mu_g.
\end{align}
\begin{proof}
See \cite{CL}; we also refer Section \ref{section 8} for the weighted version.
\end{proof}

\prop \label{ellpitic estimate}(Elliptic estimates) Let $q:\Omega\to\RR$ be a smooth function. Suppose that $|\theta|+|\frac{1}{l_0}|\leq K$, then we have, for any $r\geq2$ and $\delta>0$,
\begin{align}
||\nab^r q||_{\llb}\lesssim_{K}\sum_{s\leq r}||\Pi\nab^s q||_{\llb}+\sum_{s\leq r-1}||\nab^s\lap q||_{\lli} + ||\nab q||_{\lli},\label{ell est I}\\
||\nab^r q||_{\lli} \lesssim_{K}\delta\sum_{s\leq r}||\Pi\nab^s q||_{\llb}+\delta^{-1}\sum_{s\leq r-2}||\nab^s\lap q||_{\lli}+ \delta^{-1}||\nab q||_{\lli}.\label{ell est II}
\end{align}
\begin{proof} 
See \cite{CL}; we also refer Section \ref{section 8} for the weighted version.
\end{proof}
\rmk We recall that if $\vol\Omega< \infty$, we have $||\nab q||_{\lli}\leq C(\vol\Omega)||\lap q||_{\lli}$.

\subsection{Estimate for the projection of a tensor to the tangent space of the boundary}
The use of the projection of the tensor $\Pi\nab^sD_t^kh$ in the boundary part of energy (\ref{Er}) is essential to compensate the potential loss of regularity. A simple observation that will help us is that if $q=0$ on $\p\Omega$, then $\Pi\nab^2q$ contains only first-order derivative of $q$ and all components of the second fundamental form. To be more precise, we have
\begin{equation}
\Pi\nab^2q=\cnab^2q+\theta\nab_Nq ,\label{21}
\end{equation}
where the tangential component $\cnab^2q=0$ on the boundary. Furthermore, in $L^2$ norms, (\ref{21}) yields,
\begin{equation}
||\Pi\nab^2q||_{\llb}\leq |\theta|_{L^{\infty}(\p\Omega)}||\nab_Nq||_{\llb}.\label{22}
\end{equation}
To prove (\ref{21}), we first recall the components of the projection operator $\gamma_{i}^{j}=\delta_{i}^{j}-N_iN^j$, hence
$$
\gamma_{j}^{k}\nab_i\gamma_{k}^{l}=-\gamma_{j}^{k}\nab_i(N_kN^l)=-\gamma_j^k\theta_{ik}N^l-\gamma_j^kN_k\theta_i^l=-\theta_{ij}N^l,
$$
and so
\begin{dmath*}
\cnab_i\cnab_jq=\gamma_{i}^{i'}\gamma_{j}^{j'}\nab_{i'}\gamma_{j'}^{j''}\nab_{j''}q
=\gamma_i^{i'}\gamma_j^{j'}\gamma_{j'}^{j''}\nab_{i'}\nab_{j''}q+\gamma_i^{i'}\gamma_{j}^{j'}(\nab_{i'}\gamma_{j'}^{j''})\nab_{j''}q
=\gamma_{i}^{i'}\gamma_{j}^{j'}\nab_{i'}\nab_{j'}q-\theta_{ij}\nab_Nq.
\end{dmath*}
In general, the higher order projection formula is of the form
\begin{align}
\Pi \nab^r q = (\cnab^{r-2}\theta)\nab_N q+O(\nab^{r-1}q)+O(\cnab^{r-3}\theta). \label{tensor expression}
\end{align}
which suggests the following generalization of (\ref{22}), its detailed proof can be found in \cite{CL}.
\prop \label{tensor estimate}(Tensor estimate)
Suppose that $|\theta|+|\frac{1}{l_0}|\leq K$, and for $q=0$ on $\p\Omega$, then for $m=0,1$
\begin{multline}
||\Pi\nab^{r}q||_{\llb}\lesssim_{K} ||(\cnab^{r-2}\theta)\nab_{N}q||_{\llb}+\sum_{l=1}^{r-1}||\nab^{r-l}q||_{\llb}\\
 +(||\theta||_{\linf}+\sum_{0\leq l\leq r-2-m}||\cnab^{l}\theta||_{\llb})(\sum_{0\leq l\leq r-2+m}||\nab^l q||_{\llb}), \label{tensor est}
\end{multline}
where the second line drops for $0\leq r\leq 4$.
\subsection{Estimate for the second fundamental form}
The estimate of the second fundamental form is a direct consequence of (\ref{tensor expression}) with $q=h$ together with the physical sign condition,e.g., $|\nab_N h|\geq \epsilon>0$. 
\prop \label{theta estimate}($\theta$ estimate) 
 Suppose that $|\theta|+|\frac{1}{l_0}|\leq K$, and the physical sign condition $|\nab_N h|\geq \epsilon>0$ holds, then
\begin{multline}
||\cnab^{r-2}\theta||_{\llb}\lesssim_{K, \frac{1}{\epsilon}}||\Pi \nab^rh||_{\llb}+\sum_{s=1}^{r-1}||\nab^{r-s}h||_{\llb}\\
+(||\theta||_{\linf}+\sum_{s\leq r-3}||\cnab^s\theta||_{\llb})\sum_{s\leq r-1}||\nab^s h||_{\llb}, \label{theta est}
\end{multline}
where the second line drops for $0\leq r\leq 4$.
\section{The wave equation} \label{section 4}
In this section we study the wave equation satisfied by $h$, obtained by commutating $D_t$ with the first equation of (\ref{EE})
\begin{align}
D_t^2e(h)-\lap h=(\p_iv^j)(\p_jv^i),\q \text{in}\q [0,T]\times \Omega \label{WW}
\end{align}
with initial and boundary conditions
\begin{align}
h|_{t=0} = h_0,\q D_th|_{t=0} = h_1,
\end{align}
and
\begin{align}
h|_{\p\Omega}=0.
\end{align}
Here, $\lap = \delta^{ij} \p_i\p_j=\frac{1}{\sqrt{|\det g|}}\p_a(\sqrt{|\det g|}g^{ab}\p_b)$. In order to express the higher order tensor products in a more appealing way, we adopt the following notation introduced in \cite{CL}.
\mydef (Symmetric dot product)  Let $[D_t,\p]=-(\p v)\symdot\p$, where the symmetric dot product $(\p v)\symdot\p$ is define component-wisely by $((\p v)\symdot\p)_i=\p_i v^k \p_k$. In general, we have 
\begin{equation}
[D_t,\p^r]=\sum_{s=0}^{r-1}\p^s[D_t,\p]\p^{r-s-1}=\sum_{s=0}^{r-1}-
\begin{pmatrix}
r\\
s+1
\end{pmatrix}
(\p^{1+s}v) \symdot \p^{r-s},  \label{Dtd^r comm}
\end{equation}
where 
$$
((\p^{1+s}v) \symdot \p^{r-s})_{i_1,\cdots,i_r}=\frac{1}{r!}\sum_{\sigma\in S_r}(\p^{1+s}_{i_{\sigma_1}\cdots i_{\sigma_{1+s}}}v^k)(\p^s_{k,i_{\sigma_{s+2}\cdots i_{\sigma_{r}}}}),
$$
where $S_r$ is the $r$-symmetric group.\\
\indent In addition, the commutators between $D_t^k$ for $k\geq 2$ and spatial derivatives can be expressed as
\begin{align}
[\p, D_t^k] = \sum_{l_1+l_2=k-1}c_{l_1,l_2}(\p D_t^{l_1}v)\symdot(\p D_t^{l_2})+ \sum_{l_1+\cdots+ l_n= k-n+1, \, n\geq 3} d_{l_1,\cdots,l_n}(\p D_t^{l_1} v)\cdots (\p D_t^{l_{n-1}} v) (\p D_t^{l_n}),\label{dDt^k comm}
\end{align}
and
\begin{multline}
[\lap, D_t^{r-1}] = \sum_{l_1+l_2=r-2}c_{l_1,l_2} (\lap D_t^{l_1}v)\cdot (\p D_t^{l_2})+ \sum_{l_1+l_2=r-2}c_{l_1,l_2} (\p D_t^{l_1}v)\cdot (\p^2 D_t^{l_2})\\
+ \sum_{l_1+\cdots+ l_n= r-n, \, n\geq 3} d_{l_1,\cdots,l_n}(\p D_t^{l_3} v)\cdots (\p D_t^{l_{n}} v)\cdot (\lap D_t^{l_1} v)\cdot(\p D_t^{l_2})\\
+ \sum_{l_1+\cdots+ l_n= r-n, \, n\geq 3} e_{l_1,\cdots,l_n}(\p D_t^{l_3} v)\cdots (\p D_t^{l_{n}} v)\cdot (\p^2 D_t^{l_1} v)\cdot(\p D_t^{l_2})\\
+\sum_{l_1+\cdots+ l_n= r-n, \, n\geq 3} f_{l_1,\cdots,l_n}(\p D_t^{l_3} v)\cdots (\p D_t^{l_n} v)\cdot (\p D_t^{l_1} v)\cdot (\p^2 D_t^{l_2}), \label{lapDt comm}
\end{multline}
where the regular dot product is defined to be the trace of the symmetric dot. 
\subsection{The Energies $W_r(t)$}
By commutating $D_t^{r-1}$ on both sides of (\ref{WW}), we obtain the higher order wave equation
\begin{align}
e'(h)D_t^{r+1}h-\lap D_t^{r-1}h = f_r+g_r, \label{WWr}
\end{align}
where \begin{align}
f_r = D_t^{r-1}(\p v\cdot\p v)+[D_t^{r-1},\lap]h,
\end{align}
and $g_r$ is sum of terms of the form 
\begin{align}
e^{(m)}(h)(D_t^{i_1}h)\cdots(D_t^{i_m}h), \q 2\leq m\leq r,\q i_1+\cdots+i_m=r+1,\q 1\leq i_1\leq\cdots\leq i_m\leq r. \label{g_r}
\end{align}
Now, let us define the energy 
\begin{align}
W_r(t) = \frac{1}{2}||\sqrt{e'(h)}D_t^rh||_{\lli}+\frac{1}{2}||\nab D_t^{r-1}h||_{\lli} \label{Wr}
\end{align}
and by the standard energy estimates for the wave equations together with (\ref{dDt^k comm}), we have
\thm \label{energy est for wave eq} Let $W_r$ be defined as in (\ref{Wr}), we have
\begin{align}
\frac{dW_r}{dt} \lesssim W_r + ||f_r||_{\lli}+||g_r||_{\lli}. \label{est Wr}
\end{align}
\pf This follows from standard energy estimates for the wave equation. We refer \cite{LL} Section 4 for the detail.
\subsection{Estimates for $||f_r||_{\lli}$}
By adopting our notations used in (\ref{dDt^k comm})-(\ref{lapDt comm}), we are able to express $f_r$ as
\begin{dmath}
f_r = \sum_{l_1+l_2=r-1}c_{l_1,l_2}(\nab D_t^{l_1}v)\cdot(\nab D_t^{l_2}v)+
 \sum_{l_1+l_2=r-2}d_{l_1,l_2}(\lap D_t^{l_1}v)\cdot (\nab D_t^{l_2}h)\\
  +\sum_{l_1+l_2=r-2}e_{l_1,l_2}(\nab D_t^{l_1}v)\cdot(\nab^2 D_t^{l_2}h) +\text{error terms}, \label{fr}
\end{dmath}
where the "error terms" refer to the terms generated by the commutators, which are of the form
\begin{multline}
e_r=\sum_{l_1+\cdots+ l_n= r+1-n, \, n\geq 3} g_{l_1,\cdots,l_n}(\p D_t^{l_3} v)\cdots (\p D_t^{l_{n}} v)\cdot (\p D_t^{l_1} v)\cdot(\p D_t^{l_2}v)\\
+ \sum_{l_1+\cdots+ l_n= r-n, \, n\geq 3} e_{l_1,\cdots,l_n}(\p D_t^{l_3} v)\cdots (\p D_t^{l_{n}} v)\cdot (\p^2 D_t^{l_1} v)\cdot(\p D_t^{l_2} h)\\
+\sum_{l_1+\cdots+ l_n= r-n, \, n\geq 3} f_{l_1,\cdots,l_n}(\p D_t^{l_3} v)\cdots (\p D_t^{l_n} v)\cdot (\p D_t^{l_1} v)\cdot (\p^2 D_t^{l_2}h).
\end{multline}
We need to estimate $||f_r||_{\lli}$ and $||g_r||_{\lli}$ for $r\geq 1$. Since our estimates include mixed space-time derivatives, we would like to use the following more appealing notations.
\mydef(Mixed Sobolev norms)
let $u(t,\cdot):\RR^n \to \RR$ be a smooth function. We define
\begin{align*}
||u||_{r,\,0} = \sum_{s+k=r,\, k<r}||\nab^s D_t^ku||_{\lli},\\
||u||_r = ||u||_{r,0}+||\sqrt{e'(h)}D_t^rh||_{\lli}.
\end{align*}
We have to make sure that the $r$-th order Sobolev norms in our estimates for $||f_r||_{\lli}$, $r\geq 3$ do not include $||\nab^r h||_{\lli}$ and $||\nab ^r v||_{\lli}$. This is because that we need to control $||f_{r+1}||_{\lli}, r\geq 2$ by $\sqrt{E_r^*}$ in Section \ref{section 5.5}, but $||\nab^{r+1} h||_{\lli}$ and $||\nab ^{r+1} v||_{\lli}$ can only be controlled by $\sqrt{E_{r+1}^*}$. 
\subsubsection{When r=1,2,3}
Since $f_1 = (\nab v)\cdot (\nab v)$, we have
$$
||f_1||_{\lli} \lesssim_M ||\nab v||_{\lli}.
$$
The bounds for $||f_2||_{\lli}$ and $||f_3||_{\lli}$ are the same as what we did in our previous work \cite{LL}, and so we shall only record results here. 
$$||f_2||_{\lli}\lesssim_M ||\nab^2 v||_{\lli}+||\nab^2 h||_{\lli}+||\nab v||_{\lli}.$$
$$||f_3||_{\lli}\lesssim_M ||\nab^2 D_th||_{\lli}+||e'(h)\nab D_t^2h||_{\lli}+||h||_{2,0}+\sum_{j=1,2}||\nab^j v||_{\lli}.$$

\subsubsection{When r=4}
The bounds for the first and the third term of $f_4$ is the same as in \cite{LL}.
\begin{dmath*}
\sum_{l_1+l_2=3}||c_{l_1,l_2}(\nab D_t^{l_1}v)(\nab D_t^{l_2}v)||_{\lli}+\sum_{l_1+l_2=2}||e_{l_1,l_2}(\nab D_t^{l_1}v)(\nab^2 D_t^{l_2}h)||_{\lli}\\
\lesssim_M ||\nab^2D_t^2h||_{\lli}+\sum_{j=2,3}||h||_{j,0}+||\nab^2v||_{\lli} + ||\nab D_t^2v||_{\lli}.
\end{dmath*}
But we cannot use interpolation to bound $||\lap v\cdot \nab D_t^2h||_{\lli}$ involved in the second term of $f_4$, as $|D_t^2h|$ is no longer part of the a priori assumptions. But since
\begin{equation}
\lap v = \nab \di v+\nab\cdot\curl v, 
\end{equation}
and since $|e''(h)| \leq c_0|e'(h)|$, 
\begin{align}
|\nab \di v| \lesssim |e'(h)(\nab h)D_th|+|e'(h)\nab D_th| 
\end{align}
is bounded by a priori assumptions (\ref{v,h}) and (\ref{D_te(h)}). On the other hand \footnotemark, since $|\nab \cdot \curl v| \leq M$ as well, we conclude
\begin{dmath}
||\lap v\cdot \nab D_t^2h||_{\lli} \lesssim_M ||\nab D_t^2h||_{\lli},
\end{dmath}
and so
\begin{equation}
\sum_{l_1+l_2=2}||(\lap D_t^{l_1}v)\cdot (\nab D_t^{l_2}h)||_{\lli}\lesssim_M ||\nab D_t^2h||_{\lli}+||\nab^3 D_th||_{\lli}+\sum_{j=2,3}(||\nab^j v||_{\lli}+||\nab^j h||_{\lli}).
\end{equation}
Most of the terms in $e_4$ can be bounded by corresponding terms in $f_r$, for $r\leq 4$, and similar terms in $e_3$ times a priori assumptions, apart from terms of the form $\nab v\cdot\nab^2D_tv\cdot\nab h$, whose $L^2$ norm can be bounded by $||\nab^3h||_{\lli}$.

Therefore, we sum up and get
\begin{dmath*}
||f_4||_{\lli}\lesssim_{M} ||\nab^3 D_th||_{\lli}+||\nab^2D_t^2h||_{\lli}
+\sum_{j=2,3}(||h||_{j,0}+||\nab^j v||_{\lli}).
\end{dmath*}

\footnotetext{One could alternatively estimate $||\lap v\cdot \nab D_t^2h||_{\lli}$ by Sobolev lemma, e.g., 
\begin{align*}
||\lap v\cdot \nab D_t^2h||_{\lli} \lesssim (\sum_{j=2,3}||\nab^j v||_{\lli})(\sum_{j=1,2}||\nab^j D_t^2h||_{\lli}).
\end{align*}
However, \eqref{intro main1} then fails to be linear in $E_r^*$.
}

\subsubsection{When $r=5$ and $n\leq 4$}
The bounds for the first and the third terms of $f_5$ remains unchanged as in \cite{LL}.
 \begin{multline*}
\sum_{l_1+l_2=4}||(\nab D_t^{l_1}v)(\nab D_t^{l_2}v)||_{\lli}+\sum_{l_1+l_2=3}||(\nab D_t^{l_1}v)(\nab^2 D_t^{l_2}h)||_{\lli}\\
\lesssim_{K,M} ||\nab^2D_t^3h||_{\lli}+\sum_{1\leq i\leq 4}||v||_{i,0}+\sum_{2\leq i\leq4}||h||_{i,0}.
\end{multline*}
As for the term $\sum_{l_1+l_2=3}d_{l_1,l_2}||(\lap D_t^{l_1}v) (\nab D_t^{l_2}h)||_{\lli}$, we need the Sobolev lemma (\ref{interior sobolev}) to bound $||\lap v\cdot \nab D_t^3h||_{\lli}$ and $||\lap D_t v\cdot\nab D_t^2h||_{\lli}$ as in \cite{LL}, i.e., 
\begin{equation}
 ||\lap v\cdot \nab D_t^3h||_{\lli}\lesssim_K (\sum_{j=2,3}||\nab^jv||_{\lli})(\sum_{j=1,2}||\nab^jD_t^3h||_{\lli}) \label{mod f_5},
\end{equation}
\begin{equation*}
||\lap D_t v\cdot\nab D_t^2h||_{\lli}\lesssim_{K}(\sum_{j=3,4}||\nab^jh||_{\lli})(\sum_{j=1,2}||\nab^jD_t^2h||_{\lli}) ,
\end{equation*}
\begin{equation*}
||\lap D_t^2v\cdot \nab D_th||_{\lli}\lesssim_M ||\nab^2D_t^2 v||_{\lli}\lesssim_M ||\nab^3 D_th||_{\lli}+\sum_{j\leq 3}(||\nab^j v||_{\lli}+||\nab^jh||_{\lli}),
\end{equation*}
and
 \begin{dmath*}||\lap D_t^3v\cdot\nab h||_{\lli}\lesssim_M ||\nab\lap D_t^2h||_{\lli}+||\lap[D_t^2,\nab]h||_{\lli}\\
\lesssim_{K,M} ||\nab^3D_t^2h||_{\lli}+\sum_{j=2,3}||v||_{i,0}+\sum_{j=3,4}||h||_{i,0},
\end{dmath*}
respectively. Most of the terms in the error term $e_5$ are essentially bounded by corresponding terms in $f_r$, for $r\leq 5$, and similar terms in $e_3$ and $e_4$ times a priori assumptions, apart from the terms of the form $\nab v\cdot\nab^2D_t^2v\cdot\nab h$, which is estimated by $||\nab^2 D_t^2v||_{\lli}$.
Hence,
\begin{dmath*}
||f_5||_{\lli}\lesssim_{K,M} ||\nab^3D_t^2h||_{\lli}+(\sum_{j=2,3}||\nab^jv||_{\lli})(\sum_{j=1,2}||\nab^jD_t^3h||_{\lli})\\
+(\sum_{j=3,4}||\nab^jh||_{\lli})(\sum_{j=1,2}||\nab^jD_t^2h||_{\lli})
+\sum_{1\leq i\leq 4}||v||_{i,0}+\sum_{2\leq i\leq4}||h||_{i,0}.
\end{dmath*}
\subsection{When $r\geq 6$}
The commutator \eqref{dDt^k comm} in fact implies that
\begin{align}
D_t^k v = -\p D_t^{k-1}h + c_{\alpha'\beta'\gamma'}(\p^{\alpha_1'} v)\cdots(\p^{\alpha_m'} v)(\p^{\beta_1'}D_t^{\gamma_1'}h)\cdots(\p^{\beta_n'}D_t^{\gamma_n'}h) \label{D_t^k v},
\end{align}
where
\begin{align*}
\alpha' = (\alpha_1',\cdots,\alpha_m'),\,\, \beta' = (\beta_1',\cdots,\beta_n'),\,\,\gamma' =(\gamma_1',\cdots,\gamma_n'),\\
\alpha_1'+\cdots+\alpha_m'+(\beta_1'+\gamma_1')+\cdots+(\beta_n'+\gamma_n') = k,\\
1\leq \alpha_i'\leq k-2,\q \text{when}\q k\geq 3,\\ 
1\leq \beta_j' \leq k-2,\q \text{when}\q k\geq 4.
\end{align*}
Because of this and \eqref{fr}, we can re-express $f_r,r\geq 6$ as 
\begin{align}
f_r =c_{\alpha\beta\gamma}(\p^{\alpha_1} v)\cdots(\p^{\alpha_m} v)(\p^{\beta_1}D_t^{\gamma_1}h)\cdots(\p^{\beta_n}D_t^{\gamma_n}h) \label{fr r>6},
\end{align}
where 
\begin{align*}
\alpha = (\alpha_1,\cdots,\alpha_m),\,\, \beta = (\beta_1,\cdots,\beta_n),\,\,\gamma =(\gamma_1,\cdots,\gamma_n),\\
\alpha_1+\cdots+\alpha_m+(\beta_1+\gamma_1)+\cdots+(\beta_n+\gamma_n) = r+1,\\
1\leq \alpha_i\leq r-2, \q 1\leq i\leq m,\\
1\leq \beta_j+\gamma_j \leq r, \q 1\leq j\leq n.
\end{align*}
In addition to these, there exists at most one $i$ or $j$ such that $\alpha_i = r-2$ or $\beta_j+\gamma_j \geq r-2$, and further if $\beta_j+\gamma_j \geq r-1$, we must have $\gamma_j \geq 1$. Thus, $f_r$ never consists terms of the form $(\p^2 v)(\p^{r-1}h)$ if $r\geq 6$.

Since $f_r$ is a sum of products of the form \eqref{fr r>6}, we apply the following derivative counting method on each product to estimates $||f_r||_{\lli}$. 
\begin{itemize}
\item If $\alpha_i \geq r-2$ for some $i$ or $\beta_j+\gamma_j\geq r-2$ for some $j$, then there are at most four terms involved in the product \eqref{fr r>6}, among which at least one must satisfy a priori assumptions \eqref{geometry_bound}-\eqref{D_te(h)} if the product has more than two terms. Hence,
\begin{dmath}
||(\nab^{\alpha_1} v)\cdots(\nab^{\alpha_m} v)(\nab^{\beta_1}D_t^{\gamma_1}h)\cdots(\nab^{\beta_n}D_t^{\gamma_n}h)||_{\lli} \\
\leq C_r(K,M, \sum_{k\leq r-2}||\nab^k v||_{\lli},\sum_{k\leq r-2}||h||_{k,0})(\sum_{k\leq r-1}||\nab^k v||_{\lli}+\sum_{k\leq r-1}||h||_{k,0}+\sum_{k\leq r-1}||D_th||_{k,0}) \label{fr counting}.
\end{dmath}
Here we have used the Sobolev lemma
\begin{equation}
||u_1\cdots u_N||_{L^2} \leq C(K) ||u_1||_{H^1}\cdots||u_N||_{H^1},\q N=2,3 \label{Sobolev N leq 3}.
\end{equation}

Now, we assume $\alpha_i\leq r-3$ and $\beta_j+\gamma_j \leq r-3$ for all $i,j$.

\item If $\alpha_i <r-3$ and $\beta_j+\gamma_j < r-3$ for all $i,j$, then 
\begin{dmath}
||(\nab^{\alpha_1} v)\cdots(\nab^{\alpha_m} v)(\nab^{\beta_1}D_t^{\gamma_1}h)\cdots(\nab^{\beta_n}D_t^{\gamma_n}h)||_{\lli} \\
\leq C_r(K,M, \sum_{k\leq r-2}||\nab^k v||_{\lli},\sum_{k\leq r-2}||h||_{k,0})\label{fr counting lo}.
\end{dmath}
Here we have used the Sobolev lemma
\begin{equation}
||u_1\cdots u_N||_{L^2} \leq C(K) ||u_1||_{H^2}\cdots||u_N||_{H^2},\q N\geq 4 \label{Sobolev N geq 4}.
\end{equation}
\item If $\alpha_i = r-3$ for some $i$ and/or $\beta_j+\gamma_j =r-3$ for some $j$,  then there exists at most one $i'\neq i$ or $j'\neq j$ such that $\alpha_{i'} = r-3$ or $\beta_{j'}+\gamma_{j'} = r-3$. In this case, the product consists at most $3$ terms. Hence, \eqref{fr counting lo} remains valid in this case by Sobolev lemma. 

\end{itemize}
Therefore,  one concludes that when $r\geq 6$,
\begin{dmath}
||f_r||_{\lli}
\leq C_r(K,M, \sum_{k\leq r-2}||\nab^k v||_{\lli},\sum_{k\leq r-2}||h||_{k,0})\cdot\\
(\sum_{k\leq r-1}||\nab^k v||_{\lli}+\sum_{k\leq r-1}||h||_{k,0}+\sum_{k\leq r-1}||D_th||_{k,0}),
\end{dmath}
where $C_r$ are continuous functions.

\subsection{Estimates for $||g_r||_{\lli}$}
We recall that $e(h)= \log\rho(h)$ which satisfies 
\begin{enumerate}
\item $|e^{(k)}(h)|\leq c_0$.
\item $|e^{(k)}(h)|\leq c_0\sqrt{e'(h)}$.
\item $|e^{(k)}(h)|\leq c_0|e'(h)|^k$.
\end{enumerate}
\subsection{When $r=1,2,3,4$}
For each $r$, $g_r$ is a sum of terms of the form
\begin{align}
e^{(m)}(h)D_t^{j_1}h\cdots D_t^{j_m}h, \q j_1+\cdots+j_m=r+1,\q 1\leq j_1\leq\cdots\leq j_m\leq r \label{gr},
\end{align}
and $j_i\leq 2$ for $i\leq m-1$. Therefore, the a priori assumption (\ref{D_te(h)}) yields
$$
||e^{(m)}(h)D_t^{j_1}h\cdots D_t^{j_m}h||_{\lli}\lesssim ||(e'(h)D_t^{j_1}h)\cdots (e'(h)D_t^{j_m}h)||_{\lli} \lesssim_{M,c_0} ||e'(h)D_t^{j_m}h||_{\lli}.
$$
Hence we conclude
\begin{align}
||g_r||_{\lli}\lesssim_{M,c_0}  \sum_{j\leq r} ||e'(h)D_t^jh||_{\lli},\q r\leq 4 \label{g_r mod 1}.
\end{align}
\subsection{When $r=5$}
The only difference for estimating $g_5$ is that it contains a quadratic term $e''(h)D_t^3h\cdot D_t^3h$, whose $L^2$ norm is bounded via Sobolev lemma (\ref{interior sobolev}). We have
\begin{dmath*}
||e''(h)(D_t^3h)^2||_{\lli} \lesssim_{c_0} ||(e'(h)D_t^3h)^2||_{\lli}
\lesssim_{K,c_0} (\sum_{j=0,1}||\nab^j(e'(h)D_t^3h)||_{\lli})^2\\
\lesssim_{K,M,c_0} (|e'(h)|\cdot||\nab D_t^3h||_{\lli}+||e'(h)D_t^3h||_{\lli})^2
\end{dmath*}
Hence we conclude
\begin{align}
||g_5||_{\lli}\lesssim_{K,M,c_0} \sum_{j\leq 5} ||e'(h)D_t^jh||_{\lli} + (|e'(h)|\cdot||\nab D_t^3h||_{\lli}+||e'(h)D_t^3h||_{\lli})^2 .\label{g_r mod 2}
\end{align}

\subsection{When $r\geq 6$}
The estimates for the general case in fact follow from the case when $r=5$. Since $g_r$ is a sum of the products of the form \eqref{gr}, we apply the derivative counting method again on estimating each of the products. 
\begin{itemize}
\item 
If $j_m \geq r-2$, then the product consists of at most $4$ terms, where $j_i< r-2$ for all $i<m$, among which at least one must be of order no more than $2$, i.e., they are of the form $D_t^{j_l}h$ with $j_l \leq 2$ (and so $e'(h)D_t^{j_l}h$ satisfies \eqref{D_te(h)}). Hence, by Sobolev lemma \eqref{Sobolev N leq 3}, 
\begin{dmath}
||e^{(m)}(h)D_t^{j_1}h\cdots D_t^{j_m}h||_{\lli} \lesssim_{c_0} ||(e'(h)D_t^{j_1}h)\cdots(e'(h)D_t^{j_m}h)||_{\lli}\\
\leq C_r(K,M,c_0, \sum_{k\leq r-2}||D_th||_{k})\sum_{k\leq r-1}||D_th||_{k}. \label{gr counting}
\end{dmath}

\item If $j_m < r-3$, then by \eqref{Sobolev N geq 4} we have
\begin{dmath}
||e^{(m)}(h)D_t^{j_1}h\cdots D_t^{j_m}h||_{\lli}
\leq C_r(K,M,c_0,\sum_{k\leq r-2} ||D_th||_{k}).\label{gr counting lo}
\end{dmath}
\item If $j_m = r-3$, then there exists at most one $j_l$, where $l<m$ such that $j_l = r-3$, and the product consists of at most $3$ terms if this is the case. Hence, \eqref{gr counting lo} holds by Sobolev lemma \eqref{Sobolev N leq 3}. 
\end{itemize}

Therefore, one concludes that when $r\geq 6$, 
\begin{dmath}
||g_r||_{\lli}
\leq C_r(K,M,c_0, \sum_{k\leq r-2}||D_th||_{k})\sum_{k\leq r-1}||D_th||_{k},
\end{dmath}
where $C_r$ are continuous functions. \\

In summary, we have proved:

\thm Let $f_r$ and $g_r$ be defined as \eqref{fr} and \eqref{gr}, respectively. Then we have the estimates
\begin{dmath}
||f_r||_{\lli} \leq C(M)(||\nab^r v||_{\lli}+||\nab^r h||_{\lli}),\q r=1,2
\end{dmath}
\begin{dmath}
||f_r||_{\lli} \leq C_r(K,M, \sum_{k\leq r-2}||\nab^k v||_{\lli},\sum_{k\leq r-2}||D_th||_{k,0})\cdot\\(\sum_{k\leq r-1}||\nab^k v||_{\lli}+\sum_{k\leq r-1}||h||_{k,0}+\sum_{k\leq r-1}||D_th||_{k,0}),\q 3\leq r\leq 5\label{fr est r=3,4,5} 
\end{dmath} 
\begin{dmath}
||f_r||_{\lli} \leq C_r(K,M, \sum_{k\leq r-2}||\nab^k v||_{\lli},\sum_{k\leq r-2}||h||_{k,0})\cdot\\(\sum_{k\leq r-1}||\nab^k v||_{\lli}+\sum_{k\leq r-1}||h||_{k,0}+\sum_{k\leq r-1}||D_th||_{k,0}),\q r\geq 6 \label{fr est r>5}
\end{dmath} 
and
\begin{align}
||g_r||_{\lli} \leq C(M,c_0)\sum_{j\leq r}||\sqrt{e'(h)}D_t^jh||_{\lli}, \q 1\leq r\leq 4 \label{gr est r<4}\\
||g_r||_{\lli} \leq C_r(K,M,c_0, \sum_{k\leq r-2}||D_th||_{k})\sum_{k\leq r-1}||D_th||_{k}, \q r\geq 5\label{gr est r>4}
\end{align}

\subsection{Improved estimates for $||f_r||_{\lli}$ and $||g_r||_{\lli}$}
\mydef (Improved mixed norms) \label{improved mixed norms}
\begin{itemize}
\item $||h||_{r,1,0} := \sum_{k+s=r, k<r-1}||\nab^s D_t^k h||_{\lli} + ||\sqrt{e'(h)}\nab D_t^{r-1}h||_{\lli}$,
\item $||h||_{r,1} := ||h||_{r,1,0} + ||e'(h)D_t^rh||_{\lli}$.
\end{itemize}
Under these new norms, the estimates for $||f_r||_{\lli}$ and $||g_r||_{\lli}$ can be improved as:
\thm \label{improved fr and gr}
 Let $f_r$ and $g_r$ be defined as \eqref{fr} and \eqref{gr}, respectively. Then, 
\begin{dmath}
||f_r||_{\lli} \leq C_r(K,M, \sum_{k\leq r-2}||\nab^k v||_{\lli},\sum_{k\leq r-2}||D_th||_{k,0})\cdot\\(\sum_{k\leq r-1}||\nab^k v||_{\lli}+\sum_{k\leq r-1}||h||_{k,0}+\sum_{k\leq r-1}||D_th||_{k,1,0}),\q 3\leq r\leq 5
\end{dmath} 
\begin{dmath}
||f_r||_{\lli} \leq C_r(K,M, \sum_{k\leq r-2}||\nab^k v||_{\lli},\sum_{k\leq r-2}||h||_{k,0})\cdot\\(\sum_{k\leq r-1}||\nab^k v||_{\lli}+\sum_{k\leq r-1}||h||_{k,0}+\sum_{k\leq r-1}||D_th||_{k,1,0}),\q r\geq 6
\end{dmath} 
and
\begin{align}
||g_r||_{\lli} \leq C(M,c_0)\sum_{j\leq r}||e'(h)D_t^jh||_{\lli}, \q 1\leq r\leq 4 \\
||g_r||_{\lli} \leq C_r(K,M,c_0, \sum_{k\leq r-2}||D_th||_{k,1})\sum_{k\leq r-1}||D_th||_{k,1}, \q r\geq 5
\end{align}
\begin{proof}
It is easy to observe that the estimates for $||f_r||_{\lli}$ and $||g_r||_{\lli}$ does not include the quantity $||\nab D_t^{r-1}h||_{\lli}$, and we no longer use $e'(h)\leq c_0\sqrt{e'(h)}\leq c_0$ in the estimates for $||g_r||_{\lli}$; in other words, we keep $e'(h)$ whenever it is possible.
\end{proof}

Theorem \ref{improved fr and gr} is essential for estimating the lower order terms $||\nab D_t^k h||_{\lli}, 0\leq k\leq r-1$ without using the wave equation (See Section \ref{section 6}).

\section{Energy estimates for Euler equations with free boundary}\label{section 5}

\prop \label{a priori estimate}Let $E_r$ be defined as (\ref{Er}), then there are continuous functions $C_r$  such that, for $t\in[0,T]$,
\begin{align}
|\frac{dE_r(t)}{dt}|\leq C_r(K,\frac{1}{\epsilon},M,c_0, E_{r-1}^*)E_r^*(t), \label{main 1}
\end{align}
holds for all $r\geq 1$, where $E_r^*=\sum_{i\leq r}E_i$, provided the assumption \eqref{e_kk} and the a priori bounds (\ref{geometry_bound})-(\ref{D_te(h)}).

\subsection{Computing $\frac{d}{dt}E_r$}

Our computation for $\frac{d}{dt}E_r(t)$ is almost identical to what we had in \cite{LL}, but we can no longer use the interpolation (\ref{int interpolation}) which relies on the $L^{\infty}$ bounds of $|D_th|$ and $|D_t^2h|$. We first compute $\frac{d}{dt}E_{s,k}$, where $E_{s,k}$ is defined as \eqref{Esk}, when $s>0$.
\begin{multline}
\frac{d}{dt}E_{s,k} = \frac{1}{2}\int_{\DD_t}\rho D_t(\delta^{ij}Q(\p^sD_t^kv_i,\p^sD_t^{k}v_j)\dx
+\frac{1}{2}\int_{\DD_t}\rho D_t(e'(h)Q(\p^sD_t^k h,\p^s D_t^kh))\dx\\
+\frac{1}{2}\int_{\p\DD_t}\rho D_t(Q(\p^sD_t^k h,\p^sD_t^k h)\nu)-Q(\p^sD_t^k h,\p^sD_t^k h)\nu(\sigma v\cdot N)+\rho Q(\p^sD_t^k h,\p^sD_t^k h)D_t\nu\,dS. \label{DtEr'}
\end{multline}
The estimates (\ref{extending nab n})-(\ref{Dtgamma}) together with a priori assumptions imply\footnote{We refer Section 5 of \cite{CL}  for the detailed proof} $$|D_tq^{ij}|\lesssim M, \q |\p q^{ij}|\lesssim M+K,\q  |\sigma v\cdot N|_{L^{\infty}(\p\Omega)}\lesssim K+M,$$
$$
|D_t\nu|_{L^{\infty}(\p\Omega)} = |D_t(-\nab_N h)^{-1}|_{L^{\infty}(\p\Omega)}\lesssim 1+\frac{1}{M},
$$
and
\begin{align}
D_t\gamma^{ij}=-2\gamma^{im}\gamma^{jn}(\frac{1}{2}D_tg_{mn}) \label{bdy q}.
\end{align}
Since $|D_tq^{ij}|\lesssim M$ in the interior and on the boundary $q^{ij}=\gamma^{ij}$, 
and by \eqref{bdy q} $D_t \gamma$ is tangential, so  that (\ref{DtEr'}) can then be reduced to
\begin{multline}
\frac{d}{dt}E_{s,k} \leq \int_{\DD_t}\rho \delta^{ij}Q(D_t\p^sD_t^kv_i,\p^sD_t^{k}v_j)\dx+\int_{\DD_t}\rho e'(h)Q(D_t\p^sD_t^k h,\p^s D_t^kh)\dx\\
+\int_{\p\DD_t}\rho Q(D_t\p^sD_t^k h,\p^sD_t^k h)\nu\,dS+C(K,M)(E_r+||h||_{j}^2+||v||_{r,0}^2). \label{DtEr}
\end{multline}

Now, if $s\geq 1$, our commutators (\ref{Dtd^r comm}) and (\ref{dDt^k comm}) yield, since $D_tv_i = -\p_i h-e_n$,
\begin{align}
D_t\p^sD_t^k v_i = -\p^sD_t^k\p_i h +\sum_{0\leq m\leq s-1}c_{sr}(\p^{m+1}v)\symdot\p^{s-m}D_t^kv_i \label{comm1},\\
D_t\p^rh + (\p_jh)\p^rv^j = \p^rD_th+\sum_{0\leq m \leq r-2}d_{sr}(\p^{m+1}v)\symdot \p^{r-m}h\label{comm2},\\
D_t\p^sD_t^k h = \p^sD_t^{k+1}h+\sum_{0\leq m\leq s-1}d_{sr}(\p^{m+1}v)\symdot\p^{s-m}D_t^kh,\q \text{for}\q k\geq 1. \label{comm3}
\end{align}
We control the term $||(\p^{m+1}v)\symdot\p^{s-m}D_t^kv_i||_{\lei}$ in (\ref{comm1}) and $||(\p^{m+1}v)\symdot\p^{s-m}D_t^kh||_{\lei}$ in (\ref{comm3}) for $s+k=r$ and $s\geq 1$.
\begin{itemize}
\item  The term $||(\p^{m+1}v)\symdot\p^{s-m}D_t^kh||_{\lei}$ can be bounded by
\begin{enumerate}
\item For $k=0$,
$$
||(\p^{m+1}v)\symdot \p^{r-m}h||_{\lei}\lesssim_K |\p v|_{L^{\infty}}\sum_{j\leq r}||\p^r h||_{\lei}+|\p h|_{L^{\infty}}\sum_{j\leq r} ||\p^r v||_{\lei}.
$$
\item For $k=r-1$ (and so $m=0$), 
$$
||(\p v)\symdot(\p D_t^{r-1}h)||_{\lei} \leq |\p v|_{L^{\infty}}||\p D_t^{r-1}h||_{\lei}.
$$
\item For $1\leq k \leq r-2$ , if $m=0$ then 
$$
||(\p v)\symdot\p^{s}D_t^kh||_{\lei} \leq  |\p v|_{L^{\infty}}||\p^s D_t^{k}h||_{\lei}.
$$
On the other hand, if $m\geq 1$ and $r\geq 4$, we have
$$
||(\p^{m+1}v)\symdot\p^{s-m}D_t^kh||_{\lei} \lesssim_K \sum_{i=1,2}||\p^{m+i}v||_{\lei}\cdot\sum_{j=0,1}||\p^{s-m+j}D_t^kh||_{\lei}.
$$
Here, at most one of $m+2$ or $r-m+1$ can in fact equal to $r$ when $r\geq 4$. However, if $r=3$, then $k$ must equal to $1$, and so $s=2$. Hence,
$$
||\p^{m+1}v\symdot \p^{2-m}D_th||_{\lei} \leq |\p v|_{L^{\infty}}||\p^2D_th||_{\lei}+|\p D_th|_{L^{\infty}}||\p^2 v||_{\lei}.
$$
\end{enumerate}
\item The term $||(\p^{m+1}v)\symdot\p^{s-m}D_t^kv_i||_{\lei}$ can be bounded similarly as above with $h$ replaced by $v$.
\end{itemize}
The above anaylsis shows that the $L^2$ norm of the sum in (\ref{comm1})-(\ref{comm3})  contribute only to $||v||_{r,0}$ and $||h||_{r,0}$. Hence,
\begin{multline}
\frac{d}{dt}E_r \leq -\int_{\DD_t}\rho (\delta^{ij}Q(\p^sD_t^kv_i,\p^sD_t^{k}\p_jh)\dx+\int_{\DD_t}\rho e'(h)Q(\p^sD_t^k h,\p^s D_t^{k+1} h)\dx\\
+\int_{\p\DD_t}\rho Q(\p^sD_t^k h,D_t\p^sD_t^k h)\nu\,dS \\
+ C(K,M)(||v||_{r,0}+||h||_{r,0})(\sum_{i\leq r-1}||v||_{i,0}+||h||_{i,0})(\sum_{i\leq r}||v||_{i,0}+||h||_{i,0}). \label{DtEr''}
\end{multline}
In addition, (\ref{dDt^k comm}) and \eqref{D_t^k v} yield that for $s+k=r$
\begin{dmath}
||\p^sD_t^k\p h-\p^{s+1}D_t^{k}h||_{\lei} \lesssim  \sum_{l_1+l_2=k-1}||\p^s(\p D_t^{l_1}v\symdot\p D_t^{l_2}h)||_{\lei}\\
+ \sum_{l_1+\cdots+ l_n= k-n+1, \, n\geq 3}||\p^s(\p D_t^{l_1} v\cdots \p D_t^{l_{n-1}} v\symdot \p D_t^{l_n}h)||_{\lei}\\
 \lesssim_{K,M}(\sum_{i\leq r-1}||v||_{i,0}+||h||_{i,0})(\sum_{i\leq r}||v||_{i,0}+||h||_{i,0}). \label{error terms}
\end{dmath}

Therefore,

\begin{multline}
\frac{d}{dt}E_r \leq \int_{\DD_t}\rho (\delta^{ij}Q(\p^sD_t^kv_i,\p_j\p^{s}D_t^{k}h)\dx+\int_{\DD_t}\rho e'(h)Q(\p^sD_t^k h,\p^{s} D_t^{k+1} h)\dx\\
+\int_{\p\DD_t}\rho Q(\p^sD_t^k h,D_t\p^sD_t^k h)\nu\,dS\\
+C(K,M)(||v||_{r,0}+||h||_{r,0})(\sum_{i\leq r-1}||v||_{i,0}+||h||_{i,0})(\sum_{i\leq r}||v||_{i,0}+||h||_{i,0}). \label{DtEr final}
\end{multline}
If we integrate by parts in the first term
\begin{multline}
\int_{\DD_t}\rho \delta^{ij}Q(\p^sD_t^k \p_iv_j,\p^sD_t^kh)\dx+\int_{\DD_t}\rho e'(h)Q(\p^sD_t^k h,\p^{s} D_t^{k+1} h)\dx\\
+\int_{\p\DD_t}\rho Q(\p^sD_t^k h,D_t\p^sD_t^k h-\nu^{-1}N_i\p^sD_t^kv^i)\nu\,dS\\
+C(K,M)(||v||_{r,0}+||h||_{r,0})(\sum_{i\leq r-1}||v||_{i,0}+||h||_{i,0})(\sum_{i\leq r}||v||_{i,0}+||h||_{i,0}). \label{int by part}
\end{multline}
But since $\p^sD_t^{k+1}e(h)$ equals $e'(h)\p^s D_t^{k+1}h$ plus a sum of terms of the form
$$
e^{(m)}(h)(\p^{i_1}D_t^{j_1}h)\cdots(\p^{i_m}D_t^{j_m}h),
$$
where
$$
(i_1+j_1)+\cdots+(i_m+j_m)\leq r+1,\q 1\leq i_1+j_1\leq\cdots\leq i_m+j_m\leq r.
$$
Therefore,
\begin{multline}
\int_{\DD_t}\rho \delta^{ij}Q(\p^sD_t^k \p_iv_j,\p^sD_t^kh)\dx=\int_{\DD_t}\rho Q(\p^sD_t^k\di v,\p^sD_t^kh)\dx \\
= -\int_{\DD_t}\rho Q(\p^sD_t^{k+1}e(h),\p^sD_t^kh)\dx
\leq-\int_{\DD_t}\rho e'(h)Q(\p^sD_t^{k+1}h,\p^sD_t^kh)\dx\\
 +C(K,M)(||v||_{r,0}+||h||_{r,0})(\sum_{i\leq r-1}||v||_{i,0}+||h||_{i,0})(\sum_{i\leq r}||v||_{i,0}+||h||_{i,0}),
\end{multline}
so the first integral in (\ref{int by part}) cancels with the second term. \\

We recall $\nu = -(\p_Nh)^{-1}$, so that $\nu^{-1}N_i=\p_ih$. Hence, the boundary term in (\ref{int by part}) becomes
\begin{align}
\sum_{k+s=r,s>0}\int_{\p\DD_t}\rho Q(\p^sD_t^k h,D_t\p^sD_t^k h+(\p_ih)(\p^sD_t^kv^i)\nu\,dS. \label{reduced boundary int}
\end{align}
Now, since (\ref{comm2}) and (\ref{comm3}), (\ref{reduced boundary int}) becomes sum of the boundary inner product of $\Pi\p^sD_t^kh$ and
\begin{align}
\Pi(D_t\p^rh + (\p_jh)\p^rv^j) = \Pi\p^rD_th+\sum_{0\leq m \leq r-2}d_{mr}\Pi((\p^{m+1}v)\symdot \p^{r-m}h)\label{bdy int s=r},\\
\Pi (D_t\p^sD_t^k h+(\p_ih)(\p^sD_t^kv^i)) = \Pi\p^sD_t^{k+1}h+\Pi(\p_ih)(\p^sD_t^kv^i)+\sum_{0\leq m\leq s-1}d_{mr}\Pi((\p^{m+1}v)\symdot\p^{s-m}D_t^kh), \label{bdy int s<r}
\end{align}
for $k=0$ and $k>0$, respectively.

In addition, when $s=0$,
\begin{dmath}
\frac{d}{dt}E_{0,r} \leq -\int_{\DD_t}\rho \delta^{ij}(D_t^{r}\p_i h)(D_t^r v_j)\dx+\int_{\DD_t}\rho e'(h)(D_t^{r+1}h)(D_t^rh)\dx + C(M)||{e'(h)}D_t^rh||_{\lei}^2 \label{full time energy 1},
\end{dmath}
where we have used the fact that $|e''(h)|\leq c_0|e'(h)|$. Furthermore, since
\begin{multline}
||D_t^r\p h-\p D_t^{r}h||_{\lei} \lesssim  \sum_{l_1+l_2=r-1}||\p D_t^{l_1}v\symdot\p D_t^{l_2}h||_{\lei}\\
+ \sum_{l_1+\cdots+ l_n= r-n+1, \, n\geq 3}||\p D_t^{l_1} v\cdots \p D_t^{l_{n-1}} v\symdot \p D_t^{l_n}h||_{\lei}
  \lesssim_{K,M}(\sum_{i\leq r-1}||v||_{i,0}+||h||_{i,0})(\sum_{i\leq r}||v||_{i,0}+||h||_{i,0}),
\end{multline}
(\ref{full time energy 1}) becomes, after integrating by parts on the first integral on the RHS of (\ref{full time energy 1}),
\begin{dmath}
\frac{d}{dt}E_{0,r} \leq \int_{\DD_t}\rho \delta^{ij}( D_t^{r} h)(D_t^r \di v)\dx+\int_{\DD_t}\rho e'(h)(D_t^{r+1}h)(D_t^rh)\dx + C(M)||{e'(h)}D_t^rh||_{\lei}^2+C(M)\sum_{i\leq r}(||h||_{i,0}+||v||_{i,0})^2  \label{full time energy 2}.
\end{dmath}
But since
$$
D_t^{r}\di v = - D_t^{r+1} e(h) = -e'(h)D_t^{r+1}h - g_r,
$$
and because $||e'(h)D_t^rh||_{\lei}\leq c_0||\sqrt{e'(h)}D_t^rh||_{\lei}$, which is part of $||h||_r$, (\ref{full time energy 2}) becomes
\begin{align}
\frac{d}{dt}E_{0,r}\leq C(M)\sum_{i\leq r}(||h||_{i}+||v||_{i,0})^2.
\end{align}

Furthermore, let $K_r$ be defined as \eqref{K_r}, we have
\begin{align}
\frac{d}{dt}K_r = 2\int_{\DD_t}\rho |\p^{r-1}\curl v|\cdot|D_t\p^{r-1} \curl v|\dx.
\end{align}
But since the curl satisfies the equation
 $$D_t\curl_{ij}v=-(\p_iv^k)(\curl_{kj}v)+(\p_jv^k)(\curl_{ki}v),$$
then
\begin{dmath}
|D_t\p^{r-1}\curl v| \leq |\p^{r-1} D_t\curl v|+\sum_{0\leq m\leq r-2}e_{mr}(\p^{m+1}v)\symdot \p^{r-1-m}\,\curl v\\
 \lesssim \sum_{0\leq m\leq r-1}e_{mr}(\p^{m+1}v)\symdot \p^{r-1-m}\,\curl v.
\end{dmath}
The term $||(\p^{m+1}v)\symdot \p^{r-1-m}\,\curl v||_{\lei}$ can be bounded by
\begin{align}
|\p v|_{L^{\infty}}\sum_{j\leq r-1}||\p^{j}\curl v||_{\lei}+|\curl v|_{L^{\infty}}\sum_{j\leq r-1}||\p^{j+1}v||_{\lei}.
\end{align}

On the other hand,
\begin{dmath}
\sum_{j\leq r+1}\frac{dW_{j}^2}{dt} \lesssim \sum_{j\leq r+1}(W_j^2+W_j(||f_j||_{\lei}+||g_j||_{\lei}) \\
\lesssim E_r^*+\sum_{j\leq r}(||f_r||_{\lei}^2+||g_r||_{\lei}^2).
\end{dmath}
The first inequality comes from the energy estimates for the wave equation, e.g., Theorem \ref{energy est for wave eq}.

Summing these up, we have proved:
\thm Let $E_r$ be defined as (\ref{Er}), for all $r\geq 1$ we have
\begin{multline}
|\frac{dE_r}{dt}|\lesssim_{K,M} E_r^*+\sum_{k+s=r,k,s>0}\bigg(||\Pi\p^s D_t^kh||_{\leb} \Big(||\Pi\p^s D_t^{k+1}h||_{\leb}\\
+||\Pi(\p_ih)(\p^sD_t^kv^i)||_{\leb}+\sum_{0\leq m\leq s-1}||\Pi((\p^{m+1}v)\symdot\p^{s-m}D_t^kh)||_{\leb}\Big)\bigg)\\
+||\Pi\p^r h||_{\leb}\Big(||\Pi \p^r D_t h||_{\leb}+\sum_{0\leq m \leq r-2}||\Pi((\p^{m+1}v)\symdot \p^{r-m}h)||_{\leb}\Big)\\
+C(K,M)(\sum_{i\leq r-1}||v||_{i,0}+||h||_{i,0})(\sum_{i\leq r}||v||_{i,0}+||h||_{i})^2 +\sum_{j\leq r}(||f_r||_{\lei}^2+||g_r||_{\lei}^2)
\end{multline}

\mydef(Mixed boundary Sobolev norm)
let $u(t,\cdot):\RR^n \to \RR$ be a smooth function. We define
\begin{align*}
\lee u\ree_r = \sum_{k+s=r}||\nab^sD_t^ku||_{\llb}.
\end{align*}
\indent Now, let us get back to Lagrangian coordinate. Based on the computation we have as well as \eqref{gag-ni}, controlling $\frac{d}{dt}$ requires to bound
$$
||v||_{r,0},||h||_{r},\sum_{j\leq r-1}||\nab^jv||_{\llb}, \lee h\ree_r,
$$
and
$$\sum_{k+s=r,s\geq2}||\Pi\nab^s D_t^{k+1}h||_{\llb}.
$$
\thm \label{interior and boundary estimates thm}
With the a priori assumptions (\ref{geometry_bound})-(\ref{D_te(h)}) hold, there are continuous functions $C_r$ such that,
\begin{align}
||v||_{r,0}^2+||h||_r^2 \leq C_r(K,M,c_0, E_{r-1}^*)E_r^*.\label{interior estimates}
\end{align}
In addition to that,
\begin{align}
||D_th||_{r,1}^2+\lee h\ree _r^2 \leq C_r(K,M,c_0,\frac{1}{\epsilon}, E_{r-1}^*){E_r^*}, \label{boundary estimates}
\end{align}
where $||D_th||_{r,1}$ is given in Definition \ref{improved mixed norms}.
\subsection{Interior estimates, bounds for $||v||_{r,0}$,$||h||_{r}$}\label{section 5.2}
Our strategy is to first apply Theorem \ref{hodge} to control $||v||_{r,0}$ in terms of the energies $E_r$ and $L^2$ norm of $h$, and then we will apply our elliptic estimate (\ref{ell est II}) to control $||h||_r$.  Now, since 
\begin{equation}
||v||_{r,0} \leq ||\nab^r v||_{\lli} + \sum_{k+s=r,0< k< r}||\nab^sD_t^kv||_{\lli},
\end{equation}
and \eqref{D_t^k v} yields 
\begin{align}
\nab^s D_t^k v = -\nab^{s+1}D_t^{k-1}h+ c_{\alpha\beta\gamma}(\p^{\alpha_1} v)\cdots(\p^{\alpha_m} v)(\p^{\beta_1}D_t^{\gamma_1}h)\cdots(\p^{\beta_n}D_t^{\gamma_n}h),
\end{align}
where \footnotemark \footnotetext{The second term on the right drops when $k=1$.}
$$
\alpha_1+\cdots+\alpha_m+(\beta_1+\gamma_1)+\cdots+(\beta_n+\gamma_n) = r,
$$
$$
1\leq \alpha_i \leq r-1, \q 1\leq \beta_j +\gamma_j\leq r-1.
$$
This implies that \footnotemark \footnotetext{We remark here that we have proved in \cite{LL} that if $r\leq 4$, then $$\sum_{k+s=r,0< k< r}||\nab^sD_t^kv||_{\lli} \leq \sum_{k+s=r,0< k< r}||\nab^{s+1}D_t^{k-1}h||_{\lli}+C(K,M)(\sum_{j\leq r-1}||\nab^j v||_{\lli}+ \sum_{j\leq r-1}||h||_{j,0}).$$ }
\begin{dmath}
\sum_{k+s=r,0< k< r}||\nab^sD_t^kv||_{\lli} \leq \sum_{k+s=r,0< k< r}||\nab^{s+1}D_t^{k-1}h||_{\lli}
\\+C_r(K,M,\sum_{j\leq r-2}||\nab^j v||_{\lli}, \sum_{j\leq r-2}||h||_{j,0})
(\sum_{j\leq r-1}||\nab^j v||_{\lli}+ \sum_{j\leq r-1}||h||_{j,0}).\label{nab^s D_t^k v}
\end{dmath}
So the terms of order $r$ except for $||\nab^rv||$ can be combined with $||h||_r$, up to lower order terms. Now, Theorem \ref{hodge} yields,
\begin{align}
||\nab^r v||_{\lli}\lesssim \sqrt{E_r}+||\nab^{r-1} \di v||_{\lli}.
\end{align}
We recall that $\di v = -e'(h)D_th$, hence
\begin{align}
||\nab^r v||_{\lli}\lesssim_{M,c_0} \sqrt{E_r}+ \sum_{j\leq r}||h||_{j,0},
\end{align}
via interpolation.
Therefore, 
\begin{dmath}
||v||_{r,0} \lesssim_{K,M,c_0}\sqrt{E_r}+ \sum_{j\leq r}||h||_{j,0}\\
+C_r(\sum_{j\leq r-2}||\nab^j v||_{\lli}, \sum_{j\leq r-2}||h||_{j,0})\cdot
(\sum_{j\leq r-1}||\nab^j v||_{\lli}+ \sum_{j\leq r-1}||h||_{j,0}).\label{v_4}
\end{dmath}

To bound $||h||_r$, since (\ref{ell est II}) provides, for each $k,s$ that $k+s=r$,
\begin{align}
||\nab^sD_t^kh||_{\lli}\lesssim_{K,M}||\Pi\nab^sD_t^kh||_{\llb}+\sum_{0\leq j\leq s-2}||\nab^{j}\lap D_t^kh||_{\lli} + ||\nab D_t^k h||_{\lli},
\end{align}
for $s\geq 2$. The term $||\Pi\nab^sD_t^kh||_{\llb}$ bounded by $(||\nab h||_{L^{\infty}(\p\Omega)}E_r)^{\frac{1}{2}}$, by the construction of $E_r$, whereas $||\nab D_t^kh||_{\lli}$ is part of $\sum_{j=1}^{r}W_j$ since $k< r$. Further, by the wave equation (\ref{WWr}),
\begin{align}
\sum_{\substack{0\leq j\leq s-2\\
 2\leq s\leq r\\
 s+k=r}}||\nab^{j}\lap D_t^kh||_{\lli}\leq \sum_{\substack{0\leq j\leq s-2\\
 2\leq s\leq r\\
 s+k=r}} (||\nab^j D_t^{k+2}e(h)||_{\lli} + ||\nab^j f_{k+1}||_{\lli}+||\nab^j g_{k+1}||_{\lli})\label{nab^j lap D_t^k h}.
\end{align}
But since 
\begin{align}
\nab^{s-2}f_{k+1} = \sum_{\substack{\alpha_1+\cdots+\alpha_m+(\beta_1+\gamma_1)+\cdots+(\beta_n+\gamma_n) = r,\\
1\leq \alpha_i \leq r-1, \q 1\leq \beta_j \leq r-1}} c_{\alpha\beta\gamma}(\p^{\alpha_1} v)\cdots(\p^{\alpha_m} v)(\p^{\beta_1}D_t^{\gamma_1}h)\cdots(\p^{\beta_n}D_t^{\gamma_n}h),
\end{align}
and
\begin{align}
\nab^{s-2}g_{k+1} =\sum_{\substack{(\alpha_1+\beta_1)+\cdots+(\alpha_n+\beta_n) = r\\
1\leq \alpha_i+\beta_i\leq r-1, n\geq 2}}c_{n\alpha\beta\gamma}e^{(n)}(h) (\p^{\alpha_1}D_t^{\beta_1}h)\cdots(\p^{\alpha_n}D_t^{\beta_n}h).
\end{align}
Thus, 
\begin{multline}
\sum_{j\leq s-2, s+k=r}(||\nab^{j}f_{k+1}||_{\lli}+||\nab^{j}g_{k+1}||_{\lli}) \leq \\ C_r(K,M,\sum_{j\leq r-2}||\nab^j v||_{\lli}, \sum_{j\leq r-2}||h||_{j})
(\sum_{j\leq r-1}||\nab^j v||_{\lli}+ \sum_{j\leq r-1}||h||_{j}). \label{nab f and g}
\end{multline}
On the other hand, since $|e^{(l)}(h)|\leq c_0|e'(h)|$, and
\begin{multline}
\nab^{s-2}D_t^{k+2} e(h)
 = e'(h)\nab^{s-2}D_t^{k+2}h\\
 + \sum_{\substack{(\alpha_1+\beta_1)+\cdots+(\alpha_m+\beta_m) = r\\
1\leq \alpha_i+\beta_i\leq r-1, m\geq 2}}c_{m\alpha\beta\gamma}e^{(m)}(h) (\p^{\alpha_1}D_t^{\beta_1}h)\cdots(\p^{\alpha_m}D_t^{\beta_m}h),
\end{multline}
we have
\begin{align}
||\nab^{s-2}D_t^{k+2} e(h)||_{\lli} \leq c_0|e'(h)|\cdot||\nab^{s-2}D_t^{k+2}h||_{\lli}+ C_r(K,M,c_0,\sum_{j\leq r-2}||h||_j)\cdot\sum_{j\leq r-1}||h||_{j}. \label{nab D_t e(h)}
\end{align}

Now, \eqref{nab f and g} and \eqref{nab D_t e(h)} yield
\begin{multline}
\sum_{\substack{0\leq j\leq s-2\\
 2\leq s\leq r\\
 s+k=r}}||\nab^{j}\lap D_t^kh||_{\lli}\lesssim_{K,M,c_0}|e'(h)|\sum_{\substack{0\leq j\leq s-2\\
 2\leq s\leq r\\
 s+k=r}} ||\nab^{j}D_t^{k+2}h||_{\lli}\\
 +C_r(\sum_{j\leq r-2}||\nab^j v||_{\lli}, \sum_{j\leq r-2}||h||_{j})\cdot
(\sum_{j\leq r-1}||\nab^j v||_{\lli}+ \sum_{j\leq r-1}||h||_{j}) \label{interior estimates reduction}
\end{multline}
Furthermore, we apply \eqref{ell est II} again with $q = D_t^{k+2}h$ if $s-2\geq 2$, and then repeat the estimates \eqref{nab^j lap D_t^k h}-\eqref{nab D_t e(h)}, we get
\begin{align}
||h||_r \lesssim_{K,M,c_0}\sqrt{E_r^*}+\sum_{j\leq r}W_j+C_r(\sum_{j\leq r-2}||\nab^j v||_{\lli}, \sum_{j\leq r-2}||h||_{j})\cdot
(\sum_{j\leq r-1}||\nab^j v||_{\lli}+ \sum_{j\leq r-1}||h||_{j}).
\end{align}
But since the last term is of lower order, i.e., it can be bounded by $C_r(K,M,c_0,E_{r-2}^*)\sqrt{E_{r-1}^*}$, and so \eqref{interior estimates} follows. 

\subsection{Boundary estimates, bounds for $\sum_{j\leq r-1}||\nab^j v||_{\llb}, \lee h\ree_r$ and $||\cnab^{r-2}\theta||_{\llb}$}\label{section 5.3}
The control of $\sum_{j\leq r-1}||\nab^j v||_{\llb}$ follows directly form the estimate of $\sum_{j\leq r}||\nab^j v||_{\lli}$ by trace theorem (Theorem \ref{trace theorem}) . On the other hand,  we shall not estimate $\lee h\ree_r$ alone; instead, we estimate\footnotemark \,\,  $||D_th||_{r,1}+\lee h\ree_r$ (Definition \ref{improved mixed norms}) by (\ref{ell est II}). This has to be done since we need to estimate $||f_{r+1}||_{\lei}$ and $||g_{r+1}||_{\lei}$ by $E_r$. 

\footnotetext{The reason that we use the norm $||D_th||_{r,1}$ instead of $||h||_{r+1}$ is because the latter involves $||\nab^{r+1} h||$ which, after applying the elliptic and tensor estimates, gives $||(\cnab^{r-1}\theta)\nab_N h||_{\llb}$ but $||\cnab^{r-1}\theta||_{\llb}$ can only be controlled by $E_{r+1}$. On the other hand, we want to avoid the term $||\nab D_t^rh||_{\lli}$ (this term can not be estimated by the method given in Section \ref{section 6}.1) as well, in order to pass our estimates to the incompressible limit in Section \ref{section 6}.}

We estimate the mixed boundary $L^2$ norm $\lee h\ree_r$ by \eqref{ell est I}, we have
\begin{align}
\lee h \ree_r \lesssim_{K,M,c_0} \sum_{k+s=r}||\Pi\nab^s D_t^k h||_{\llb} + \sum_{\substack{k+s=r\\j\leq s-1}}||\nab^j\lap D_t^k h||_{\lli}+ \sum_{j\leq r-1}||\nab D_t^j h||_{\lli}. 
\end{align}
In addition, for $0< \delta<1$, we have
\begin{dmath}
||D_th||_{r,1} \lesssim_{K,M,c_0} \delta\sum_{\substack{k+s=r\\s\geq 2}}||\Pi\nab^sD_t^{k+1}h||_{\llb}\\
+\delta^{-1}(\sum_{\substack{k+s=r\\s\geq 2, j\leq s-2}}||\nab^j\lap D_t^{k+1}h||_{\lli}+W_{r+1}+\sum_{j\leq r-2}||\nab D_t^{j+1} h||_{\lli}),\label{est D_th}
\end{dmath}
via \eqref{ell est II}.
Moreover, 
\begin{dmath}
 \sum_{\substack{k+s=r\\s\geq 1, j\leq s-1}}||\nab^j\lap D_t^k h||_{\lli}+\sum_{\substack{k+s=r\\s\geq 2, j\leq s-2}}||\nab^j\lap D_t^{k+1}h||_{\lli} \lesssim_{K,M,c_0}\\
  \delta\sum_{\substack{k+s=r\\s\geq 2}}||\Pi\nab^sD_t^{k+1}h||_{\llb}+\delta^{-1}\sum_{j\leq r+1}W_j\\
  +C_r(\sum_{j\leq r-1}||\nab^j v||_{\lli}, \sum_{j\leq r-1}||h||_{j})\cdot
(\sum_{j\leq r}||\nab^j v||_{\lli}+ \sum_{j\leq r}||h||_{j}).
\end{dmath}
This in fact follows from the analysis we had for \eqref{interior estimates reduction}. Therefore, by \eqref{interior estimates}, together with \eqref{est D_th},  and since $\sum_{j\leq r+1}W_j$ is part of $\sqrt{E_r^*}$, we obtain
\begin{dmath}
||D_th||_{r,1}+\lee h\ree_r \lesssim_{K,M,c_0}\delta\sum_{\substack{k+s=r\\s\geq 2}}||\Pi\nab^sD_t^{k+1}h||_{\llb}
+\delta^{-1}C_r(K,M,c_0,E_{r-1}^*)\sqrt{E_r^*}.
\end{dmath}

On the other hand, applying \eqref{tensor est} to $||\Pi\nab^sD_t^{k+1}h||_{\llb}$ with $q=D_t^{k+1}h$, then for $s+k=r$ and $s\geq 2$, we have
\begin{multline}
\delta||\Pi\nab^sD_t^{k+1}h||_{\llb} \lesssim  \delta||(\cnab^{s-2}\theta)\nab_N D_t^{k+1}h||_{\llb}
+\delta\sum_{j\leq s-1}||\nab^j D_t^{k+1}h||_{\llb}\\
+ \delta(||\theta||_{\linf}+\sum_{0\leq l\leq s-2}||\cnab^{l}\theta||_{\llb})(\sum_{0\leq l\leq s-2}||\nab^l D_t^{k+1}h||_{\llb})\\
+ \delta(||\theta||_{\linf}+\sum_{0\leq l\leq s-3}||\cnab^{l}\theta||_{\llb})(\sum_{0\leq l\leq s-1}||\nab^l D_t^{k+1}h||_{\llb}) .
\end{multline}
Now, we assume inductively that \eqref{boundary estimates} holds for lower orders \footnote{In fact, we have proved in \cite{LL} that \eqref{boundary estimates} holds for $r\leq 4$.}, i.e., 
\begin{align}
||D_th||_{r',1}+\lee h\ree_{r'} \leq C_{r'}(K,M,c_0,\frac{1}{\epsilon},E_{r'-1}^*)\sqrt{E_{r'}^*} \label{boundary estimates inductive},
\end{align}
whenever $r'\leq r-1$. 
Then \eqref{theta est} yields that 
\begin{align}
\sum_{2\leq s\leq r}||\cnab^{s-2}\theta||_{\llb} \leq C_r(K,M,c_0,\frac{1}{\epsilon}, E_{r-1}^*)\sqrt{E_r^*} \label{theta bound}.
\end{align}
This, together with \eqref{boundary estimates inductive} implies that
\begin{dmath}
\delta\sum_{s+k=r}||\Pi\nab^sD_t^{k+1}h||_{\llb} \lesssim  \delta\sum_{\substack{s+k=r\\s\geq 2}}||(\cnab^{s-2}\theta)\nab_N D_t^{k+1}h||_{\llb}\\
+\delta C_r(K,M,c_0,\frac{1}{\epsilon}, E_{r-1}^*)\cdot\sum_{j\leq r}\lee h \ree_j
+C_{r}(K,M,c_0,\frac{1}{\epsilon},E_{r-1}^*)\sqrt{E_{r}^*}\label{higer tensor}.
\end{dmath}
Now, since $2\leq s\leq r$, we have
\begin{multline}
\delta\sum_{\substack{s+k=r\\s\geq 2}}||(\cnab^{s-2}\theta)\nab_N D_t^{k+1}h||_{\llb}\lesssim_K
\delta||\theta||_{\linf}||\nab D_t^{r-1}h||_{\llb}+\delta||\nab D_th||_{\linf}||\cnab^{r-2}\theta||_{\llb}\\
+\sum_{\substack{s+k=r\\3\leq s\leq r-1}}\delta||\cnab^{s-2}\theta||_{\llb}^{\frac{1}{2}}||\nab D_t^{k+1}h||_{\llb}^{\frac{1}{2}}||\cnab^{s-2}\theta||_{H^1(\p\Omega)}^{\frac{1}{2}}||\nab D_t^{k+1}h||_{H^1(\p\Omega)}^{\frac{1}{2}} \label{higer tensor in 3 dim},
\end{multline}
via \eqref{gag-ni} when $\Omega\in\RR^3$. Furthermore, when $\Omega\in\RR^2$ we have
\begin{dmath}
\delta\sum_{\substack{s+k=r\\3\leq s\leq r-1}}||(\cnab^{s-2}\theta)\nab_N D_t^{k+1}h||_{\llb} \leq \sum_{\substack{s+k=r\\3\leq s\leq r-1}}\delta||\nab D_t^{k+1}h||_{\linf}||\cnab^{s-2}\theta||_{\llb}\\
\lesssim_K \sum_{\substack{s+k=r\\3\leq s\leq r-1}}\delta||\nab D_t^{k+1}h||_{H^1(\p\Omega)}||\cnab^{s-2}\theta||_{\llb} \label{higher tensor in 2 dim}.
\end{dmath}

Moreover, applying \eqref{theta bound} to \eqref{higer tensor in 3 dim} and \eqref{higher tensor in 2 dim} implies that 
\begin{dmath}
\delta\sum_{\substack{s+k=r\\s\geq 2}}||(\cnab^{s-2}\theta)\nab_N D_t^{k+1}h||_{\llb}\leq  C_r(K,M,c_0,\frac{1}{\epsilon},E_{r-1}^*)\sqrt{E_r^*}\\
+\delta C_r(K,M,c_0,\frac{1}{\epsilon},E_{r-1}^*)\lee h \ree_r.
\end{dmath}
Thus, \eqref{higer tensor} becomes
\begin{equation}
\delta\sum_{s+k=r}||\Pi\nab^s D_t^{k+1}h||_{\llb} \leq C_r(K,M,c_0,\frac{1}{\epsilon}, E_{r-1}^*)(\sqrt{E_r^*}+\delta\lee h\ree_r)\label{higher tensor bound}.
\end{equation}
Therefore, 
\begin{align}
||D_th||_{r,1}+\lee h\ree_r \leq C_r(K,M,c_0,\frac{1}{\epsilon}, E_{r-1}^*)(\sqrt{E_{r}^*}+\delta\lee h\ree_r),
\end{align}
where the term $$\delta C_r(K,M,c_0,\frac{1}{\epsilon}, E_{r-1}^*)\lee h\ree_r $$ can be moved to the LHS when $\delta = \delta(K,M,c_0,\frac{1}{\epsilon}, E_{r-1}^*)$ is chosen sufficiently small, and so \eqref{boundary estimates} follows.
\subsection{Bounds for $\sum_{k+s=r}||\Pi\nab^s D_t^{k+1}h||_{\llb}$ }
This is in fact \eqref{higher tensor bound} with $\delta=1$. But since now \eqref{boundary estimates} has been proved, we obtain
\begin{align}
\sum_{k+s=r}||\Pi\nab^s D_t^{k+1}h||_{\llb} \leq C_r(K,M,c_0,\frac{1}{\epsilon},E_{r-1}^*)\sqrt{E_r^*}.
\end{align}

\subsection{Bounds for $\sum_{j\leq r+1}\frac{dW_{j}^2}{dt}$}\label{section 5.5}
We recall that we have
\begin{align}
\sum_{j\leq r+1}\frac{dW_{j}^2}{dt} \lesssim E_r^* + \sum_{j\leq r+1}W_j(||f_{j}||_{\lli}+||g_{j}||_{\lli}).
\end{align}
Therefore, it suffices to bound $\sum_{j=1}^{r+1}||f_{j}||_{\lli}$ and $\sum_{j=1}^{r+1}||g_{j}||_{\lli}$. However, we have
\begin{align}
||f_{r+1}||_{\lli}+||g_{r+1}||_{\lli} \leq  C_r(K,M,c_0,\frac{1}{\epsilon},E_{r-1}^*)\sqrt{E_r^*},
\end{align}
this is because 
$$
\sum_{j\leq r}(||\nab^j v||_{\lli}+||h||_{j}+||D_th||_{j}) \leq C_r(K,M,c_0,\frac{1}{\epsilon},E_{r-1}^*)\sqrt{E_r^*},
$$
as a consequence of Theorem \ref{interior and boundary estimates thm}.
\subsection{The energy estimates } \label{section 5.6}
We are now ready to prove Proposition \ref{a priori estimate}. Since we have showed that our energies $E_r$ control the interior and boundary Sobolev norms of $v$ and $h$, the only thing left is to control the product of the projected tensors, i.e.,
\begin{align}
\sum_{s+k=r,s>0}\Big(\sum_{0\leq m\leq s-1}\Pi((\nab^{m+1}v)\symdot\nab^{s-m}D_t^kh)\Big),\q \text{for}\,\, k>0 \label{tensor k>0}\\
\sum_{0\leq m \leq r-2}\Pi((\nab^{m+1}v)\symdot \nab^{r-m}h),\q \text{for}\,\, k=0 \label{tensor k=0}
\end{align}
\begin{align}
\sum_{s+k=r,s>0}\Pi((\nab h)\symdot (\nab^sD_t^kv)).\q \text{for}\,\, k>0 \label{tensor special}
\end{align}

We cannot use interpolation (\ref{bdy interpolation}) here since it only applies to tangential derivative $\cnab$. Our strategy is to apply Gagliardo-Nirenberg inequality (i.e., \eqref{gag-ni}) to control terms that involving mixed derivatives\footnote{We want our estimates to be linear in the highest order. One can use Sobolev lemma only to control mixed Sobolev norms as well but the highest order energy would appear quadratically that way.}.  By letting $\alpha=\nab^{s-1}v$ in (\ref{trace}) we get
$$
||\nab^{s-1}v||_{\llb}\lesssim_K \sum_{j\leq s}||\nab^j v||_{\lli}.
$$

Now,  when $\Omega\in\RR^3$, each term of (\ref{tensor k>0}) is bounded as
\begin{itemize}
\item If $m=0$, then
\begin{align}
||\Pi((\nab v)\symdot\nab^{s}D_t^kh)||_{\llb} \leq ||\nab v||_{L^{\infty}} ||\nab^sD_t^k h||_{\llb}.
\end{align}
\item If $m\geq 1$, since $k\geq 1$, we must have $1\leq m\leq r-2$. But if $m=r-2$, then $k=1$ and so $s=r-1$, hence
\begin{align}
||\Pi((\nab^{r-1} v)\symdot\nab D_t h)||_{\llb} \leq ||\nab D_th||_{L^{\infty}}||\nab^{r-1}v||_{\llb}.
\end{align}
Otherwise, since $1\leq m\leq r-3$, we have
\begin{multline}
||\Pi((\nab^{m+1}v)\symdot\nab^{s-m}D_t^kh)||_{\llb}\lesssim_K \\
 ||\nab^{m+1}v||_{\llb}^{1/2}||\nab^{s-m}D_t^kh||_{\llb}^{1/2}||\nab^{m+1}v||_{H^1(\p\Omega)}^{1/2}||\nab^{s-m}D_t^kh||_{H^1(\p\Omega)}^{1/2} \label{product tensor n=3}.
\end{multline}
\end{itemize}
On the other hand, if $\Omega\in\RR^2$, then \eqref{product tensor n=3} can instead be bounded via Sobolev lemma, i.e., 
\begin{align}
||\Pi((\nab^{m+1}v)\symdot\nab^{s-m}D_t^kh)||_{\llb}\lesssim_K ||\nab^{m+2}v||_{\llb}||\nab^{s-m}D_t^kh||_{\llb}.
\end{align}
Therefore, the boundary estimates \eqref{boundary estimates} yields
\begin{align}
\sum_{s+k=r}\sum_{\substack{s>0\\0\leq m\leq s-1}}||\Pi((\nab^{m+1}v)\symdot\nab^{s-m}D_t^kh)||_{\llb}\leq C_r(K,M,c_0,\frac{1}{\epsilon},E_{r-1}^*)\sqrt{E_r^*}.
\end{align}

Similarly, \eqref{tensor k=0} can be bounded by
\begin{itemize}
\item If $m=0$ or $m=r-2$, we have
\begin{align}
||\Pi((\nab v)\symdot \nab^rh||_{\llb} \leq ||\nab v||_{L^{\infty}}||\nab^rh||_{\llb},\\
||\Pi((\nab^{r-1} v)\symdot \nab^2 h||_{\llb} \leq ||\nab^2 h||_{L^{\infty}}||\nab^{r-1}v||_{\llb}.
\end{align}
\item If $1\leq m \leq r-3$, we have
\begin{multline}
||\Pi((\nab^{m+1}v)\symdot\nab^{r-m}h)||_{\llb}\lesssim_K \\
 ||\nab^{m+1}v||_{\llb}^{1/2}||\nab^{r-m}h||_{\llb}^{1/2}||\nab^{m+1}v||_{H^1(\p\Omega)}^{1/2}||\nab^{r-m}h||_{H^1(\p\Omega)}^{1/2}.
\end{multline}
\end{itemize}

 As for (\ref{tensor special}), we recall (e.g., \cite{LL}) that when $r\leq 4$, we have
\begin{align}
\sum_{s+k=r}||\Pi((\nab h)\symdot \nab^sD_t^kv)||_{\llb} \lesssim_{K,M} \lee h\ree_r+ \sum_{j\leq r-1}||\nab^j v||_{\llb}.
\end{align}
However, when $r\geq 5$, since
\begin{align}
\nab^s D_t^k v = -\nab^{s+1} D_t^{k-1}h + c_{\alpha\beta\gamma}(\p^{\alpha_1} v)\cdots(\p^{\alpha_m} v)(\p^{\beta_1}D_t^{\gamma_1}h)\cdots(\p^{\beta_n}D_t^{\gamma_n}h)
\end{align}
where
\begin{align*}
\alpha_1+\cdots+\alpha_m+(\beta_1+\gamma_1)+\cdots+(\beta_n+\gamma_n) =r,\\
1\leq \alpha_i \leq r-1,\q 1\leq \beta_j+\gamma_j\leq r-1.
\end{align*}
Then there must be at most one $\alpha_i\geq r-2$ and further if $\alpha_i = r-1$, the other term must satisfy the a priori assumption \eqref{v,h}. Moreover, there are at most three $i$'s such that $\alpha_i \geq r-3$. Hence,
\begin{multline}
||(\p^{\alpha_1} v)\cdots(\p^{\alpha_m} v)(\p^{\beta_1}D_t^{\gamma_1}h)\cdots(\p^{\beta_n}D_t^{\gamma_n}h)||_{\llb} \leq\\
 C_r(K,M, \sum_{k\leq r-2}||\nab^k v||_{\llb},\sum_{k\leq r-1}\lee h\ree_k)(\sum_{k\leq r-1}||\nab^k v||_{\llb}+\sum_{k\leq r}\lee h\ree_k).
\end{multline}

\prop \label{a priori estimate strong}
Let $r\geq r_0>\frac{n}{2}+\frac{3}{2}$, there is a continuous function $\mathcal{T}_r>0$ such that if
$$
0<T\leq \mathcal{T}_r(c_0, K,\mathcal{E}(0),E_r^*(0)),
$$
where
\begin{align}
\mathcal{E}(t)=|(\nab_N h(t,\cdot))^{-1}|_{\linf}. \label{control ep}
\end{align}
Then any smooth solution of (\ref{EE}) for $0\leq t\leq T$ satisfies
\begin{align}
E_r^*(t)\leq 2E_r^*(0) \label{after gronwall},\\
\mathcal{E}(t)\leq 2\mathcal{E}(0)\label{check E 2},\\
g_{ab}(0,y)Z^aZ^b\lesssim g_{ab}(t,y)Z^aZ^b\lesssim g_{ab}(0,y)Z^aZ^b \label{eqv metric 2},
\end{align}
there exists a $\eta>0$ such that
\begin{align}
|N(x(t,\bar{y}))-N(x(0,\bar{y}))| \lesssim \eta,\q\,\,\bar{y}\in \p\Omega \label{cont N},\\
|x(t,y)-x(0,y)| \lesssim \eta,\q\,\,y\in \Omega \label{cont x},\\
|\frac{\p x(t,\bar{y})}{\p y}-\frac{\p x(0,\bar{y})}{\p y}|\lesssim \eta,\q\,\,\bar{y}\in \p\Omega, \label{cont dx}
\end{align}
hold.
To prove Proposition \ref{a priori estimate strong}, we will be using Sobolev lemmas. But then we must make sure that we can control the Sobolev constants. By Lemma \ref{interior sobolev} and \ref{boundary soboolev}, the Sobolev constants depend on $K = \frac{1}{l_0}$, in fact we are allowed to pick a $K$ depending only on initial conditions, which is proved in \cite{CL}. On the other hand, the change of the Sobolev constants in time are controlled by a bound for the time derivative of the metric in Lagrangian coordinate. We also need to control the constant $\frac{1}{\epsilon}$ appears to be in the physical sign condition (\ref{RT sign}).

\lem \label{check a priori} Assume the conditions in Proposition \ref{a priori estimate strong} hold. Then there are continuous functions $C_{r_0}$ such that
\begin{align}
||\nab v||_{L^{\infty}(\Omega)}+||\nab h||_{L^{\infty}(\Omega)}\leq C_{r_0}(K,c_0,E_0,\cdots,E_{r_0}),\label{check 1}\\
||\nab^2 h||_{L^{\infty}(\Omega)}+||\nab D_t h||_{L^{\infty}(\Omega)}\leq C_{r_0}(K,c_0, E_0,\cdots,E_{r_0}),\label{check 2-deriv-h} \\
||\nab \cdot \curl v||_{L^{\infty}(\Omega)} \leq C_{r_0}(K, E_0,\cdots,E_{r_0})\label{check curl},\\
||e'(h)D_th||_{L^{\infty}(\Omega)}+||e'(h)D_t^2h||_{L^{\infty}(\Omega)} \leq C_{r_0}(K,c_0, E_0,\cdots,E_{r_0}),\label{check D_te(h)}\\
||\theta||_{\linf}\leq C_{r_0}(K,c_0, \mathcal{E},E_0,\cdots,E_{r_0}) ,\label{check theta}\\
|\frac{d}{dt}\mathcal{E}|\leq C_{r_0}(K,\mathcal{E},E_0,\cdots,E_{r_0}). \label{check E}
\end{align}

\begin{proof}
By Sobolev lemmas, we have
$$
||\nab v||_{L^{\infty}(\Omega)}\lesssim_K \sum_{ j\leq 3} ||\nab^j v||_{\lli},
$$
$$
||\nab h||_{L^{\infty}(\Omega)}\lesssim_K \sum_{j\leq 3} ||\nab^j h||_{\lli},
$$
and
$$
||\nab^2 h||_{L^{\infty}(\Omega)}+||\nab D_t h||_{L^{\infty}(\Omega)}\lesssim_K \sum_{j\leq 4} ||\nab^j h||_{\lli}+\sum_{j\leq 3}||\nab^j D_th||_{\lli}.
$$
So, as a consequence of our interior and boundary estimates, (\ref{check 1})-(\ref{check 2-deriv-h}) follows. In addition to these, we have
\begin{align*}
||\nab \cdot \curl v||_{L^{\infty}(\Omega)} \lesssim_K \sum_{j\leq 3}||\nab^j \cdot \curl v||_{\lli},
\end{align*}
and so \eqref{check curl} follows. Now, since the assumptions on $e'(h)$ yield
$$
\sum_{j=1,2}|\nab^j e'(h)|\leq C(M,c_0),
$$
thus
\begin{align*}
||e'(h)D_th||_{L^{\infty}(\Omega)}+||e'(h)D_t^2h||_{L^{\infty}(\Omega)} \lesssim_{K,M,c_0} \sum_{j=1,2}(W_j+||\nab^j D_th||_{\lli}+||\nab^j D_t^2h||_{\lli}),
\end{align*}
and so \eqref{check D_te(h)} follows. 
On the other hand, since $$|\nab^2 h|\geq |\Pi \nab^2 h| = |\nab_N h||\theta|\geq \mathcal{E}^{-1}|\theta|, $$ so (\ref{check theta}) follows from (\ref{check 2-deriv-h}). Lastly, (\ref{check E}) is a consequence of
$$
\frac{d}{dt}||(-\nab_N h(t,\cdot))^{-1}||_{\linf}\lesssim ||(-\nab_N h(t,\cdot))^{-1}||_{\linf}^2||\nab_N h_t(t,\cdot)||_{\linf} \label{change of e},
$$
and (\ref{check 2-deriv-h}).
\end{proof}

\subsubsection{Proof of Proposition \ref{a priori estimate strong}}

\indent Since when $r\geq r_0>\frac{n}{2}+\frac{3}{2}$, we have
$$
|\frac{d}{dt}E_r|\leq C_r(c_0, K,\mathcal{E},E_0,\cdots,E_{r_0})E_r^*,
$$
and the RHS is in fact a polynomial of $E_r^*$ with positive coefficients, we get (\ref{after gronwall}) from Lemma \ref{check a priori} and Gronwall's lemma if $\mathcal{T}_r(c_0, K,\mathcal{E}_0,E_r^*(0))>0$ is sufficiently small. (\ref{check E 2}) is a direct consequence of (\ref{check E}). In addition, we get from (\ref{after gronwall}) and Lemma \ref{check a priori} that
 \begin{align}
||\nab v||_{L^{\infty}(\Omega)}+||\nab h||_{L^{\infty}(\Omega)}\leq C(c_0,K,\mathcal{E}(0),E_0(0),\cdots,E_{r_0}(0)), \label{nab v and h at 0}\\
||\nab^2 h||_{L^{\infty}(\Omega)}+||\nab D_th||_{L^{\infty}(\Omega)}+||\theta||_{\linf}\leq C(c_0,K,\mathcal{E}(0),E_0(0),\cdots,E_{r_0}(0)). \label{nab^2 h and theta at 0}
\end{align}
It follows from these that, when $0<T\leq
\mathcal{T}_r(c_0,K,\mathcal{E}(0),E_r^*(0))$ with $\mathcal{T}_r$ chosen to be sufficiently small,
\begin{align}
||\nab v(t,\cdot)||_{L^{\infty}(\Omega)}+||\nab h(t,\cdot)||_{L^{\infty}(\Omega)}\lesssim ||\nab v(0,\cdot)||_{L^{\infty}(\Omega)}+||\nab h(0,\cdot)||_{L^{\infty}(\Omega)},\\
||\nab\cdot\curl v(t,\cdot)||_{L^{\infty}(\Omega)} \lesssim ||\nab\cdot\curl v(0,\cdot)||_{L^{\infty}(\Omega)},
\end{align}
and
\begin{multline}
||\nab^2 h(t,\cdot)||_{L^{\infty}(\Omega)}+||\nab D_th(t,\cdot)||_{L^{\infty}(\Omega)}+||\theta(t,\cdot)||_{\linf}\\
\lesssim||\nab^2 h(0,\cdot)||_{L^{\infty}(\Omega)}+||\nab D_th(0,\cdot)||_{L^{\infty}(\Omega)}+||\theta(0,\cdot)||_{\linf},
\end{multline}
where $0<t\leq T$.

On the other hand, we have
\begin{align}
||v(t,\cdot)||_{L^{\infty}(\Omega)}\lesssim ||v(0,\cdot)||_{L^{\infty}(\Omega)}, \label{control v in the interior}\\
||\rho(t,\cdot)||_{L^{\infty}(\Omega)} \lesssim ||\rho(0,\cdot)||_{L^{\infty}(\Omega)}.\label{control h}
\end{align}
In fact, (\ref{control v in the interior}) follows since $D_tv=-\p h-\e_n$ and (\ref{nab v and h at 0}), whereas (\ref{control h}) follows since $|D_t\rho|\leq |\rho\,\di v|$. Now,
  (\ref{eqv metric 2}) follows because $D_tg$ behaves like $\nab v$. Furthermore, (\ref{cont N}) follows from
$$
D_t N_a = -\frac{1}{2}N_a(D_t g^{cd})N_cN_d,
$$
and (\ref{eqv metric 2}). On the other hand, since by the definition of the Lagrangian coordinate, we have
$$
D_t x(t,y) = v(t,x(t,y)),
$$
and so (\ref{cont x}) follows since (\ref{control v in the interior}). Lastly, because
$$
D_t\frac{\p x}{\p y}=\frac{\p v(t,x)}{\p x}\frac{\p x}{\p y},
$$
(\ref{cont dx}) follows since (\ref{check 1}). \\
We close this section by briefly going over the idea which shows that one can choose $K$ depends only on the initial conditions.

\lem \label{l0 l1} Let $0\leq \eta\leq 2$ be a fixed number, define $l_1=l_1(\eta)$ to be the largest number such that
$$
|N(\bar{x}_1)-N(\bar{x}_2)|\leq \eta, \q\text{whenever}\,\, |\bar{x}_1-\bar{x}_2|\leq l_1, \bar{x}_1,\bar{x}_2\in\p\DD_t.
$$
Suppose $|\theta|\leq K$, we recall that $l_0$ is the injective radius defined in Section 1.4, then
\begin{align*}
l_0\geq \min(l_1/2,1/K),\\
l_1\geq \min(2l_0, \eta/K).
\end{align*}
\begin{proof}
See Lemma 3.6 of \cite{CL}
\end{proof}
In fact, Theorem \ref{l0 l1} shows that $l_0$ and $l_1$ are comparable to each other as long as the free surface is regular.
\lem \label{K initial data}
Fix $\eta>0$ sufficiently small, let $\mathcal{T}$ be in Proposition \ref{a priori estimate strong}. Pick $l_1>0$ such that, whenever $|x(0,y_1)-x(0,y_2)|\leq 2l_1$,
\begin{align}
|N(x(0,y_1))-N(x(0,y_2))|\leq \frac{\eta}{2} \label{K initial}.
\end{align}
Then if $t\leq \mathcal{T}$ we have
\begin{align}
|N(x(t,y_1))-N(x(t,y_2))|\leq \eta, \label{K at t}
\end{align}
whenever $|x(t,y_1)-x(t,y_2)|\leq l_1$.
\begin{proof}
We have
\begin{multline}
|N(x(t,y_1))-N(x(t,y_2))|\\
\leq |N(x(t,y_1))-N(x(0,y_1))|+|N(x(0,y_1))-N(x(0,y_2))|+|N(x(0,y_2))-N(x(t,y_2))|,
\end{multline}
and so (\ref{K at t}) follows from (\ref{cont N}) and (\ref{cont x}).
\end{proof}
Theorem \ref{K initial data} allows us to pick
$$
l_1(t)\leq \frac{l_1(0)}{2},
$$
in other words, we have if $\frac{1}{l_1(0)}\leq \frac{K}{2}$, then
$$
\frac{1}{l_1(t)}\leq K.
$$
Therefore, Theorem \ref{l0 l1} yields
$$
\frac{1}{l_0(t)}\leq K.
$$

\section{Incompressible limit}\label{section 6}
We consider the Euler equations depending on the sound speed $\kk$, i.e., 
\begin{align}
\begin{cases}
D_t v_{\kk} = -\p h_{\kk} - \e_n,\\
\di v_{\kk} = -D_te_{\kk}(h). \label{comp ww equation}
\end{cases}
\end{align}
Here, we recall that the sound speed $\kk$ is defined by 
$$
\kk : = p_{\kk}'(\rho)|_{\rho=1}.
$$

The purpose of this section is to provide an uniform energy estimate that can be carried over to the incompressible water wave and prove the convergence of $v_{\kk}, h_{\kk}$ as $\kk\to\infty$. In order to do so, we impose the following conditions on $\rho_{\kk}(h)$:
\begin{align}
\rho_{\kk}(h)\to 1,\q \text{and hence}\,\,e_{\kk}(h)\to 0,\q\text{as}\,\,\kk\to\infty,
\label{e_kk unifrom a}
\end{align}
and for each fixed $r\geq 1$, there exists a constant $c_0$ such that
\begin{align}
|e_{\kk}^{(k)}(h)|\leq c_0,\q \text{and}\,\, |e_{\kk}^{(k)}(h)| \leq c_0|e_{\kk}'(h)|^k \leq c_0|e_{\kk}'(h)|, \q \text{if}\,\, k\leq r+1.
\label{e_kk unform b}
\end{align}
\rmk 
We remark here that the conditions \eqref{e_kk unifrom a}-\eqref{e_kk unform b} are satisfied if the equation of state is of the form $p(\rho)= C_\gamma \kk (\rho^\gamma-1)$.

\prop \label{uniform kk}
Let $\widetilde{E}_{r,\kk}$ be defined as
\begin{align}
\widetilde{E}_{r,\kk} = \sum_{s+k=r}E_{s,k}+K_r+\sum_{1\leq j\leq r+1}\widetilde{W}_j,\label{tilde E_r}
\end{align}
where $E_{s,k}$ and $K_r$ are defined as \eqref{Esk}-\eqref{K_r}, and
$$
\widetilde{W}_j = \frac{1}{2}||e_{\kk}'(h)D_t^jh_{\kk}||_{\lli}+\frac{1}{2}||\sqrt{e_{\kk}'(h)}\nab D_t^{j-1}h_{\kk}||_{\lli}.
$$
 For $r\geq r_0>\frac{n}{2}+\frac{3}{2}$, there is a continuous function $\mathcal{T}_r>0$ such that if
$
0<T\leq \mathcal{T}_r
$,
then any smooth solution of (\ref{EE}) for $0\leq t\leq T$ satisfies
\begin{align}
\widetilde{E}_{r,\kk}^*(t)\leq 2\widetilde{E}_{r,\kk}^*(0),
\label{uniform energy bounds}
\end{align}
provided the physical sign condition $-\nab_{N}h_{\kk}|_{\p\Omega}\geq \epsilon>0$ holds.

Based on the analysis we have in Section \ref{section 5.5}, Proposition \ref{uniform kk} is a direct consequence of the next theorem.
\thm\label{uniform energybounds}
Let $\widetilde{E}_{r,\kk}$ be defined as (\ref{tilde E_r}), then there are continuous functions $C_r$ such that, for each fixed $r$
\begin{align}
|\frac{d\widetilde{E}_{r,\kk}(t)}{dt}|\leq C_r(K,\frac{1}{\epsilon},M,\mathfrak{h}_\Omega,c_0,\widetilde{E}_{r-1,\kk}^*)\widetilde{E}_{r,\kk}^*(t),\q t\in[0,T] \label{uniform kk energy est}
\end{align}
holds for all $\kk$ ($\mathfrak{h}_\Omega$ is defined as \eqref{def of h_Omega}), provided that the assumptions \eqref{e_kk unifrom a}-\eqref{e_kk unform b} and
\begin{align}
|\theta_{\kk}|+\frac{1}{l_0}\leq K,\q\text{on}\,\p\Omega,\label{boot kk 1}\\
-\nab_{N}h_{\kk}\geq \epsilon>0,\q\text{on}\,\p\Omega,\label{REkk}\\
1\leq \rho_{\kk}\leq M,\q\text{in}\,\Omega,\\
|\nab\cdot\curl v_{\kk}|\leq M,\q\text{in}\,\Omega,\\
|e_{\kk}'(h)D_th_{\kk}|+|e_{\kk}'(h)D_t^2h_{\kk}|\leq M,\q\text{in}\,\Omega,\label{boot kk extra}\\
|\nab v_{\kk}|+|\nab h_{\kk}|+|\nab^2 h_{\kk}|+|\nab D_th_{\kk}|\leq M,\q\text{in}\,\Omega. \label{boot kk v,h}
\end{align}

\indent  It is easy to see that under a priori assumptions (\ref{boot kk 1})-(\ref{boot kk v,h}), the estimates for $||f_r||_{\lli}$ and $||g_r||_{\lli}$ (Theorem \ref{improved fr and gr}) stay unchanged. To prove Proposition \ref{uniform kk}, the analysis we had in the Section 5 implies that it suffices to prove an analogous version of Theorem \ref{interior and boundary estimates thm}:
\thm \label{The key lemma}
Let 
$$
||h||_{r,1} = \sum_{k+s=r, k\leq r-2}||\nab^s D_t^kh||_{\lli}+||\sqrt{e'(h)}\nab D_t^{r-1}h||_{\lli}+||e'(h)D_t^rh||_{\lli}.
$$
Then under the a priori assumptions (\ref{boot kk 1})-(\ref{boot kk v,h}), there are continuous functions $C_r$ such that,
\begin{align}
||v||_{r,0}^2+||h||_{r,1}^2 \leq C_r(K,M,c_0, \mathfrak{h}_\Omega, \widetilde{E}_{r-1}^*)\widetilde{E}_r^*,\label{new interior estimates}
\end{align}
where $\mathfrak{h}_\Omega$ is defined as \eqref{def of h_Omega}. In addition to that,
\begin{align}
||D_th||_{r,1}^2+\lee h\ree _r^2 \leq C_r(K,M,c_0,\frac{1}{\epsilon}, \widetilde{E}_{r-1}^*){\widetilde{E}_r^*}.
\label{new boundary estimates}
\end{align}
\begin{proof}
\eqref{new interior estimates} follows from the arguments in Section \ref{section 5.2}, apart from terms of the form $||\nab D_t^kh_{\kk}||_{\lli}$, $0\leq k\leq r-1$. There would be no problem to control $||\nab D_t^kh_{\kk}||_{\lli}$ via $||\lap D_t^kh_\kk||_{\lli}$ if $\Omega$ were bounded. This bound is a consequence of the classical Poincar\'e's inequality, i.e., $||D_t^k h_\kk||_{\lli}\leq C(\vol\Omega)||\nab D_t^k h_\kk||_{\lli}$. In the case when $\Omega$ is unbounded, we can still bound 
$||D_t^k h_\kk||_{\lli}$ via $||\nab D_t^k h_\kk||_{\lli}$, given $D_t^k h_\kk$ decays fast enough at infinity. Also, we need the following Poincar\'e' inequality in unbounded domains:
\lem \label{poincare for unbounded domain}
Let $\Sigma\in \RR^n$ be a strip with boundary $\p \Sigma = \Gamma_1\cup\Gamma_2$, then there exists a constant $C=C(\mathfrak{h}_\Sigma)$ such that
\begin{align}
||u||_{L^2(\Sigma)} \leq C(\mathfrak{h}_\Sigma)||\nab u||_{L^2(\Sigma)},\q \text{for each}\,\, u\in H^1(\Sigma), \q u|_{\Gamma_1}=0, \label{poincare 618}
\end{align}
where $\mathfrak{h}_\Sigma$ is the ``height" of $\Sigma$ in the bounded direction.
\begin{proof}
 Without loss of generality, we assume $ \Gamma_1\subset \{x_n=0\}$ and $\Sigma$ is bounded in the $x_n$-direction. Since $u|_{\Gamma_1}=0$, we have
\begin{align}
|u(x',x_n)| \leq \int_{0}^{x_n} |\p_n u(x',\tau)|\,d\tau,
\end{align}
and so 
\begin{align}
|u(x',x_n)|^2 \leq \mathfrak{h}_\Sigma\int_0^{\mathfrak{h}_\Sigma}|\p_n u(x',\tau)|^2\,d\tau.
\end{align}
Hence, \eqref{poincare 618} follows by integrating this with respect to $x=(x',x_n)$ with $C=\mathfrak{h}_\Sigma^2$.
\end{proof}
Now, since $D_t^k h_{\kk}|_{t=0}\in L^2_w(\Omega)$ with $w(x)=(1+|x|^2)^{\mu}, \mu\geq 2$ whenever $\kk$ is sufficient large (Theorem \ref{thm data}), there exists strip $\Omega_{\bar{\epsilon}}\subset\Omega$, chosen independent of $\kk$ and bounded in $x_n$-direction, such that
\begin{align}
\int_{\Omega-\Omega_{\bar{\epsilon}}}|D_t^kh_{\kk}|^2\,dy\leq \bar{\epsilon},\q \text{for each}\,\, 0\leq k\leq r-1, \label{uniform localized solution}
\end{align}
for all $\kk$ sufficiently large, and this yields
\begin{align}
||D_t^k h_{\kk}||_{\lli} \lesssim ||D_t^k h_{\kk}||_{L^2(\Omega_{\bar{\epsilon}})}.
\end{align}
Next, Lemma \ref{poincare for unbounded domain} implies
\begin{dmath}
||D_t^k h_{\kk}||_{\lli} \leq C(\mathfrak{h}_{\Omega})||\nab D_t^k h_{\kk}||_{L^2(\Omega_{\bar{\epsilon}})}\\
\leq ||\nab D_t^k h_{\kk}||_{\lli},
 \text{where}\,\, \mathfrak{h}_\Omega \,\,\text{is the height of}\,\, \Omega_{\bar{\epsilon}}, \label{def of h_Omega}
\end{dmath}
and hence, we have
\begin{dmath}
||\nab D_t^k h_{\kk}||_{\lli}^2 \leq ||D_t^kh_{\kk}||_{\lli}||\lap D_t^k h_{\kk}||_{\lli}
\leq \\
C(\mathfrak{h}_{\Omega})||\nab D_t^kh_{\kk}||_{\lli}(||e_{\kk}'(h)D_t^{k+2}h_{\kk}||_{\lli}+||f_{k+1}||_{\lli}+||g_{k+1}||_{\lli}). \label{6.17}
\end{dmath}
Therefore
\begin{align}
\sum_{0\leq k\leq r-1}||\nab D_t^k h_{\kk}||_{\lli} 
\leq C(\mathfrak{h}_{\Omega})(\sum_{1\leq j\leq r+1}\widetilde{W}_j +\sum_{0\leq j \leq r}\big(||f_j||_{\lli}+||g_j||_{\lli})\big).
\end{align}
The only thing that we have to check at this point is that the estimates for $||f_j||_{\lli}$ and $||g_j||_{\lli}$, $2\leq j\leq r$ does not rely on $||\nab D_t^{j-1}h||_{\lli}$, but this is just Theorem \ref{improved fr and gr}. In addition, \eqref{new boundary estimates} is just \eqref{boundary estimates} and it is proved in Section \ref{section 5.3}. This concludes the proof for Theorem \ref{The key lemma}, and hence Theorem \ref{uniform energybounds}. Proposition \ref{uniform kk} then follows from repeating the arguments in Section \ref{section 5.6}.
\end{proof}

\subsection{Passing $(v_\kk, h_\kk)$ to the limit}
The uniform energy bound $\widetilde{E}_{r,\kk}^*(t)\leq 2\widetilde{E}_{r,\kk}^*(0)$ allows us to prove that the solution of the compressible Euler equations \eqref{comp ww equation} converges as $\kk\to\infty$. To be more precise, we prove:
\thm \label{thm convergen}
 Let $u_0$ be a divergence free vector field such that its corresponding pressure $p_0$, defined by $\lap p_0 = -(\p_i u_0^k)(\p_k u_0^i)$ and $p_0\big|_{\p\DD_0}=0$, satisfies the physical condition $-\nab_N p_0\big|_{\p\DD_0} \geq \epsilon >0$. Let $(u, p)$ be the solution of the incompressible free boundary Euler equations with data $u_0$, i.e.
$$
\rho_0 D_t u = -\p p,\q \di u=0,\qquad p|_{\p\DD_0}=0,\quad u|_{t=0}=u_0
$$
with the constant density $\rho_0 =1$. Furthermore, let $(v_{\kk}, h_{\kk})$ be the solution for the compressible Euler equations \eqref{comp ww equation}, with the density function $\rho_{\kk} :h\rightarrow \rho_{\kk}(h)$, and the initial data $v_{0\kk}$ and $ h_{0\kk}$, satisfying the compatibility condition up to order $r+1$, as well as the physical sign condition \eqref{REkk}. Suppose that $\rho_{\kk} \to \rho_0=1$, $v_{0\kk}\to u_0$ and $h_{0\kk}\to p_0$ as $\kk\to\infty$, such that $\widetilde{E}_{r,\kk}^*(0), r\geq 4$  is bounded uniformly independent of $\kk$, then
$$
(v_{\kk},h_{\kk})\to (u,p)\q \text{in}\,\,C^{r-2}([0,T],\Omega),
$$
where 
$C^l([0,T],\Omega)$
consists all functions $u(t,x)$ with
$\nab^s D_t^k u(t,\cdot),\,\,s+k\leq l$
continuous in $\Omega$.
\begin{proof}
Theorem \ref{The key lemma} together with Sobolev lemma would imply 
\begin{dmath}
\sum_{s+k=r-2}||\nab^s D_t^k v_{\kk}||_{L^{\infty}(\Omega)}+ \sum_{s+k=r-2}||\nab^s D_t^kh_{\kk}||_{L^{\infty}(\Omega)} \lesssim_K \\
\sum_{s+k\leq r}||\nab^s D_t^k v_{\kk}||_{\lli}+ \sum_{s+k\leq r}||\nab^s D_t^kh_{\kk}||_{\lli} \leq 2\widetilde{E}_{r,\kk}^*(0).
\end{dmath}
Furthmore, this implies that when $s+k=r-2$, we have
\begin{align}
\nab^s D_t^k v_{\kk}, \nab^s D_t^kh_{\kk} \in C^{0,\frac{1}{2}}(\Omega),
\end{align}
where $C^{0,\frac{1}{2}}(\Omega)$ is the H\"{o}lder space.  Now, Arzela-Ascoli theorem shows that the solution $(v_{\kk},h_{\kk})$ is uniformly bounded and equicontinuous in $C^{r-2}([0,T]\times\Omega)$. Therefore, $(v_{\kk},h_{\kk})$ converges in $C^{r-2}$ with $r\geq 4$, after possibly passing to a subsequence. Finally, $(v_\kk, h_\kk)\to (u,p)$ follows from $D_t^2e_\kk(h)\to 0$, which is a consequence of $||h_\kk||_{C^2([0,T],\Omega)}$ is bounded independent of $\kk$.
\end{proof}

\section{The initial data and the physical condition}\label{section 7}
The purpose of this section is to prove that there exist initial data satisfying the compatibility condition in some weighted Sobolev spaces that satisfies the conditions in Theorem \ref{thm convergen}. In \cite{LL}, we proved that given a smooth initial domain $\DD_0$ with $\vol\DD_0 <\infty$, and let $u_0$ be a divergence free vector field such that its corresponding pressure $p_0$, defined by $-\lap p_0 = (\p_i u_0^k)(\p_k u_0^i)$ and $p_0|_{\p\DD_0}=0$, satisfies the physical condition $-\nab_N p_0 |_{\p\DD_0} \geq \epsilon>0$. We are able to construct a sequence of initial data $(v_0, h_0)=(v_{0\kk}, h_{0\kk})$, satisfying the compatibility condition up to order $5$, i.e., $D_t^k h|_{\{0\}\times \p\DD_0}=0$ for $0\leq k\leq 5$, such that $(v_{0\kk}, h_{0\kk})\to (u_0,p_0)$ and $E^*_{r,\kk}(0)$ are uniformly bounded for $r\leq 4$.

In the case of a water wave, we also need to prove that the data $(v_{0\kk}, h_{0\kk})$ admits sufficient decay. This can be achieved by constructing data in weighted Sobolev spaces. Our goal is to prove:
\thm \label{thm data}
For each fixed $r\geq 4$, given the initial domain $\DD_0$ is unbounded, diffeomorphic to the half space $\{x\in\RR^n:x_n\leq 0\}$, and any divergence free $u_0 \in H^s_w, s\geq r+1$ (see Definition \ref{weighted sob sp}), where $w=(1+|x|^2)^\mu,\mu\geq 2$, such that its corresponding pressure $p_0$ (defined as above) verifies the physical sign condition. Then there exist data $v_0 = v_{0,\kk}$ and $h_0=h_{0,\kk}$, satisfying the compatibility condition up to order $r+1$, i.e.,
$$
h_k|_{\p\DD_0}=h_{k,\kk}|_{\p\DD_0}=0,\q 0\leq k\leq r+1,
$$
such that the quantities
$$||v_{0,\kk}||_{H_w^{s}(\DD_0)} \q \text{and}\q \sum_{k=0}^r||h_{k,\kk}||_{H_w^{s-k}(\DD_0)},\q s\geq r+1
$$
 are uniformly bounded independent of $\kk$. In addition, we have $v_{0\kk}\to u_0$ and $h_{0\kk}\to p_0$ in $C^1(\DD_0)$.

\subsection{Existence of initial data in weighted Sobolev spaces}

\mydef(The weighted Sobolev spaces) \label{weighted sob sp}

Let $\mathbb{D}\subset\RR^n$ be a domain and $w(x) = (1+|x|^2)^{\mu}$, $\mu\geq 2$. For $p\in[1,\infty)$, we let $L_w^p(\mathbb{D})$ be the Banach space consists of functions $u$ such that
\begin{align*}
||u||_{L^p_w(\mathbb{D})} := (\int_{\mathbb{D}}|u(x)|^p w(x)\dx)^{1/p}<\infty.
\end{align*}
In addition, for any positive integer $s$, we let
$
W_w^{s,p}(\mathbb{D})
$
to be the corresponding weighted Sobolev spaces with the norm
$$
||u||_{W_w^{s,p}(\mathbb{D})} = \sum_{|\alpha|\leq s}||\nab^\alpha u||_{L^p_w(\mathbb{D})}.
$$
Finally, $H^s_w(\mathbb{D}):=W_w^{s,2}(\mathbb{D})$ by convention.

We assume $\DD_0=\Omega$ by a slight abuse of notation.  Theorem \ref{thm data} can be achieved by solving: \footnotetext{The system \eqref{full sys} is generated by re-writing the wave equations \eqref{WWr} as Laplace equations. We refer \cite{LL} for the detailed construction. Also, in \cite{LL}, we only considered the case when $e_{\kk}(h)=\kk^{-1}h$, and we shall solve the general case here as well.}
\begin{align}
\begin{cases}
v_{0} = u_0 +\p\phi,\q \text{in}\,\,\Omega,\\
\lap \phi = -e_{\kk}'(h_0)h_1,\q \text{in}\,\,\Omega,\,\, \text{and}\,\, \phi|_{\p\Omega}=0,\\
\lap h_k = e_{\kk}'(h_0) h_{k+2}+ F_k+G_k,\q \text{in}\,\,\Omega,\,\,\text{and},\,\, 0\leq k\leq r-1,\\
\label{full sys}
h_r=h_{r+1}=0,\q \text{in}\,\,\Omega.
\end{cases}
\end{align}
 Here, $F_k:= f_{k+1}|_{t=0}$ and $G_k=g_{k+1}|_{t=0}$, and hence
\begin{align}
F_k =c_{\alpha\beta}^{\gamma,k}(\p^{\alpha_1} v_0)\cdots(\p^{\alpha_m} v_0)(\p^{\beta_1}h_{\gamma_1})\cdots(\p^{\beta_n}h_{\gamma_n}) \label{Fk},
\end{align}
where 
\begin{enumerate}
\item[I.]$\alpha_1+\cdots+\alpha_m+(\beta_1+\gamma_1)+\cdots+(\beta_n+\gamma_n) = k+2$
\item[II.] $ 1\leq \alpha_1\leq\cdots\leq \alpha_m\leq k+1$
\item[III.]$1\leq \beta_1+\gamma_1\leq\cdots\leq \beta_n+\gamma_n \leq k+1$ and $\beta_j\geq 1$, $\gamma_j\leq k-1$ for all $j$
\end{enumerate}
and
\begin{align}
G_k = c^{\gamma_1\cdots\gamma_m,k}e_\kk^{(m)}(h_0)h_{\gamma_1}\cdots h_{\gamma_m},\label{Gk}
\end{align}
where 
\begin{enumerate}
\item[I.]$\gamma_1+\cdots+\gamma_m = k+2$
\item[II.] $ 1\leq \gamma_1\leq\cdots\leq \gamma_m\leq k+1$
\end{enumerate}
We show the existence of solution for \eqref{full sys} via successive approximation starting from the solution $(h_0^0,h_1^0,\cdots, h_{r-1}^0)$ that solves
\begin{align}
\lap h_k^0 = F_k(\p^{\alpha} u_0, \p^{\beta_0}h_0^0,\cdots,\p^{\beta_{k-1}}h_{k-1}^0),\q 0\leq k\leq r-1 \label{full-system nu=0}
\end{align}
and we define $(h_0^{\nu},\cdots,h_{r-1}^{\nu})$ inductively by solving
\begin{align}
\begin{cases}
v_{0}^\nu = u_0 +\p\phi^\nu,\\
\lap \phi^\nu = -e_\kk'(h_0^{\nu-1})h_1^{\nu-1},\\
\lap h_k^\nu = e_\kk'(h_0^{\nu-1}) h_{k+2}^{\nu-1}+ F_k^{\nu}+G_k^{\nu-1},\q 0\leq k\leq r-3\\
\lap h_k^{\nu} = F_k^{\nu}+G_k^{\nu-1},\q k=r-2,r-1\\
\label{full sys nu}
\phi^\nu|_{\p\Omega}=h_k^{\nu}|_{\p\Omega}=0.
\end{cases}
\end{align}
Here, 
\begin{align*}
F_k^{\nu} = F_k(\p^{\alpha} v_0^\nu, \p^{\beta_0}h_0^\nu,\cdots,\p^{\beta_{k-1}}h_{k-1}^\nu),\\
G_k^{\nu-1} = G_k(h_0^{\nu-1},\cdots,h_{k+1}^{\nu-1}).
\end{align*}
Now, we define that for $0\leq k\leq r-1$,
\begin{align*}
m_k^{\nu} := ||h_k^{\nu}||_{H_w^{s-k}(\Omega)},\q s\geq r+1\\
m_*^{\nu} := \sum_{k\leq r-1} m_k^{\nu}+||v_0^\nu||_{H_w^s(\Omega)}.
\end{align*}
According to the elliptic estimate (Theorem \ref{improve ell est thm}), we have
\begin{align}
||\nab \phi^\nu||_{H_w^s(\Omega)}\leq C||e_\kk'(h_0^{\nu-1})h_1^{\nu-1}||_{H^{s-1}_w(\Omega)},\label{bounds for phi}\\
||h_k^\nu||_{H^{s-k}_w(\Omega)}\leq C(||e_\kk'(h_0^{\nu-1}) h_{k+2}^{\nu-1}||_{H^{s-k-2}_w(\Omega)}+||F_k^\nu||_{H^{s-k-2}_w(\Omega)}+||G_k^{\nu-1}||_{H_w^{s-k-2}(\Omega)}),\q 0\leq k\leq r-3. \label{bounds for h_k^nu}\\
||h_k^\nu||_{H^{s-k}_w(\Omega)}\leq C(||F_k^\nu||_{H^{s-k-2}_w(\Omega)}+||G_k^{\nu-1}||_{H_w^{s-k-2}(\Omega)}),\q k=r-2,r-1.
\end{align}
On the other hand, the Sobolev inequalities are still valid for weighted Sobolev spaces; in other words, we have
\lem (Weighted Sobolev inequalities) \label{sob lemma ww}
Let $w(x)=(1+|x|^2)^\mu$, and let $\Omega$ be a domain with $C^1$ boundary, then
\begin{enumerate}
\item[(a)]$ ||u||_{L_w^{np/(n-sp)}(\Omega)} \leq C||u||_{W^{s,p}_w(\Omega)}$, if $sp<n$.
\item[(b)]$||u||_{L^{\infty}(\Omega)} \leq C||u||_{W^{s,p}_w(\Omega)}$, if $sp>n$.
\end{enumerate}
\begin{proof}
Part (a) follows from the proof given by Evans \cite{Ev} with a slight modification. Part (b) is a direct consequence of the standard Sobolev inequality.
\end{proof}
\rmk The above lemma can be generalized to a much larger class of weighted Sobolev spaces. We refer Turesson \cite{T} Chapter 3 for the details.

Lemma \ref{sob lemma ww} allows us to get the bounds for $||F_k^{\nu}||_{H_w^{s-k-2}}$ and $||G_k^{\nu-1}||_{H_w^{s-k-2}}$. 
\subsubsection{Bounds for $||F_k^{\nu}||_{H_w^{s-k-2}}$}
Since $F_k^{\nu}$ is a sum of products of the form \eqref{Fk},  we have
\begin{itemize}
\item If the product involves less than $4$ terms, i.e., $m+n\leq 3$, then
\begin{dmath}
||(\p^{\alpha_1} v_0^\nu)\cdots(\p^{\alpha_m} v_0^\nu)(\p^{\beta_1}h_{\gamma_1}^{\nu})\cdots(\p^{\beta_n}h_{\gamma_n}^{\nu})||_{H_w^{s-k-2}}\\ \leq C||\p^{\alpha_1}v_0^\nu||_{H_w^{s-k-1}}\cdots ||\p^{\alpha_m}v_0^\nu||_{H_w^{s-k-1}}||\p^{\beta_1}h_{\gamma_1}^{\nu}||_{H_w^{s-k-1}}\cdots||\p^{\beta_n}h_{\gamma_n}^{\nu}||_{H_w^{s-k-1}}\\
\leq p(||v_0^\nu||_{H_w^{s}}, m_0^{\nu},\cdots, m_{k-1}^{\nu}), \label{aa}
\end{dmath}
for some polynomial $p$, where the last inequality is because $\beta_j\leq k+1-\gamma_j$ and $\gamma_j\leq k-1$.
\item If the product involves at least $4$ terms, i.e., $m+n\geq 4$. Then we must have $1\leq\alpha_i\leq \alpha_m\leq k-1$ and $1\leq \beta_j+\gamma_j\leq \beta_n+\gamma_n\leq k-1$. But since $\beta_j\geq 1$, we have
\begin{dmath}
||(\p^{\alpha_1} v_0^\nu)\cdots(\p^{\alpha_m} v_0^\nu)(\p^{\beta_1}h_{\gamma_1}^{\nu})\cdots(\p^{\beta_n}h_{\gamma_n}^{\nu})||_{H_w^{s-k-2}}\\ \leq C||\p^{\alpha_1}v_0^\nu||_{H_w^{s-k}}\cdots ||\p^{\alpha_m}v_0^\nu||_{H_w^{s-k}}||\p^{\beta_1}h_{\gamma_1}^{\nu}||_{H_w^{s-k}}\cdots||\p^{\beta_n}h_{\gamma_n}^{\nu}||_{H_w^{s-k}}\\
\leq p(||v_0^{\nu}||_{H_w^{s}}, m_0^{\nu},\cdots, m_{k-2}^{\nu}). 
\end{dmath}
\end{itemize}
\subsubsection{Bounds for $||G_k^{\nu-1}||_{H_w^{s-k-2}}$ and $||e_\kk'(h_0^{\nu-1})h_{k+2}^{\nu-1}||_{H_w^{s-k-2}}$}
Since $|e^{(m)}(h_0)|\leq c|e_\kk'(h_0)|^m$, \eqref{Gk} together with the weighted Sobolev inequalities imply
\begin{align}
||G_k^{\nu-1}||_{H_w^{s-k-2}} \leq q(e_\kk'(h_0^{\nu-1})m_0^{\nu-1}, \cdots, e_\kk'(h_0^{\nu-1})m_{k+1}^{\nu-1}),\q 0\leq k\leq r-2 \label{7.11}
\end{align}
and since $h_r=0$, 
\begin{align}
||G_{r-1}^{\nu-1}||_{H_w^{s-r-1}} \leq q(e_\kk'(h_0^{\nu-1})m_0^{\nu-1}, \cdots, e_\kk'(h_0^{\nu-1})m_{r-1}^{\nu-1}),
\end{align}
for some polynomial $q$. On the other hand, we have
\begin{dmath}
||e_\kk'(h_0^{\nu-1})h_{k+2}^{\nu-1}||_{H_w^{s-k-2}}\leq \widetilde{q}(e_\kk'(h_0^{\nu-1})m_0^{\nu-1}, e_\kk'(h_0^{\nu-1})m_{k+2}^{\nu-1}), \q 0\leq k\leq r-3 \label{bb}
\end{dmath}
for some polynomial $\widetilde{q}$.

\subsubsection{Bound for $||v_0^\nu||_{H^s_w}$}
The first equation of \eqref{full sys nu} yields
\begin{align}
||v_0^\nu||_{H_w^s} \leq ||u_0||_{H_w^s}+||\nab \phi^\nu||_{H_w^s}.
\end{align}
But since 
\begin{dmath}
||e_\kk'(h_0^{\nu-1})h_{1}^{\nu-1}||_{H_w^{s-1}}\leq q(e_\kk'(h_0^{\nu-1})m_0^{\nu-1}, e_\kk'(h_0^{\nu-1})m_{1}^{\nu-1}),
\end{dmath}
and so \eqref{bounds for phi} implies
\begin{align}
||v_0^\nu||_{H_w^s} \leq ||u_0||_{H_w^s}+q(e_\kk'(h_0^{\nu-1})m_0^{\nu-1}, e_\kk'(h_0^{\nu-1})m_{1}^{\nu-1}).\label{cc}
\end{align}
\subsubsection{A priori bounds for the full system \eqref{full sys nu}}

We conclude from \eqref{aa}-\eqref{bb} that
\begin{dmath}
m_k^{\nu}\leq Ce'm_{k+2}^{\nu-1}+P(||v_0^\nu||_{H_w^{s}}, m_0^{\nu},\cdots, m_{k-1}^{\nu}, e'm_0^{\nu-1},\cdots, e'm_{k+1}^{\nu-1}),
\end{dmath}
for $0\leq k\leq r-3$ and
\begin{align}
m_{r-2}^{\nu}\leq P(||v_0^\nu||_{H_w^{s}}, m_0^{\nu},\cdots, m_{r-3}^{\nu},e'm_0^{\nu-1},\cdots, e'm_{r-1}^{\nu-1}),\\
m_{r-1}^{\nu}\leq P(||v_0^\nu||_{H_w^{s}}, m_0^{\nu},\cdots, m_{r-2}^{\nu},e'm_0^{\nu-1},\cdots, e'm_{r-1}^{\nu-1}).
\end{align}
Summing these up, we get
\begin{align}
m_*^{\nu} \leq P(e_\kk' m_*^{\nu-1}, ||v_0^\nu||_{H_w^{s}})\leq  Q(e_\kk' m_*^{\nu-1}, ||u_0||_{H_w^{s}})  \label{A priori data},
\end{align}
for some polynomials $P$ and $Q$ via \eqref{cc}. In particular, this implies that $m_*^{\nu}$ is uniformly bounded for all $\nu$ by induction whenever $e_\kk'$ (and hence $\kk^{-1}$) is sufficiently small. Finally, the existence follows from subtracting two successive systems of \eqref{full-system nu=0}-\eqref{full sys nu} and the a priori bound \eqref{A priori data}, which is identical to what is in \cite{LL}. 

On the other hand, since $s\geq r+1\geq 5$ and $\Omega\in \RR^n$, $n=2,3$, we have
\begin{align}
||v_{0,\kk}-u_0||_{C^1(\Omega)} \leq C||v_{0,\kk}-u_0||_{H_w^s(\Omega)} \leq ||\nab \phi_{\kk}||_{H_w^s(\Omega)} \leq q(e_{\kk}'||h_{0,\kk}||_{H_w^s(\Omega)}, e'_{\kk}||h_{1,\kk}||_{H_w^{s-1}(\Omega)}),
\end{align}
where $q(0)=0$.
This implies $v_{0,\kk}\to u_0$ in $C^1(\Omega)$ since $||h_{k,\kk}||_{H_w^{s-k}}$ are uniformly bounded independent of $\kk$ and $e'_{\kk}\to 0$. Similarly, we have
\begin{align}
||h_{0,\kk}-p_0||_{H_w^s} \leq C\Big(||e'_\kk h_{2,\kk}||_{H^{s-2}_w}+ ||e''_\kk h_{1,\kk}^2||_{H_w^{s-2}}+ (||u_0||_{H_w^s}+||\phi_\kk||_{H_w^s})||\phi_\kk||_{H_w^s}\Big).
\end{align} 
The RHS $\to 0$ as $\kk\to \infty$ via \eqref{bounds for phi} and \eqref{7.11}. This concludes the proof of Theorem \ref{thm data}.

\cor 
$\widetilde{E}_{r,\kk}^*(0)$ in Proposition \ref{uniform kk} is uniformly bounded whenever $\kk^{-1}$ is sufficiently small. 
\begin{proof}
Since $w(x)\geq 1$, we have
\begin{align}
\sum_{k+s\leq r}\int_{\Omega}\rho_0Q(\p^sh_{k,\kk},\p^sh_{k,\kk})\dx \lesssim \sum_{ k\leq r}||h_{k,\kk}||_{H^{r-k}(\Omega)}^2 \leq \sum_{k\leq r}||h_{k,\kk}||_{H^{r-k}_w(\Omega)}^2 ,
\end{align}
and by the trace lemma,
\begin{dmath}
\sum_{k+s\leq r}\int_{\p\Omega}\rho_0Q(\p^sh_{k,\kk},\p^sh_{k,\kk})\,dS \lesssim \sum_{k\leq r}||h_{k,\kk}||_{H^{r-k}(\p\Omega)}^2 \lesssim  \sum_{k\leq r}||h_{k,\kk}||_{H_w^{r+1-k}(\Omega)}^2.
\end{dmath}
In addition to these, we have for $r\geq 2$ that
\begin{dmath}
\sum_{k\leq r}||(\p^{r-k}D_t^k v_\kk)\big|_{t=0}||_{L^2} \lesssim ||v_{0,\kk}||_{H_w^r}+P(||v_{0,\kk}||_{H_w^{r-1}},\sum_{k\leq r-1}||h_{k,\kk}||_{H_w^{r-1-k}}),
\end{dmath}
for some polynomial $P$ via \eqref{D_t^k v}. This shows when $s+k\geq 2$,
\begin{align}
\int_{\Omega}\rho_0Q(\p^s D_t^k v_\kk\big|_{t=0},\p^s D_t^kv_\kk\big|_{t=0})\dx
\end{align}
is uniformly bounded. Finally, since $h_r = h_{r+1}=0$ in $\Omega$, we have
\begin{align*}
\sum_{k=1}^{r+1} \widetilde{W}_k(0)
\end{align*}
 is uniformly bounded as well.
 \end{proof}

\subsection{The physical condition}
It is plausible to assume the physical condition \eqref{taylor-sign cond} on the initial data. We are able to show that for a slight compressible (i.e., $\kk^{-1}$ is small), irrotational water wave, the quantity 
$
-\nab_N h_0
$
is pointwisely greater than a positive constant depending only on the geometry of the free surface, as long as the free surface is not self-intersecting. This can be shown via the maximum principle since $h_0$ is superharmonic in the case of a slightly compressible and irrotational water wave. The original version of our proof is given by Wu \cite{W2}.

In particular, Theorem \ref{thm data} together with Lemma \ref{sob lemma ww} yield that for $r=4$, there exists a constant $C$ such that
\begin{align}
\sum_{k=1,2}||h_k||_{L^{\infty}(\DD_0)}\leq \frac{C}{(1+|x|^2)^{\mu}},\q \mu\geq 2 \label{uniform bounds for D_th and D_t^2h}
\end{align}
whenever $\kk^{-1}$ is sufficiently small. In addition, since $\curl v=0$, we have $\p_iv_j = \p_jv_i$ for each $i,j$, and so $h_0$ and $x_n$ satisfies
\begin{align}
-\lap (h_0+x_n) = |\p v_0|^2 - (e'_{\kk}(h_0)h_2+e''_{\kk}(h_0)h_1^2) \label{superharmonic}.
\end{align}
Now, \eqref{uniform bounds for D_th and D_t^2h} guarantees that the right hand side of \eqref{superharmonic} is positive pointwisely whenever $\kk$ is large (and so $e'_{\kk}$ and $e''_{\kk}$ are small); in other words, $h_0+x_n$ is superharmonic in the case of a slightly compressible, irrotational liquid. For any $\psi\in C_c^1(\p\DD_0)$, $\psi\geq 0$, let $\phi$ be the harmonic extension of $\psi$ in $\DD_0$, i.e., $\phi$ solves
\begin{align}
\begin{cases}
\lap \phi = 0,\q \text{in}\,\,\DD_0\\
\phi |_{\p\DD_0} = \psi.
\end{cases}
\end{align}
In fact, it is easy to see that 
\begin{align}
\phi(x) = o(|x|^{2-n}),\q \nab\phi  = o(|x|^{1-n}) \label{little o},
\end{align}
as $|x|\to\infty$. 

Now, applying the Green's identity \footnotemark to $\phi$ and $h_0+x_n$, we get
\begin{align}
\int_{\p\DD_0}(h_0+x_n)\nab_N \phi - \phi \nab_N(h_0+x_n)\,dS = \int_{\DD_0}\phi (|\nab v_0|^2-e'_{\kk}(h_0)h_2-e''_{\kk}(h_0)h_1^2))\dx.
\end{align}
\footnotetext{Green's identity holds here on unbounded domains because of the decay properties and the $L^2$ integrability of our functions involved. }
But since $h_0=0$ on $\p\DD_0$, we have
\begin{align}
\int_{\p\DD_0}-\phi\nab_N h_0\,dS = \int_{\p\DD_0}(\phi\nab_N x_n-x_n\nab_N\phi)\,dS+\int_{\DD_0}\phi (|\nab v_0|^2-e'_{\kk}(h_0)h_2-e''_{\kk}(h_0)h_1^2))\dx. \label{rmk 2}
\end{align}
On the other hand, applying the Green's identity  again to $\phi$ and $x_n$ on the strip region between $\p\DD_0$ and $\{x\in\RR^n:x_n=b\}$ (with the upward unit normal $N_b=e_n$), we get
\begin{dmath}
\int_{\p\DD_0}(\phi\nab_N x_n-x_n\nab_N\phi)\,dS = \int_{x_n=b}(\phi\nab_{N_b} x_n-x_n\nab_{N_b}\phi)\,dS
= \int_{x_n=b} \phi \,dS- b\int_{x_n=b}\nab_{N_b}\phi\,dS = \int_{x_n=b} \phi \,dS.
\end{dmath}
The integral $\int_{x_n=b}\nab_{N_b}\phi\,dS=0$ is a direct consequence of \eqref{little o} and the Gauss-Green's formula when $n\geq 3$. Therefore,
\begin{dmath}
\int_{\p\DD_0}-\phi\nab_N h_0\,dS = \int_{x_n=b} \phi \,dS+\int_{\DD_0}\phi (|\nab v_0|^2-e'_{\kk}(h_0)h_2-e''_{\kk}(h_0)h_1^2))\dx\geq \int_{x_n=b} \phi \,dS.
\end{dmath}
Let $G=G(x,y)$ be the Green's function for the region $\DD_0$, then by Green's representation formula we have
$$
\phi(x) = \int_{\p\DD_0}\psi(y)\nab_N G(x,y)\,dS(y),\q \text{for}\,\, x\in\DD_0.
$$
But this then implies
\begin{dmath}
\int_{\p\DD_0}-\psi(y)\nab_N h_0(y)\,dS(y)\geq \int_{x_n=b}\phi(x)\,dS(x) = \int_{\p\DD_0}\psi(y)\int_{x_n=b}\nab_N G(x,y)\,dS(x)\,dS(y).
\end{dmath}
Since $\psi\in C_c^1(\p\DD_0)$, $\psi\geq 0$ is arbitrary, we must have that for each $y\in\p\DD_t$, 
\begin{align}
-\nab_N h_0(y)\geq \int_{x_n=b}\nab_N G(x,y)\,dS(x). \label{rmk 3}
\end{align}
From the maximum principle, we know that there exists $\epsilon>0$ such that
$
\int_{x_n=b}\nab_N G(x,y)\,dS(x) \geq \epsilon,
$
for every $y\in\p\DD_0$.

Therefore, the following theorem is justified for a slightly compressible, irrotational liquid.
\thm \label{thm 7.2} Assume that at time $0$, the water region $\DD_0\subset\RR^n$, $n\geq 3$ is unbounded, diffeomorphic to $\{x\in\RR^n:x_n\leq 0\}$, whose boundary $\p\DD_0$ satisfies $|\theta|+|1/l_0|\leq K$. Then there exists a positive constant $\epsilon$, depending only on $\p\DD_0$ such that
\begin{align*}
-\nab_N h_0(y)\geq \epsilon>0
\end{align*}
holds for each $y\in\p\DD_0$.
\rmk 
In the original proof given by Wu \cite{W2}, the pressure $p_0$ is automatically superharmonic, since $v_0$ is divergence free implies 
$$
-\lap p_0 = |\nab v_0|^2>0.
$$
But we need to put extra effort to make sure that $h_0$ is superharmonic in the case of a slightly compressible liquid.

\rmk 
The presence of the gravity is essential for proving that $-\nab_N h_0$ is bounded uniformly below by a positive constant. Since otherwise the term $\int_{\p\DD_0}(\phi\nab_N x_n-x_n\nab_N\phi)\,dS$ on the right of \eqref{rmk 2} would be $0$. In this case we can only conclude $-\nab_N h_0\geq 0$.

\section{The weighted a priori estimates for the Euler equations}\label{section 8}
The purpose of this section is to generalize Proposition \ref{main 1} to weighted $L^2$ Sobolev spaces. In the previous section, we have shown that for each fixed $r$, there exist data in $H^{r+1}_w$ that satisfying the compatibility condition, and we shall prove that the corresponding weighted energies for the compressible Euler equations remain bounded within short time.  This will follow from the analysis we have in Section \ref{section 5} given the estimates in Section 3 remain valid in weighed Sobolev spaces; in other words, we need to establish the Christodoulou-Lindblad type elliptic estimates (Proposition \ref{ellpitic estimate}), as well as the tensor estimate (Proposition \ref{tensor estimate}) in the case of weighted spaces. Throughout this section, the weight function $w(x) = (1+|x|^2)^{\mu},\mu\geq 2$.

\subsection{The weighted Christodoulou-Lindblad type elliptic estimates}
We adopt the notations used in Section 3. Let $\Omega$ be a general domain in $\RR^n$ and let $\nab $ be the covariant differentiation with respect to the metric $g_{ij}$ in $\Omega$, and $\cnab$ will refer to the covariant differentiation on $\p\Omega$ with respect to the induced metric $\gamma_{ij}=g_{ij}-N_iN_j$. We will also assume that the normal $N$ to $\p\Omega$ is extended to a vector field of $\Omega$ via the geometric normal coordinate satisfying $g_{ij}N^iN^j\leq 1$ (e.g., Lemma \ref{trace 1}).

\lem \label{Hodge ww} Let $u:\Omega\to\RR^n$ be a vector field and let $\beta_k = \nab_{i_1}\cdots\nab_{i_r}u_k:=\nab_I^r u_k$. If $|\theta|+\frac{1}{l_0}\leq K$, then 
\begin{dmath}
\int_{\Omega}|\nab \beta|^2 w\,d\mu_g \leq 
C(K)\int_{\Omega}(N^iN^jg^{kl}\gamma^{IJ}\nab_k\beta_{Ii}\nab_l\beta_{Jj}+|\di \beta|^2+|\curl \beta|^2+|\beta|^2)w\,d\mu_g \label{CL 5.17}.
\end{dmath}
Here, $\gamma^{IJ} = \gamma^{i_1j_1}\cdots\gamma^{i_rj_r}$.\label{CL 5.18}
\begin{proof}
We follow the proof given in Christodoulou-Lindblad \cite{CL}. Since $g^{ij} = \gamma^{ij}+N^iN^j$, we have 
 $$|\nab\beta|^2 = g^{IJ}g^{kl}\nab_k\beta_I\nab_l\beta_J$$ 
 can be written as a sum of terms of the form (that is, the normal-tangential form)
\begin{align}
N^{i_1}N^{j_1}\cdots N^{i_s}N^{j_s}\gamma^{i_{s+1}j_{s+1}}\cdots \gamma^{i_rj_r}g^{kl}\nab_k\beta_I\nab_l\beta_J,
\end{align}
and if we control the right hand side of \eqref{CL 5.17}, then we have the bounds for integral of \eqref{CL 5.18} for $s=1,2$. However, the following Hodge-type decomposition holds (e.g., \cite{CL}): let $q^{IJ}$ be any product of factors $q^{ij}$ of the form $g^{ij}, \gamma^{ij}$ or $N^iN^j$, then
\begin{dmath}
g^{ij}g^{kl}q^{IJ}\nab_i\beta_{Ik}\nab_j\beta_{Jl}\leq \Big(2(N^iN^jg^{kl}+g^{ij}N^kN^l)+2g^{ik}g^{jl}\\
+(\gamma^{ij}\gamma^{kl}-\gamma^{ik}\gamma^{jl})\Big)q^{IJ}\nab_i\beta_{Ik}\nab_j\beta_{Jl} \label{CL 5.3}.
\end{dmath}
In addition to this, if $R^{ijklIJ}:=(\gamma^{ij}\gamma^{kl}-\gamma^{ik}\gamma^{jl})q^{IJ}$, then
\begin{dmath}
\int_{\Omega}R^{ijklIJ}\nab_k\alpha_{Ii}\nab_j\beta_{Jl}w\,d\mu_g = -\int_{\Omega}(\nab_kR^{ijklIJ})\alpha_{Ii}\nab_j\beta_{Jl}w\,d\mu_g-\int_{\Omega}(R^{ijklIJ})\alpha_{Ii}\nab_j\beta_{Jl}(\nab_k w)\,d\mu_g,\label{CL 5.12}
\end{dmath}
 via integrating by parts, since $(\gamma^{ij}\gamma^{kl}-\gamma^{ik}\gamma^{jl})\nab_j\nab_k\beta=0$ and $N_k(\gamma^{ij}\gamma^{kl}-\gamma^{ik}\gamma^{jl})=0$. 
 
Now, by \eqref{CL 5.3} and \eqref{CL 5.12}, and since the weight satisfies $|\nab w|\leq \frac{Cw}{1+|x|}$, the bounds for integral of \eqref{CL 5.18} for $s=1,2$ gives us the integral of \eqref{CL 5.18} also for $s=0$. This is because
\begin{dmath}
\Big|(\gamma^{ij}\gamma^{kl}-\gamma^{ik}\gamma^{jl})q^{IJ}(\nab_i\beta_{Ik}\nab_j\beta_{Jl}-\nab_k\beta_{Ii}\nab_j\beta_{Jl})\Big| \leq C|\curl \beta|\cdot |\nab \beta|,
\end{dmath}
and
\begin{align}
\Big|g^{ik}g^{jl}q^{IJ}\nab_i\beta_{Ik}\nab_j\beta_{Jl}\Big|\leq C|\di \beta|^2.
\end{align} 
But then we can use \eqref{Hodge} to get \eqref{CL 5.17}.
\end{proof}

\lem Let $\beta$ be defined as in the previous lemma. If $|\theta|+\frac{1}{l_0}\leq K$, then
\begin{align}
||\beta||_{L^2_w(\p\Omega)}^2\leq C(K) \Big(||\nab \beta||_{L^2_w(\Omega)}^2+||\beta||_{L^2_w(\Omega)}^2\Big),\label{prototype 1}\\
||\beta||_{\llbw}^2\leq C||\Pi\beta||_{\llbw}^2+C(K)\Big(||\di\beta||_{\lliw}^2+||\curl\beta||_{\lliw}^2+||\beta||_{\lliw}^2\Big),\label{prototype 2}\\
||\nab\beta||_{\lliw}\leq C||\Pi\nab\beta||_{\llbw}||\beta||_{\llbw}+C(K)\Big(||\di\beta||_{\lliw}^2+||\curl\beta||_{\lliw}^2+||\beta||_{\lliw}^2\Big).\label{prototype 3}
\end{align}
\begin{proof}
Inequality \eqref{prototype 1} is just \eqref{trace ww}. Let $g^{IJ}=g^{i_1j_1}\cdots g^{i_kj_k}$, \eqref{prototype 2} follows by induction from
\begin{dmath*}
\int_{\p\Omega}g^{IJ}g^{ij}\beta_{Ii}\beta_{Jj}w\,d\mu_{\gamma}
=\int_{\Omega}\nab_k(N^kg^{IJ}(N^iN^j+\gamma^{ij})\beta_{Ii}\beta_{Jj}w)\,d\mu_g\\
=\int_{\Omega}(\nab_k N^k)g^{IJ}(N^iN^j+\gamma^{ij})\beta_{Ii}\beta_{Jj}w\,d\mu_g+\int_{\Omega}N^kg^{IJ}(N^iN^j+\gamma^{ij})\beta_{Ii}\beta_{Jj}(\nab_k w)\,d\mu_g\\
+2\int_{\Omega}N^kg^{IJ}(N^iN^j+\gamma^{ij})\beta_{Ii}\nab_k\beta_{Jj}w\,d\mu_g
\end{dmath*}
On the other hand, we have
\begin{dmath*}
2\int_{\Omega}N^kg^{IJ}(N^iN^j+\gamma^{ij})\beta_{Ii}\nab_k\beta_{Jj}w\,d\mu_g
= 2\int_{\Omega}N^kg^{IJ}N^iN^j\beta_{Ii}\nab_k\beta_{Jj}w\,d\mu_g
+2\int_{\Omega}N^kg^{IJ}\gamma^{ij}(\beta_{Ii}\nab_k\beta_{Jj}-\beta_{Ii}\nab_j\beta_{Jk})w\,d\mu_g
+2\int_{\Omega}N^kg^{IJ}\gamma^{ij}\beta_{Ii}\nab_j\beta_{Jk}w\,d\mu_g.
\end{dmath*}
However, 
\begin{dmath*}
2\int_{\Omega}N^kg^{IJ}\gamma^{ij}\beta_{Ii}\nab_j\beta_{Jk}w\,d\mu_g
= -2\int_{\Omega}\nab_j(N^kg^{IJ}\gamma^{ij})\beta_{Ii}\beta_{Jk}w\,d\mu_g
-2\int_{\Omega}N^kg^{IJ}\gamma^{ij}\beta_{Ii}\beta_{Jk}\nab_j w\,d\mu_g
-2\int_{\Omega}N^kg^{IJ}\gamma^{ij}\nab_j\beta_{Ii}\beta_{Jk}w\,d\mu_g,
\end{dmath*}
since $N_j\gamma^{ij}=0$.
Hence,
\begin{dmath}
\int_{\p\Omega}g^{IJ}N^{ij}\beta_{Ii}\beta_{Jj}w\,d\mu_{\gamma}=-\int_{\p\Omega}g^{IJ}\gamma^{ij}\beta_{Ii}\beta_{Jj}w\,d\mu_{\gamma}+2\int_{\Omega}N^kg^{IJ}\gamma^{ij}\beta_{Ii}(\nab_k\beta_{Jj}-\nab_j\beta_{Jk})w\,d\mu_g\\
+ 2\int_{\Omega}N^kg^{IJ}N^iN^j\beta_{Ii}\nab_k\beta_{Jj}w\,d\mu_g-2\int_{\Omega}N^kg^{IJ}\gamma^{ij}\nab_j\beta_{Ii}\beta_{Jk}\,d\mu_g\\
-2\int_{\Omega}\nab_j(N^kg^{IJ}\gamma^{ij})\beta_{Ii}\beta_{Jk}w\,d\mu_g+\int_{\Omega}(\nab_k N^k)g^{IJ}g^{ij}\beta_{Ii}\beta_{Jj}w\,d\mu_g+\int_{\Omega}N^kg^{IJ}g^{ij}\beta_{Ii}\beta_{Jj}(\nab_k w)\,d\mu_g-2\int_{\Omega}N^kg^{IJ}\gamma^{ij}\beta_{Ii}\beta_{Jk}(\nab_j w)\,d\mu_g
\end{dmath}
The last four terms are bounded by $||\beta||_{\lliw}^2$ since $|\nab N|\leq K$ and $|\nab w|\leq Cw/(1+|x|)$, whereas the terms on the first line are contributed to $||\Pi\beta||_{\llbw}^2$ and $||\curl\beta||_{\lliw}^2$. Finally, the terms on the second line are contributed to $||\di \beta||_{\lliw}^2$, and so this finishes proving \eqref{prototype 2}. \eqref{prototype 3} is just \eqref{CL 5.17} after integrating by parts.
\end{proof}

\thm(Christodoulou-Lindblad type elliptic estimates)\label{Christodoulou-Lindblad type elliptic est}
Let $q:\Omega\to\RR$ be a function and suppose $|\theta|+\frac{1}{l_0}\leq K$, we have, for any $r\geq 2$ and $\delta>0$, 
\begin{align}
||\nab^r q||_{\llbw}+||\nab^r q||_{\lliw} \lesssim_K \sum_{s\leq r}||\Pi\nab^s q||_{\llbw}+\sum_{s\leq r-1}||\nab^s\lap q||_{\lliw}+||\nab q||_{\lliw},\label{ell est I ww}\\
||\nab^{r-1} q||_{\llbw}+||\nab^r q||_{\lliw} \lesssim_K \delta\sum_{s\leq r}||\Pi\nab^s q||_{\llbw}+\delta^{-1}\Big(\sum_{s\leq r-2}||\nab^s\lap q||_{\lliw}+||\nab q||_{\lliw}\Big)\label{ell est II ww}.
\end{align}
\begin{proof}
It suffices to prove \eqref{ell est I ww} and \eqref{ell est II ww} for $r=2$.  By \eqref{prototype 3}, we have
\begin{dmath}
||\nab^2 q||_{\lliw} \leq C(K)\Big(||\Pi\nab^2 q||_{\llbw}||\nab q||_{\llbw}+||\lap q||_{\lliw}\Big)\\
\leq \delta C(K)||\Pi\nab^2 q||_{\llbw}+C(\delta^{-1},K)||\nab q||_{\llbw}+||\lap q||_{\lliw}.\label{8.13}
\end{dmath}
On the other hand, by \eqref{prototype 2}, we have
\begin{dmath}
||\nab^2 q||_{\llbw} \leq C(K)\Big(||\Pi\nab^2 q||_{\llbw}+||\lap\nab q||_{\lliw}+||\nab q||_{\lliw}\Big).\label{8.14}
\end{dmath}
Then \eqref{ell est I ww} follows from \eqref{8.13}-\eqref{8.14} and induction with $\delta=1$. To prove \eqref{ell est II ww}, we have via \eqref{prototype 1} that
\begin{align}
||\nab q||_{\llbw}\leq C(K)\Big(||\nab^2 q||_{\lliw}+||\nab q||_{\lliw}\Big).
\end{align}
\eqref{ell est II ww} then follows from \eqref{8.13} and induction.
\end{proof}

\subsection{The weighted tensor estimate}
\thm \label{tensor estimate ww} Suppose that $|\theta|+|\frac{1}{l_0}|\leq K$, and for $q=0$ on $\p\Omega$, then for $m=0,1$
\begin{multline}
||\Pi\nab^{r}q||_{\llbw}\lesssim_{K} ||\Big(\sum_{s\leq r-2}(\cnab^{s}\theta)\Big)\nab_{N}q||_{\llbw}+\sum_{l=1}^{r-1}||\nab^{r-l}q||_{\llbw}\\
 +(||\theta||_{\linf}+\sum_{0\leq l\leq r-2-m}||\cnab^{l}\theta||_{\llbw})(\sum_{0\leq l\leq r-2+m}||\nab^l q||_{\llbw}), \label{tensor est ww}
\end{multline}
where the second line drops if $0\leq r\leq 4$.
\begin{proof}
The proof follows from the interpolation inequalities on the boundary, e.g., Theorem \ref{bdy interpolation}. We refer \cite{CL} Proposition 4.7 for the detailed proof. 
\end{proof}

In addition, the weighted estimate for the second fundamental form $\theta$ is then a immediate consequence.
\thm \label{theta estimate ww}
 Suppose that $|\theta|+|\frac{1}{l_0}|\leq K$, and the physical sign condition $|\nab_N h|\geq \epsilon>0$ holds, then
\begin{multline}
||\cnab^{r-2}\theta||_{\llbw}\lesssim_{K, \frac{1}{\epsilon}}||\Pi \nab^rh||_{\llbw}+\sum_{s=1}^{r-1}||\nab^{r-s}h||_{\llbw}\\
+(||\theta||_{\linf}+\sum_{s\leq r-3}||\cnab^s\theta||_{\llbw})\sum_{s\leq r-1}||\nab^s h||_{\llbw}, \label{theta est ww}
\end{multline}
where the second line drops for $0\leq r\leq 4$.

\subsection{The weighted energy estimates for Euler equations}
The higher order weighted energies for the compressible Euler equations are
\begin{equation}
E_{w,r} = \sum_{s+k=r} E_{w,s,k}+K_{w,r}+\sum_{j\leq r+1}W_{w,j}^2\label{Er ww},\q r\geq 2,\q E_{w,r}^*=\sum_{r'\leq r} E_{w,r'},
\end{equation}
 where
\begin{dmath}
E_{w,s,k}(t)=\frac{1}{2}\int_{\DD_t}\rho \delta^{ij}Q(\p^sD_t^kv_i,\p^sD_t^{k}v_j)w\dx+\frac{1}{2}\int_{\DD_t}\rho e'(h)Q(\p^sD_t^k h,\p^s D_t^kh)w\dx\\
+\frac{1}{2}\int_{\p\DD_t}\rho Q(\p^sD_t^k h,\p^sD_t^k h)\nu w\,dS,  \label{Esk ww}
\end{dmath}
where $\nu = (-\nab_{N}h)^{-1}$ and
\begin{align}
K_{w,r}(t) &= \int_{\DD_t}\rho|\p^{r-1}\curl v|^2w\dx \label{K_r ww},\\
W_{w,r}(t) &= \frac{1}{2}||\sqrt{e'(h)}D_t^rh||_{L^2_w(\DD_t)}+\frac{1}{2}||\nab D_t^{r-1}h||_{L^2_w(\DD_t)} .\label{Er time ww}
\end{align}

Using Theorem \ref{Christodoulou-Lindblad type elliptic est} -- Theorem \ref{theta estimate ww}, the weighted Sobolev lemmas as well as the interpolation inequalities (e.g., Lemma \ref{interior sobolev}, Lemma \ref{boundary soboolev}, Theorem \ref{int interpolation}, Theorem \ref{bdy interpolation} and Theorem \ref{gag-ni thm}), and the fact that our weight $w$ satisfies $|\p^rw|\leq C_rw$, we are able to repeat the analysis we have done in Section \ref{section 4} and Section \ref{section 5} to obtain the weighted elliptic bounds:
\begin{align}
||v||_{w,r,0}^2+||h||_{w,r}^2 \leq C_r(K,M,c_0, E_{w,r-1}^*)E_{w,r}^* \label{intro int  est ww},\\
||D_th||_{w,r}^2+\lee h\ree_{w,r}^2 \leq C_r(K,M,c_0,\frac{1}{\epsilon},  E_{w,r-1}^*)E_{w,r}^* \label{intro bdy est ww},
\end{align}
where 
\begin{align*}
||v||_{w,r,0} := \sum_{k+s=r, k<r} ||\p^s D_t^k v||_{L^2_w(\DD_t)},\\
||h||_{w,r} := \sum_{k+s=r, k<r}||\p^s D_t^k h||_{L^2_w(\DD_t)} + ||\sqrt{e'(h)}D_t^rh||_{L^2_w(\DD_t)},\\
\lee h\ree_{w,r} :=\sum_{k+s=r} ||\p^s D_t^kh||_{L^2_w(\p\DD_t)}.
\end{align*}
But these  yield the analogous energy estimates for $E_{w,r}$.

\thm
Let $E_{w,r}$ be defined as \eqref{Er ww}, then there are continuous functions $C_r$ such that for each fixed $r\geq 1$, we have
\begin{align}
|\frac{dE_{w,r}(t)}{dt}|\leq C_r(K,\frac{1}{\epsilon},M, c_0, E_{w, r-1}^*)E_{w,r}^*(t),
\end{align}
provided the \eqref{e_kk} and a priori assumptions \eqref{geometry_bound}-\eqref{D_te(h)}.

\begin{appendix}
\section{Appendix}
\noindent {\bf List of notations:}

\begin{itemize}
\item $D_{t}$: the material derivative
\item $\p_i$: partial derivative with respect to Eulerian coordinate $x_i$
\item $\DD_t\in\RR^n$: the domain occupied by fluid particles at time $t$ in Eulerian coordinate
\item $\Omega\in\RR^n$: the domain occupied by fluid particles in Lagrangian coordinate
\item $\p_a = \frac{\p}{\p y_a}$: partial derivative with respect to Lagrangian coordinate $y_a$
\item $\nab_a$: covariant derivative with respect to $y_a$
\item $\Pi S$: projected tensor $S$ on the boundary
\item $\cnab,\cp$: projected derivative on the boundary
\item $N$: the outward unit normal of the boundary
\item $\theta=\cnab N$: the second fundamental form of the boundary
\item $\sigma = tr(\theta)$: the mean curvature
\item $\kk=\kk(x)$: the sound speed
\item $L^{p}_w(\Omega)$: The weighted $L^p$ space
\item $W^{s,p}_w(\Omega)$: The weighted Sobolev space
\end{itemize}
{\bf Mixed norms}
\begin{itemize}
\item $\lee\cdot\ree_r = \sum_{k+s=r}||\nab^sD_t^k\cdot||_{\llb}$
\item $||\cdot||_{r,0} = \sum_{s+k=r,k<r}||\nab^s D_t^k\cdot||_{\lli}$
\item $||\cdot||_{r} = ||\cdot||_{r,0}+||\sqrt{e'(h)}D_t^r\cdot||_{\lli}$
\item $||\cdot||_{r,1,0} = \sum_{k+s=r, k<r-1}||\nab^s D_t^k \cdot||_{\lli} + ||\sqrt{e'(h)}\nab D_t^{r-1}\cdot||_{\lli}$,
\item $||\cdot||_{r,1} = ||\cdot||_{r,1,0} + ||e'(h)D_t^r\cdot||_{\lli}$.

\end{itemize}
{\bf Weighted norms}
\begin{itemize}
\item $||u||_{L^p_w(\Omega)} = (\int_{\Omega}|u(x)|^p w(x)\dx)^{1/p}$
\item $||u||_{W^{s,p}_w(\Omega)} = \sum_{|\alpha|\leq s}||\nab^{\alpha}u||_{L^p_w(\Omega)}$
\end{itemize}
\subsection{The geometry of the boundary, extension of normal to the interior and the geodesic normal coordinate}
\indent The definition of our energy (\ref{Er}) relies on extending the normal to the interior, which is done by foliating the domain close to the boundary into the surface that do not self-intersect.  We also want to control the time evolution of the boundary, which can be measured by the time derivative of the normal in the Lagrangian coordinate. We conclude the above statements by the following two lemmas, whose proof can be found in \cite{CL}.\\
\lem \label{trace 1} let $l_0$ be the injective radius (\ref{inj rad}), and let $d(y)=dist_g(y,\p\Omega)$ be the geodesic distance in the metric $g$ from $y$ to $\p\Omega$. Then the co-normal $n=\nab d$ to the set $S_a=\p\{y\in\Omega:d(y)=a\}$ satisfies, when $d(y)\leq \frac{l_0}{2}$ that
\begin{align}
|\nab n| \lesssim |\theta|_{L^{\infty}(\p\Omega)}, \label{extending nab n}\\
|D_t n| \lesssim |D_t g|_{L^{\infty}(\Omega)},
\end{align}
where we have used the convention that $A\lesssim B$ means $A\leq CB$ for universal constant $C$.
\lem \label{trace 2} let $l_0$ be the injective radius (\ref{inj rad}),and let $d_0$ be a fixed number such that $\frac{l_0}{16}\leq d_0\leq \frac{l_0}{2}$. Let $\eta$ be a smooth cut-off function satisfying $0\leq \eta(d)\leq1$, $\eta(d)=1$ when $d\leq\frac{d_0}{4}$ and $\eta(d)=0$ when $d>\frac{d_0}{2}$. Then the psudo-Riemannian metric $\gamma$ given by
$$
\gamma_{ab}=g_{ab}-\tilde{n}_a\tilde{n}_b,
$$
where $\tilde{n}_c=\eta(\frac{d}{d_0})\nab_cd$
satisfies
\begin{align}
|\nab\gamma|_{L^{\infty}(\Omega)}\lesssim(|\theta|_{L^{\infty}(\p\Omega)}+\frac{1}{l_0})\\
|D_t\gamma(t,y)|\lesssim |D_t g|_{L^{\infty}(\Omega)}. \label{Dtgamma}
\end{align}
\rmk The above two lemmas yield that the quantities $|D_t n|$ and $|D_t\gamma(t,y)|$ involved in the $Q$-inner product is controlled by the a priori assumptions (\ref{geometry_bound})-(\ref{D_te(h)}),since $D_tg$ behaves like $\nab v$ by (\ref{Dtg}). Hence, the time derivative on the coefficients of the $Q$-inner product generates only lower-order terms. In addition, by (\ref{geometry_bound}) , $|\nab n|$ and $|\nab\gamma|$ are controlled by $K$, which is essential when proving the Christodoulou-Lindblad type elliptic estimates.

The next lemma introduces the partition of unity $\{\chi_i\}$ in a domain with sufficient regular boundary.
\lem \label{cut off function} Suppose that $\Omega\in\RR^n$ is a domain whose boundary satisfying the condition $|\theta|+\frac{1}{l_0}\leq K$. Then there are functions $\chi_i\in C_c^{\infty}(\RR^n), i=1,2,\cdots$, such that
\begin{align}
0\leq \chi_i\leq 1,\q \sum \chi_i =1,\q \sum |\p^{\alpha}\chi_i|\leq C_{\alpha}K^{|\alpha|},\q diam(supp(\chi_i))\leq K^{-1},
\end{align}
and for each $x\in\RR^n$, there are at most $16^n$ $i$'s such that $\chi_i(x)\neq 0$. Furthermore, either $supp(\chi_i)\cup \p\Omega$ is empty or is part of a graph contained in $\p\Omega$, for which (possibly after a rotation) is given by
\begin{align}
x_n = f_i(x'),\q |\p f_i|\leq c_1,\q N(x_i)=e_n,\q \text{for}\,\, |x'-x_i'|\leq l_0. 
\end{align}
\begin{proof}
See \cite{CL}.
\end{proof}

\subsection{Sobolev lemmas}
Let us now state some Sobolev lemmas in a domain with boundary, whose proofs are standard and can be found in \cite{CL},\cite{Ev} and \cite{T}.
\lem \label{interior sobolev} (Interior Sobolev inequalities)
Suppose $\frac{1}{l_0}\leq K$ and $\alpha$ is a $(0,r)$ tensor, then
\begin{align}
||\alpha||_{L^{\frac{2n}{n-2s}}(\Omega)}\lesssim_{K} \sum_{l=0}^s||\nab^l\alpha||_{\lli},\q 2s<n,\label{A7}\\
||\alpha||_{L^{\infty}(\Omega)}\lesssim_{K} \sum_{l=0}^{s}||\nab^l\alpha||_{\lli},\q 2s>n.\label{A8}
\end{align}
These inequalities remains valid in weighted spaces $L^p_w(\Omega)$ if the weight satisfies $|\p^r w|\leq C_r w/(1+|x|)^r$.
\begin{proof}
See \cite{CL}. 
\end{proof}
\noi Similarly, on the boundary $\p\Omega$, we have
\lem \label{boundary soboolev} (Boundary Sobolev inequalities)
\begin{align}
||\alpha||_{L^{\frac{2(n-1)}{n-1-2s}}(\Omega)}\lesssim_{K} \sum_{l=0}^s||\nab^l\alpha||_{\llb},\q 2s<n-1,\\
||\alpha||_{L^{\infty}(\Omega)}\lesssim_{K} \delta||\nab^s\alpha||_{\llb}+\delta^{-1}\sum_{l=0}^{s-1}||\nab^l\alpha||_{\llb},\q 2s>n-1,
\end{align}
for any $\delta>0$. These inequalities remain valid in weighted spaces $L^p_w(\Omega)$ as well. In addition, for the boundary we can also interpret the norm be given by the inner product $\langle \alpha, \alpha\rangle=\gamma^{IJ}\alpha_I\alpha_J$, and the covariant derivative is then given by $\cnab$. 

\subsection{Interpolation on spatial derivatives}
We shall first record spatial interpolation inequalities. Most of the results are are standard in $\RR^n$, but we must control how it depends on the geometry of our evolving domain. The coefficients involved in our inequalities depend on $K$, whose reciprocal is the lower bound for the injective radius $l_0$. 
\thm \label{int interpolation} (Interior interpolation)
Let $u$ be a $(0,s)$ tensor, and suppose $\frac{1}{l_0}\leq K$, we have
\begin{align}
\sum_{j=0}^{l}||\nab^ju||_{L^{\frac{2r}{k}}(\Omega)}\lesssim ||u||_{L^{\frac{2(r-l)}{k-l}}(\Omega)}^{1-\frac{l}{r}}(\sum_{i=0}^{r}||\nab^iu||_{\lli}K^{r-i})^{\frac{l}{r}}. \label{interpolation}
\end{align}
In particular, if $k=l$,
\begin{align}
\sum_{j=0}^{k}||\nab^ju||_{L^{\frac{2r}{k}}(\Omega)}\lesssim ||u||_{L^{\infty}(\Omega)}^{1-\frac{k}{r}}(\sum_{i=0}^{r}||\nab^iu||_{\lli}K^{r-i})^{\frac{k}{r}}. \label{31}
\end{align}
These inequalities remains valid when $L^p(\Omega)$ is replaced by $L^p_w(\Omega)$ if $w\geq 0$ satisfies $|\p^r w|\leq C_r w / (1+|x|)^r$.
\begin{proof}
It suffices to prove \eqref{interpolation} with $s=0$, i.e., when $u$ is a function, since $u$ can be replaced by its magnitude $|u|$. Furthermore, since \eqref{interpolation} is equivalent to
\begin{align}
\sum_{j\leq l}||\nab^j u||_{L^s(\Omega)}\leq C(K)||u||_{L^q(\Omega)}^{1-a}(\sum_{i\leq r}||\nab^i u||_{L^p(\Omega)})^a,\label{interpolation CL form}
\end{align}
where $a=l/r$ and $\frac{r}{s}=\frac{l}{p}+\frac{r-l}{q}$. We can further reduce \eqref{interpolation CL form} to the case when $r=2$ and $s=1$, because the general cases follow from the logarithmic convexity. 

Using the partition of unity $\{\chi_i\}$ defined in Lemma \ref{cut off function}, we write $u=\sum u_i$, where $u_i=\chi_iu$. In a neighbourhood of $supp(\chi_i)$, we can then write $\Omega$ as a graph after a rotation:
\begin{align*}
x_n = f(x'),\q |\p f|\leq C.
\end{align*}
We now define the reflection
\begin{align*}
\tilde{u}_i(x)=
\begin{cases}
u_i(x),\q \text{when}\,x\in\Omega\\
u_i(\tilde{x}),\q \text{when}\, x\in \Omega^c
\end{cases}
\end{align*}
Here, $\tilde{x}=(x',x_n-2(x_n-f(x'))$. Then by the interpolation in $\RR^n$, we have
\begin{align*}
||\nab \tilde{u}_i||_{L^s(\RR^n)}^2\leq ||\tilde{u}_i||_{L^q(\RR^n)}||\nab^2\tilde{u}_i||_{L^p(\RR^n)}.
\end{align*}
But since for every $1\leq p'\leq \infty$ and $|{\p\tilde{x}^i}/{\p x^j}|\leq C$, 
\begin{align*}
||\nab^{\alpha}\tilde{u}_i||_{L^{p'}(\RR^n)}\leq C(||\nab^{\alpha}u_i||_{L^{p'}(\Omega)}+||\nab^{\alpha}\tilde{u}_i||_{L^{p'}(\Omega^c)})
\leq C||\nab^{\alpha}u_i||_{L^{p'}(\Omega)},
\end{align*}
for $|\alpha|\leq 2$. Furthermore,  we have
\begin{align*}
||\nab u_i||_{L^{p'}(\Omega)} \leq C||(\nab\chi_i)u||_{L^{p'}(\Omega)}+C||\chi_i \nab u||_{L^{p'}(\Omega)},\\
||\nab^2 u_i||_{L^{p'}(\Omega)} \leq C||(\nab^2\chi_i)u||_{L^{p'}(\Omega)}+C||(\nab\chi_i)\nab u||_{L^{p'}(\Omega)}+C||\chi_i \nab^2 u||_{L^{p'}(\Omega)}
\end{align*}
and this gives \eqref{interpolation CL form} via Lemma \ref{cut off function} for $l=1$ and $r=2$. The general case follows by letting $M_k=\sum_{i\leq k}||\nab^i u||_{L^{s(k)}}$, and so far we have proven $M_1\lesssim M_0M_2$, and hence we get $M_k^2\lesssim M_{k-1}M_{k+1}$ follows from this special case. But the logarithmic convexity then gives $M_k\lesssim M_0^{(r-l)/r}M_r^{l/r}$. Finally, the weighted case follow from the non-weighted case since $|\p^r w|\lesssim |w|/(1+|x|)^r$.
\end{proof}

\subsection{Interpolation on $\p\Omega$}

\thm \label{bdy interpolation} (Boundary interpolation)
Let $u$ be a $(0,s)$ tensor, then
\begin{equation}
||\cnab^l u||_{L^{\frac{2r}{k}}(\p\Omega)}\lesssim ||u||_{L^{\frac{2(r-l)}{k-l}}(\p\Omega)}^{1-\frac{l}{r}}||\cnab^ru||_{\llb}^{\frac{l}{r}}. \label{boundary interpolation ineq}
\end{equation}
In particular, if $k=l$,
\begin{align}
||\cnab^k u||_{L^{\frac{2r}{k}}(\p\Omega)}\lesssim ||u||_{L^{\infty}(\p\Omega)}^{1-\frac{k}{r}}||\cnab^ru||_{\llb}^{\frac{k}{r}}.
\end{align}
Furthermore, if $w\geq 0$ satisfies $|\p^r w|\leq C_rw/(1+|x|)^r$, then
\begin{equation}
||\cnab^l u||_{L^{\frac{2r}{k}_w}(\p\Omega)}\lesssim ||u||_{L^{\frac{2(r-l)}{k-l}_w}(\p\Omega)}^{1-\frac{l}{r}}(\sum_{i\leq r}||\cnab^i u||_{L^2_w(\p\Omega)}^{\frac{l}{r}}).\label{boundary interpolation eq ww}
\end{equation}

\begin{proof}
The proof for \eqref{boundary interpolation ineq} can be found in \cite{CL}, and \eqref{boundary interpolation eq ww} follows from the same proof and the lower order terms on the RHS is generated when the derivatives fall on the weight function $w$.
\end{proof}
\subsection{Elliptic estimates in weighted Sobolev spaces} \label{section A5}
This section is devoted to set up the elliptic estimates in weighted Sobolev spaces $H^{s}_{w}(\Omega)$ (Definition \ref{weighted sob sp}) with weight $w(x) = (1+|x|^2)^\mu$, $\mu\geq 2$, where $\Omega\subset\RR^n, n=2,3$ be a smooth domain, diffeomorphic to the half space $\{x\in\RR^n:x_n\leq 0\}$. Consider the Dirichlet boundary value problem
\begin{align}
\begin{cases}
\lap u = f,\q \text{in}\,\, \Omega\\ \label{ell eq ww}
u=0,\q \text{in}\,\, \p\Omega
\end{cases}
\end{align}
then the following $L^2$ elliptic estimate holds.
\thm (Boccia-Salvato-Transirico \cite{BST})
Fix $s\geq 2$ and $p\in(0,\infty)$, then
\begin{align}
||u||_{W^{s,p}_w(\Omega)} \leq C(||f||_{W^{s-2, p}_w(\Omega)}+||u||_{L^p_w(\Omega)}), \label{ell est ww}
\end{align}
holds for all $u\in W^{s,p}_w(\Omega)$ that solves \eqref{ell eq ww}.

Now we show that the $||u||_{L^p_w(\Omega)}$ on the RHS of \eqref{ell est ww} can in fact be dropped. It is worth to mention here that we have no problem to drop this term if $\Omega$ were bounded, since $\lambda=0$ is not an eigenvalue of $\lap$ in this case (e.g, chapter 6.2 in Evans \cite{Ev}). However, it is in general impossible to drop the term $||u||_{L^2}$ in elliptic estimates when $\Omega$ is unbounded, unless $u$ is sufficiently smooth and decays fast enough at infinity. 

\thm (Rellich-Kondrachov embedding for weighted spaces)\label{RK ww}
The spaces $H^1_{0,w}(\Omega)$ (the space consists of $u\in H^{1}_w(\Omega)$ with $u|_{\p\Omega}=0$) are compactly embedding in the spaces $L^q(\Omega)$ for any $q<2n/(n-2)$.
\begin{proof}
We follow the proof given by Gilbarg-Trudinger \cite{GT} with some modifications. We initially assume $q=1$. Let $\mathcal{A}$ be a bounded subset in $H^{1}_{0,w}(\Omega)$. Without loss of generality we assume that $\mathcal{A}\in C_c^1(\Omega)$ and that $||u||_{H^{1}_w(\Omega)}\leq 1$. For fixed $\delta>0$, let $\mathcal{A}_\delta:=\{u_{\delta}:u\in \mathcal{A}\}$, where $u_\delta$ is the mollification of $u$, i.e., $u_{\delta} = \eta_{\delta}*u$, where $\eta(x)$ is a smooth bump function supported in the unit ball satisfying $\int \eta(x)\dx =1$, and $\eta_{\delta} = \delta^{-n}\eta(\delta^{-1}x)$. 

For each $u\in\mathcal{A}$, we have
\begin{align*}
||u_{\delta}(x)||_{L^{\infty}(\Omega)}\leq \delta^{-n}||\eta||_{L^{\infty}(\Omega)}||u||_{H^{1}_w(\Omega)},\\
||\nab u_{\delta}(x)||_{L^{\infty}(\Omega)}\leq \delta^{-n-1}||\nab \eta||_{L^{\infty}(\Omega)}||u||_{H^{1}_w(\Omega)},
\end{align*}
and so $\mathcal{A}_\delta$ is a bounded, equicontinuous subset of $C_c(\Omega)$ and hence precompact in $C_c(\Omega)$, and consequently also precompact in $L^1(\Omega)$. Nevertheless, we have

\begin{dmath*}
|u(x) - u_{\delta}(x)|\leq \int_{|z|\leq 1}\eta(z)|u(x)-u(x-\delta z)\,dz\\
\leq \int_{|z|\leq 1}\eta(z)\int_0^{\delta |z|}|\nab_r u(x-r\frac{z}{|z|})\,dr\,dz,
\end{dmath*}
and hence
\begin{dmath*}
\int_{\Omega}|u(x) - u_{\delta}(x)|\dx \leq \delta\int_{\Omega}|\nab u|\dx \\
\leq \delta(\int_{\RR^n}\frac{1}{(1+|x|^2)^\mu}\dx)^{1/2}||u||_{H^1_w(\Omega)}.
\end{dmath*}
But since $\int_{\RR^n}\frac{1}{(1+|x|^2)^\mu}\dx<\infty$ when $n\leq 3$ and so $u_\delta$ is uniformly close to $u$ in $L^1(\Omega)$. It then follows that $\mathcal{A}
$ is precompact in $L^1(\Omega)$. Now, for any $q<2n/(n-2)$, we have
\begin{align*}
||u||_{L^q(\Omega)} \lesssim ||u||_{L^1}^a||u||_{L^{2n/(n-2)}}^{1-a}
\end{align*}
for some $0<a<1$ via interpolation. In addition, we have
$$
||u||_{L^{2n/(n-2)}}\lesssim ||u||_{H^1_w},
$$
by Sobolev lemma and the fact that $w(x)\geq 1$. This concludes that a bounded set in $H^1_{0,w}(\Omega)$ must be precompact in $L^q(\Omega)$. 
\end{proof}

\rmk The classical Rellich-Kondrachov embedding theorem yields that $H^1(\Omega)$ is  compactly embedding in the spaces $L^q(\Omega)$ when $\Omega$ is bounded. 

\thm(Improved elliptic estimates) \label{improve ell est thm}
Let $u\in H^s_w(\Omega)\cap H^1_{0,w}(\Omega)$ be a function that solves \eqref{ell eq ww}, and if $f\in H^{s-2}_w(\Omega)$ then
\begin{align}
||u||_{H^{s}_w(\Omega)} \leq C||f||_{H^{s-2}_w(\Omega)}. \label{improve ell est ww}
\end{align}
\begin{proof}
It suffices to prove \eqref{improve ell est ww} when $s=2$. If \eqref{improve ell est ww} is not true, then there exists a sequence $\{u_m\}\subset H^2_w(\Omega)\cap H^1_{0,w}(\Omega)$ satisfying
\begin{dmath*}
||u_m||_{L^2_w(\Omega)}=1,\q  ||u_m||_{L^2(\Omega)}\leq 1,\q ||\lap u_m||_{L^2_w(\Omega)}\to 0.
\end{dmath*}
By virtue of the apriori estimate \eqref{ell est ww}, Theorem \ref{RK ww}, and the weakly compactness of bounded subsets in $H^2_w(\Omega)$, there exists a subsequence, relabelled as $\{u_m\}$, converging weakly to a function $u\in  H^2_w(\Omega)\cap H^1_{0,w}(\Omega) $ satisfying $||u||_{L^2_w(\Omega)}= 1$. However, for any $\phi\in L_w^2(\Omega)$, we must have
\begin{equation*}
\int_{\Omega} \phi (\lap u)w = 0.
\end{equation*}
Hence, $\lap u =0$ and so $u=0$ by the uniqueness assertion (e.g. G-T \cite{GT}, Theorem 8.9 or maximum principle since $u$ decays to $0$ at $\infty$). But this implies $||u||_{L^2_w}=0$, a contradiction. 
\end{proof}

\subsection{Gagliardo-Nirenberg interpolation inequality}
\thm \label{gag-ni thm}
Let $u$ be a $(0,r)$ tensor defined on $\p\Omega\in\RR^2$ and suppose $\frac{1}{l_0}\leq K$, we have
\begin{align}
||u||_{L^4{(\p\Omega})}^2 \lesssim_K ||u||_{\llb}||u||_{H^1(\p\Omega)},\label{gag-ni}
\end{align}
where $H^1(\p\Omega)$ is defined via tangential derivative $\cnab$. Furthermore, \eqref{gag-ni} remains valid in the case of weighted Sobolev spaces.
\begin{proof}
It suffices for us to work in the local coordinate charts $\{U_i\}$ of $\p\Omega$. We consider the corresponding partition of unity $\{\chi_i\}$, where each $\chi_i$ is supported in $U_i$ and vanishing on the boundary of $U_i$. As proved in Lemma \ref{cut off function}, $\chi_i$ can be chosen to satisify 
$$
\sum_i |\cnab \chi_i| \leq C(K).
$$
Now by the result of Constantin and Seregin \cite{Cons}, we have
$$
||u_i||_{L^4(U_i)}^2 \lesssim ||u_i||_{L^2(U_i)}||\cnab u_i||_{L^2(U_i)}, 
$$
where $u_i=\chi_i u$. But since 
$$
||\cnab u_i||_{L^2(U_i)} = ||\cnab (\chi_i u)||_{L^2(U_i)} \leq |\cnab \chi_i|_{L^{\infty}}||u||_{L^2(U_i)}+||\chi_i\cnab u||_{L^2(U_i)} \label{gag local} .
$$
Hence, (\ref{gag-ni}) follows by summing up (\ref{gag local}). This proof remains valid with $L^p$ being replaced by $L^p_w$, where $w$ is defined in Section \ref{section A5}.
\end{proof}

\subsection{The trace theorem}
\thm \label{trace theorem} (The trace theorem)
Let $\alpha$ be a $(0,r)$ tensor, and assume that $|\theta|_{\linf}+\frac{1}{l_0}\leq K$, then
\begin{align}
||\alpha||_{\llb}\lesssim_{K} \sum_{j\leq 1}||\nab^j\alpha||_{\lli}. \label{trace}
\end{align}
Furthermore, 
\begin{align}
||\alpha||_{L^2_w(\p\Omega)}\lesssim_{K} \sum_{j\leq 1}||\nab^j\alpha||_{L^2_w(\Omega)}. \label{trace ww}
\end{align}
Here, $w$ is defined in Section \ref{section A5}.
\begin{proof}
It suffices to show \eqref{trace ww} only, since the proof for \eqref{trace} is almost identical. Let $N$ be the extension of the normal in the interior of $\Omega$ given by the geodesic normal coordinate (i.e., Lemma \ref{trace 1}). Then
\begin{align}
\int_{\p\Omega}|\alpha|^2w\,d\mu_{\gamma} = \int_{\Omega}\nab_k(N^k|\alpha|^2w)\,d\mu_{g}.
\end{align}
But since $|\nab N|\leq K$ and $|\nab w|\leq Cw$, \eqref{trace ww} follows.
\end{proof}

\paragraph*{Acknowegement}
I would like to express my deepest thanks to my advisor Hans Lindblad for many useful suggestions and comments. I would like to thank Marcelo Disconzi, Theo Drivas, Dan Ginsberg, Chris Kauffman, Yannick Sire, Qingtang Su, Shengwen Wang, Yi Wang, Yakun Xi and Hang Xu for many long and insightful discussions. In addition, I thank the anonymous referee for careful reading and helpful comments.

\end{appendix}

\end{document}